\ifcase2
\or
\documentclass[10pt,a5paper]{article}
    \let\Title\title
    \let\Author\author
    \def\Abstract#1{\def\myAbstract{#1}}
    \def\atPreamble{\usepackage{mybiblatex}%
        \let\inPlot\input
    	}
\or
\documentclass[10pt,a4paper]{article}
    \AtEndDocument{\bibliography{bibexport}}
    \let\Title\title
    \let\Author\author
    \def\Abstract#1{\def\myAbstract{#1}}
    \def\atPreamble{%
        \let\inPlot\includegraphics
        \usepackage{natbib}
        \bibliographystyle{agsm}
        \let\harvardurl\url
        \usepackage{twoopt}
        \newcommandtwoopt{\footcite}[3][][]%
            {\footnote{\cite[##1][##2]{##3}}}
        }
\or
\documentclass[namedate,webpdf,imatrm]{ima-authoring-template}
    \title[Accurate multi-continuum micromorphic
        homogenisations] {Accurate families of multi-continuum
        micromorphic homogenisations in multi-D space-time via
        dynamical systems theory}
    \author{A.~J. Roberts*
        \address{\orgdiv{School of Mathematical Sciences}, 
        \orgname{University of Adelaide}, 
        \orgaddress{
        \state{South Australia}, \country{Australia}}}}
    \def\JournalCase{}
    \AtEndDocument{\bibliography{bibexport}}
    \def\Title#1{}
    \def\Author#1{}
    \let\Abstract\abstract
    \def\atPreamble{%
        \let\inPlot\includegraphics
        \usepackage{natbib}
        \bibliographystyle{agsm}
        \let\harvardurl\url
        \usepackage{twoopt}
        \newcommandtwoopt{\footcite}[3][][]%
            {\footnote{\cite[##1][##2]{##3}}}
        }
    \AtBeginDocument{
        \IfFileExists{tmacInfo.sty}{\usepackage{tmacInfo}}{}
        \let\TMACpara\paragraph  \def\paragraph#1{\TMACpara{\textbf{#1}}\quad}
    }
\fi

\newif\ifJ 
\ifdefined\JournalCase\Jtrue\else\Jfalse\fi

\Title{Construct accurate multi-continuum micromorphic
homogenisations in multi-D space-time with computer
algebra}
\Author{A.~J. Roberts
\thanks{Mathematical Sciences, University of Adelaide, South Australia. 
\protect\url{http://orcid.org/0000-0001-8930-1552},
\protect\url{mailto:profajroberts@protonmail.com}}
}
\date{\today}

\ifcase0     
   \let\include\input	
\or
\or
\or
\or
\or
\or
\else
\fi

\ifJ\usepackage{times,mathptmx}
    \usepackage{url,doi,microtype,graphicx,hyperref}
    \def\notJbreak{}
\else
    \usepackage{ajr}
    \usepackage{pgfplots}\pgfplotsset{compat=newest}
    \usepgflibrary{arrows.meta}
    \usepgfplotslibrary{external}
    \def\notJbreak{\nonumber\\&\quad}
\fi
\usepackage[leftcaption,raggedright]{sidecap}
\usepackage{amsmath,amssymb}
\usepackage[defaultlines=3,all]{nowidow}
\usepackage{defns}    \RaisedNamesfalse
\usepackage{reducecode,mycleveref}

\atPreamble

\allowdisplaybreaks
\numberwithin{equation}{section}

\renewcommand{\vec}[1]{\text{\boldmath$#1$}}
\newcommand{\alphaD}{\ensuremath{\partial_t^\alpha}}
\newcommand{\alphahat}{\ensuremath{\hat\alpha}}
\def\rve{\textsc{rve}}
\def\im{\textsc{im}}
\def\Span{\operatorname{span}}

\newcommand{\ov}{\vec 0}
\Frak K\Frak L\Frak u\Frak d 
\newcommand{\vf}{\ensuremath{\mathfrak v}}
\Vec j\Vec k\Vec f\Vec q\Vec x 
\def\zv{\ensuremath{z}}
\Vec F\Vec G\Vec X\Vec U
\Vec{phi}\Vec{ell}\Vec{theta}
\Cal A\Cal D\Cal E\Cal J\Cal K\Cal L\Cal M
\Bb E\Bb H\Bb J\Bb M\Bb R\Bb T\Bb U\Bb X\Bb Z
\newcommand{\tU}{\tilde U}
\newcommand{\tv}{\tilde v}
\newcommand{\tG}{\tilde G}
\newcommand{\tGv}{\tilde\Gv}
\def\sumM{\sum_{m=0}^{M-1}}
\def\dag{\ensuremath{\dagger}}

\newcommand\dom{{\ensuremath{\mathbb D}}}
\def\x{{\ensuremath{\RaisedName{\detokenize{\x}}{\mathcal{%
    \mathchoice{\scriptstyle X}{\scriptstyle X}%
    {\scriptscriptstyle X}{\scriptscriptstyle X}}}}}}
\def\fx{{\ensuremath{\RaisedName{\detokenize{\fx}}\vec{\mathcal{%
    \mathchoice{\scriptstyle X}{\scriptstyle X}%
    {\scriptscriptstyle X}{\scriptscriptstyle X}}}}}}
\renewcommand{\divv}{\grad^T}
\newcommand{\gradx}{\vec\nabla_{\!\!\fx}}
\newcommand{\divx}{\gradx^T}
\newcommand{\gradq}{\vec\nabla_{\!\!\thetav}}
\newcommand{\divq}{\gradq^T}
\newcommand{\deltab}{\fd}
\renewcommand{\E}[1]{E^{\ifx#1+\else#1\fi\frac12}}
\def\Res{\ensuremath{\operatorname{Res}}}
\let\ellm\cE
\newcommand{\erfc}{\operatorname{erfc}}
\def\fde{\textsc{fde}}
\newcommand{\GM}[1]{(\textsc{gm}, #1)}
\newcommand{\Lap}[1]{\ensuremath{\operatorname{L}\left\{#1\right\}}}

\newcounter{i}

\def\symBox#1{#1}
\def\oSym{\symBox{$\color{green!70!black}\circledcirc$}}
\def\xSym{\symBox{$\color{green!70!black}\otimes$}}
\def\uSym{\symBox{$\color{blue}\blacktriangleright$}}
\def\vSym{\symBox{$\color{red}\blacktriangle$}}

\ifJ
    \newcommand{\ajr}[1]{}
\else
    \usepackage[colorinlistoftodos,textsize=footnotesize]{todonotes}
      \setuptodonotes{inline,color=orange!20}
      \ifx\undefined
        \newcommand{\Xtodo}[2][]{\todo[#1]{#2}} 
      \else
        \newcommand{\Xtodo}[2][]{\tikzexternaldisable%
          \todo[#1]{#2}%
          \tikzexternalenable} 
      \fi 
    \newcommand{\ajr}[1]{\Xtodo[author=AJR,color=yellow!30]{#1}}
\fi

\Abstract{Homogenisation empowers the efficient macroscale system level prediction of physical scenarios with intricate microscale structures.
Here we develop an innovative powerful, rigorous and flexible framework for asymptotic homogenisation of dynamics at the \emph{finite} scale separation of real physics, with proven results underpinned by modern dynamical systems theory.
The novel systematic approach removes most of the usual assumptions, whether implicit or explicit, of other methodologies.
By no longer assuming averages the methodology constructs so-called multi-continuum or micromorphic homogenisations systematically informed by the microscale physics.
The developed framework and approach enables a user to straightforwardly choose and create such homogenisations with clear physical and theoretical support, and of highly controllable accuracy and fidelity.   
}

\begin{document}
\maketitle

\ifJ
\YouMustUseMyOverideOfMathptmx
\else
\begin{abstract}\myAbstract\end{abstract}
\tableofcontents
\nocite{Roberts2024a}%
\fi


\section{Introduction}

\begin{figure}
\centering
\caption{\label{FegDeform} 4\% compression of a material with circular inclusions, seen near either end, may develop a nontrivial microscale structure, seen near the middle \protect\cite[from Figure~4(2nd) of][]{TGuo2024}.  
Macroscale homogenisation of such microscale structures may best be via a multi-continuum model.}
\includegraphics[width=10cm]{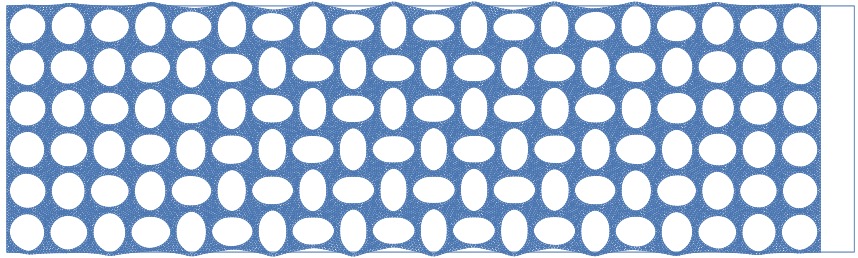}
\end{figure}

In most standard constitutive models for the mechanical behaviour of materials, the stress at a given point uniquely depends on the current deformation gradients \cite[e.g.,][]{Bazant2002}.
However, modern designed metamaterials have exceptional mechanical properties that depend largely on the intricate underlying complicated microstructure, rather than the bulk properties of their constituent materials \cite[e.g.,][]{Sarhil2024}. 
So-called multi-continuum or micromorphic homogenisations, or generalized Cosserat theories, are enhanced models which aim to capture the significant macroscale effects of such heterogenous microstructure without fully-resolving the microstructures.\footcite[e.g.,][]{Alavi2023, Liupekevicius2024, Lakes1995}
\cref{FegDeform} shows one example of a microstructured material whose macroscale deformation can only adequately be predicted by the generalised modelling of non-standard microscale modes \cite[]{Rokos2019}. 
\cite{Altenbach2011} gathered a wide cross-section of applications and approaches to such generalised continuum mechanics.
The pioneering vision of \cite{Muncaster83b} was that the dynamical systems concept of invariant manifolds would provide a general rigorous route to coarse-scale modelling of nonlinear dynamics in mechanics.
Cognate to this vision, we develop a flexible, systematic and practical approach to the multi-continuum homogenisation modelling of heterogeneous systems (such as the material deformation of \cref{FegDeform}) using modern developments in the rigorous theory and application of invariant manifolds\footcite[e.g.,][]{Aulbach2000, Potzsche2006, Haragus2011, Roberts2013a, Bunder2018a} \text{(often abbreviated to~\im{}s).}

\begin{definition}[multi-continuum models] \label{Dmulticont}
Consider the scenario of a physical field of interest~\(u\), evolving in time~\(t\), on a macroscale spatial domain of size~\(L\), but with significant microscale structure on the much smaller size~\(\ell\). 
\begin{itemize}
\item An \emph{\(M\)-mode model} is written in terms of \(M\)~variables, say~\(\Uv:=(U_0,\ldots,U_{M-1})\), each defined over the microscale length~\(\ell\), and evolving according to a closed system of \ode{}s\slash\pde{}s \(\Dn t\alpha{U_i} = G_i(\Uv)\) that are valid for describing the evolution of macroscale structures (that is, significant to the system size~\(L\)) of the physical problem.

\item A \emph{multi-continuum}, or \emph{micromorphic}, model is one where some of the \(M\)~parameters \emph{do not} reflect conserved characteristics of the field~\(u\). 

\item A modelling scheme\slash method is \emph{transitive} if for every system~(\dag) in some class, and for every suitable \(M_1<M_2\), an \(M_1\)-mode model of an \(M_2\)-mode model of~(\dag) is the same as the direct \(M_1\)-mode model of~(\dag).
\end{itemize}
\end{definition}

\begin{example}[tri-continuum homogenisation] \label{EG1D}
\cref{SoneDintro} details homogenisation of a material with \(\ell\)-periodic microstructure in 1-D space. 
Let the physical `displacement', `material' or `heat'~\(u(t,x)\) evolve according to the \text{heterogeneous \pde}
\begin{equation}
\Dn t\alpha u=\D x{}\left\{\kappa(x)\D xu\right\},
\quad 0<x<L\,,
\label{Ehdifpde}
\end{equation}
where \(\alpha=1\) is a diffusion problem, and \(\alpha=2\) is a wave problem.
We focus on modelling interesting dynamics in the interior of the domain, say a connected open set \(\XX\subset[0,L]\), between the physical boundary layers near \(x=0,L\).

Throughout \cref{SoneDintro} we revisit this example, non-dimensionalised in space and time so that microscale period \(\ell=2\pi\) and \(\kappa\approx1\)\,.
We choose the specific heterogeneity coefficient 
\begin{equation}
\kappa(x):=1/(1+a\cos x).
\label{EGkappa}
\end{equation}
over non-dimensional parameters~\(a\), \(|a|<1\)\,.\footnote{%
The computer algebra of \ifJ Supplementary Code \protect\cite[Appendix~A,][]{Roberts2024a} \else\cref{cas} \fi works for this case and for a wide variety of \(2\pi\)-periodic heterogeneity provided that \(\kappa(x)=1\) when \(a=0\)\,.}
Every good homogenisation method would then construct the \emph{one-mode} macroscale homogenised~\pde
\begin{equation}
\alphaD U_0 = U_{0xx}+\cdots\,,\quad\text{for } x\in\XX\,,
\label{Eedifpde}
\end{equation}
for some effective mean field~\(U_0(t,x)\), and with the macroscale effective material constant to be one in this case.
Some homogenisation methods additionally construct higher order corrections denoted by the ellipsis~\(\cdots\) to improve accuracy of the homogenisation (e.g., \cref{EGmfg}).
The model~\cref{Eedifpde} is not termed a multi-continuum model because the macroscale variable~\(U_0\), usually defined as an average of the field~\(u\), is closely connected to the conservation property of the field~\(u\).

Alternatively, \cref{SoneDintro} introduces a novel unified framework for developing a three-mode, tri-continuum, homogenisation expressed in terms of three macroscale quantities~\(U_0,U_1,U_2\), defined by the microscale physics, that evolve within~\XX\ according to the three coupled \text{macroscale \pde{}s}
\begin{align}&
\alphaD U_{0}=(1+\tfrac12a^2)U_{0xx}-\tfrac12aU_{1x}+\cdots,
&&
\alphaD U_{1}=-U_1+aU_{0x}-2U_{2x}+\cdots,
\nonumber\\&
\alphaD U_{2}=-U_2+2U_{1x}+\cdots
&&\text{for }x\in\XX.
\label{Epde3}
\end{align}
This is a multi-continuum model (\cref{Dmulticont}) because although~\(U_0\) is connected to the conservation of~\(u\) through its definition as an unweighted average, the variables \(U_1,U_2\) are not connected to conserved quantities.
Instead \(U_1,U_2\) are connected to specific, physics-informed, microscale structures (such as~\cref{EG1Dv1}) within each microscale \text{period of~\(\kappa(x)\).}

Our multi-continuum modelling is transitive  (\cref{Dmulticont}) since the adiabatic quasi-equilibrium approximation of the second mode in~\cref{Epde3} recovers the classic homogenised one-mode \pde~\cref{Eedifpde}. 
This adiabatic approximation gives \(U_1\approx aU_{0x}\) (upon neglecting~\(U_{2x}\) as being of higher-order).
Consequently, the first \pde\ of~\cref{Epde3} reduces to the usual homogenised macroscale \pde\ \(\alphaD U_{0}\approx(1+\tfrac12a^2-\tfrac12a^2)U_{0xx}\)\,, that is, \text{\(\alphaD U_{0}\approx U_{0xx}\) as in~\cref{Eedifpde}.}
\qed
\end{example}

The following seven subsubsections overview the structure of the article, and highlight some of this article's main contributions to multi-continuum homogenisation in comparison to much extant theory and practice.
Throughout, the article compares and contrasts with other approaches: details of most points of advantage for this novel dynamical systems approach may be found by searching for ``in contrast''---there are~31 such explicitly identified advantageous contrasts.

\subsubsection{Overview of the article's structure}

This article greatly extends, in two major parts, a novel dynamical systems method \cite[][]{Bunder2018a} to a class of multi-continuum\slash micromorphic homogenisations.
The first part, \cref{SoneDintro,Shceq}, introduce the key ideas and methods in the more accessible and commonly addressed scenario of the linear self-adjoint homogenisation of a scalar field in space-time with microscale periodic heterogeneity in one spatial dimension.
The second part, \cref{Sgentheory,Selastic2d}, greatly extends the framework and theory of the methodology to a wide class of nonlinear non-autonomous dynamics of a general field in space-time with arbitrary number of `large spatial' dimensions and heterogeneity on both the macroscale and a periodic or quasi-periodic microscale.
For an example in 2-D elasticity: \cref{SSShoh} derives a homogenisation in terms of \emph{high-order gradients} of the usual two mean displacements; whereas \cref{SSStchm} derives a \emph{tri-continuum} homogenisation in three degrees of freedom as determined by the sub-cell physics of \text{the heterogeneous elasticity.}

Three broad families of methods enhance standard homogenisation.%
\footcite[e.g.,][]{Bazant2002, Forest2011, Sarhil2024}
One may:  firstly, introduce extra fields, the multi-continua, that provide supplementary information on the small-scale kinematics;
secondly, improve the resolution of the standard displacement-field modelling by incorporating higher-order gradients;
and thirdly, introduce nonlocal effects into the homogenisation.
Our systematic dynamical systems framework and results unifies the approach and connects to all three of these families.
In order: 
firstly, multi-continuum Cosserat-like theories \emph{are} the central theme;  
secondly, higher-order gradient models are invoked for all cases of both basic homogenisation (\cref{Shcegsmm,SSShoh}) and also multi-continua homogenisation (\cref{Shcegbmm,secsrh}); and lastly,
examples of accurate nonlocal models are created by regularising some of the higher-order models (\cref{Shcegsmm,SSScsw}).
\cite{Craster2015} commented that ``Homogenization theory is typically limited to static or quasi-static low frequency situations and the purpose of this contribution is to briefly review how to tackle dynamic situations''. 
In contrast, this contribution remedies this limitation by focussing upon dynamics---statics thereby being the special case of \text{no time variations.}

\subsubsection{Introduce the novel methodology}
\label{Sinm}

In a plenary lecture \cite{Forest2011} discussed that ``in most cases \ldots\ proposed extended homogenization procedures [for micromorphic continua] remain heuristic''.  
In contrast, the novel approach here is created by innovatively synthesising three separate ingredients from the mathematics of dynamical systems.
Let's summarise the three ingredients.
\begin{itemize}
\item Firstly, it seems a paradox that macroscale homogenised models have translationally invariant symmetry in space, say~\xv, given that the underling microscale is spatially heterogeneous in~\xv.
Here we resolve the paradox by embedding any given heterogeneous problem in the ensemble of all phase-shifted versions (\cref{secpse,Spse}).
\cite{Smyshlyaev2000} [p.1339] also commented that such an ensemble provides a useful view of homogenisation.
Then the ensemble of problems is translationally symmetric in~\xv, with the heterogeneity in an orthogonal `thin space' dimension.
Consequently, homogenisation then preserves the embedded system's translational \text{symmetry in~\xv.}

Interestingly, \cite{Auriault2009} assert that ``homogenization of a medium with a high density of heterogenities is \emph{only possible} if we consider regions containing a large number of these heterogenities'' (my italics).
In contrast, by using an ensemble of all phase-shifts we \emph{prove} (\cref{Pgenft,Pgenabt,Pgensm}) that homogenised models exist at a finite scale separation that may only contain as little as, say, three or so periods \text{of the heterogeneity.}

This embedding also encompasses extensions (\cref{Spse}) to arbitrary quasi-periodic microscale heterogeneity \cite[]{Roberts2022a}.

\item Secondly, the orthogonal `thin space' dimension gives rise to the classic cell-problem (\rve\ problem) with periodicity naturally required (\cref{secsmh,Simmmh}), instead of being assumed.
Then \im~theory supports choosing multi-continuum modes from the eigenmodes corresponding to small eigenvalues because for a wide range of circumstances the theory then guarantees that the multi-continuum modelling is emergent (\cref{Dim}).
Such multi-mode modelling in spatially extensive systems was initiated by \cite{Watt94b} for the shear dispersion in flow along a channel.

\item Thirdly, recently \cite{Bunder2018a} developed an \im\ model at each spatial locale, weakly coupled to neighbouring locales, via a multivariate Taylor series in space.  
A remarkable new remainder term  for the series \cite[expression~(44)]{Bunder2018a} quantifies the error of the model so that the homogenised model is valid simply wherever and whenever the remainder term is small enough for the \text{purposes at hand.}

\end{itemize}

Recently, \cite{Alavi2023} [p.2166] noted one limitation of other approaches is that ``proper elaboration of the macroscopic kinematic and static quantities that pertain to the micromorphic continuum is a problematic issue \ldots\ [previous] works did not relate macroscopic to microscopic kinematic quantities and they did not formulate a boundary value problem at the microscale.''
In contrast,  embedding the system provides, via dynamical systems theory, a sound boundary value problem (e.g., \cref{Eemdifpde,Eempde}), and the constructed \im{}s explicitly connect the microscale to the selected macroscale kinematic \text{quantities (e.g., \cref{EG1Dv1,EhetnonIM,EhcegManifold2u}).}

Further, Maugin asserted \cite[p.14]{Altenbach2011} ``any relationship that can be established with a sub-level degree of physical description is an asset that no true physicist can discard.''
Herein we develop a new clear relationship between macro- and microscale dynamics that thus should be acquired, especially for nonlinear out-of-equilibrium problems.

\subsubsection{Strong theoretical support}

For classic homogenisation, much research has established rigorous convergence and error bounds \cite[e.g.,][]{Dohnal2015, Zhikov2016}.
However, in our context of multi-continua\slash micromorphic homogenisation, and akin to the plenary lecture by \cite{Forest2011}, \cite{Alavi2023} [p.2165] commented that ``Periodic homogenization methods \ldots\ are most of the time lacking a clear mathematical proof of convergence of the microscopic solution towards the limit solution for vanishing values of the scale parameter''.
Complementing this comment, in their review \cite{Fish2021} discuss that [p.775] ``Mathematical homogenization theory based on the multiple-scale asymptotic expansion assumes scale separation.'' 
Such scale separation is usually defined as the mathematical limit \(\ell/L\to0\) for microscale lengths~\(\ell\) and macroscale lengths~\(L\).  
In contrast, here there is \emph{no} such limit assumed or imposed: there is \emph{no} imposition of parameters having to be scaled; there is \emph{no}~\(\epsilon\) used in this work.
Consequently, herein a reader has to put aside much of the previous traditional of \text{classic homogenisation.}
  
Instead, the established dynamical systems theory of \im{}s provides requisite theorems in the two main cases\footcite[e.g.,][]{Aulbach2000, Potzsche2006, Haragus2011, Roberts2013a, Bunder2018a}.  
A first class of systems are those with some dissipation for which \cref{SSSsd} creates a framework with associated mathematical proofs of the existence of homogenised models with exact closure, with convergence in time to the macroscale homogenisation, and that is controllably approximated.
A second class is that of undamped, wave-like, systems where \cref{SSSwls} argues that we can prove that constructible nearby systems possess the exact guiding-centre macroscale homogenised models constructed. 
Moreover, this theory and its proofs are not just for the vanishing scale ratio limit, but for the finite scale separations \text{of real applications.}

\cite{Alavi2023} [pp.2164--5] discussed how ``generalized continuum models \ldots\ face some limitations in their capability to predict microstructures' supposedly intrinsic mechanical properties accurately.''
Similarly,  \cite{Fish2021} [p.776] commented that ``In upscaling methods, the fine-scale response is approximated or idealized, and only its average effect is captured.'' 
In contrast,  \im\ theory guarantees the existence of an \emph{exact} closure for its multi-continuum homogenisations in some domain (\cref{Simmmh}), a closure that in principle exactly captures the fine-scale response---this in-principle exactness follows from the upcoming \cref{Dim} that an \im\ is composed of exact solutions of the given physical system.
Importantly, our approach does \emph{not} invoke assumed averaging.  
Further,  theory \cite[]{Potzsche2006} guarantees we approximate the exact closure to \text{a controllable accuracy (see \cref{CordAcc}).}

\subsubsection{Quantifying errors and uncertainty}

In a broad-ranging discussion of multiscale simulation of materials, \cite{Chernatynskiy2013} [p.160] asserted that ``quantifying the errors at each scale is challenging, quantifying how those errors propagate across scales is an even more daunting task'' which has two aspects rectified here.
Firstly, as \cref{Sinm} mentioned, our approach connects with innovative quantitive estimates for the remainder error incurred by the macroscale modelling: in the case of one large spatial dimension via expression~(23) of \cite{Roberts2013a}; in the case of multiple large dimensions via expression~(44) of \cite{Bunder2018a}.  
In contrast to the usual error estimates for homogenisation, such as error\({}< C_0\epsilon\) for some~\(C_0\) \cite[e.g.,][Thms.~2.2--2.4]{Dohnal2015}, these expressions~(23) and~(44), and appropriate generalisations, give the exact remainder error, in every spatial locale, for any chosen approximate \text{\im~homogenised model.} 

Secondly, the dynamical systems framework established here comes with the methodology to construct (in future research) an associated systematic accurate projection of initial conditions \cite[e.g.,][]{Roberts89b, Watt94b}.  

The same projection also rationally propagates errors and uncertainty from micro- to macro-scales.
Thus error propagation no longer need be such a ``daunting task''.
Such projection should then also be informative in optimising microscale topology towards desired macroscale behaviour \cite[e.g.,][]{Herrmann2024, Shahbaziana2022}.

\subsubsection{Clarify and guide multi-continua choices}

\cite{Rizzi2021} [p.2253] concluded that ``It is thus an an important future task \ldots\ to derive guidelines for a favorable choice among the numerous available micromorphic continuum models.''
The dynamical systems approach established here clarifies and guides this choice:
\cref{secsmh,secsee,Simmmh,Se2cmmh} organise choices into a set of possibilities, and then guides a choice in the set depending upon the spatio-temporal scales of interest in the envisaged \text{scenarios of application.}

\cite{Alavi2023} [p.2164] also raises the ``difficulty to choose a priori an appropriate model for a given microstructure.
Enriched continuum theories are required in such situations to capture the effect of spatially rapid fluctuations at the mesoscopic and macroscopic levels. \ldots\ micromorphic homogenization raises specific difficulties compared to higher gradient and higher-order theories''  
The dynamical systems framework herein, as discussed by \cref{SSSisr}, empowers us to also guide the selection of multi-continuum modes aimed to improve spatial resolution, irrespective of \text{the temporal resolution.}

\cite{Efendiev2023} also use a spectral decomposition to guide multi-continuum models.  
However, they require macroscale variables~\(U_i\) to be \emph{averages} and hence deduce [p.3] ``unless there is some type of high contrast, the averages \ldots\ will become similar'' which leads them to conclude [p.16] that for ``multi-continuum models, \ldots\ one needs high contrast to have different average values''.  
In contrast, the physics-informed approach here is \emph{not} wedded to averaging, so is more flexible, and need not be limited to high contrast (although \cref{Shceq} details a high contrast example, the examples of \cref{SoneDintro,Selastic2d} are not).  
By measuring and being based upon whatever structures the physical equations tell us about sub-cell dynamics, we form and justify a wide family \text{of multi-continuum models.}

\cite{Eggersmann2019} discuss that observed history dependent elasticity is often equivalent to having ``internal variables''---variables  that we call micromorphic, like~\(U_1,U_2\) in the model~\cref{Epde3} of \cref{EG1D}.  
\cite{Eggersmann2019} [p.9] comment that ``the essential conceptual drawback of the internal variable formalism is that the internal variable set is often not known or is the result of modeling assumptions.''
In contrast, here we determine the internal variable set in a rigorous framework from the sub-cell dynamics, and systematically choose the specific internal variables from the desired space-time resolution \text{for the model (\cref{Simmmh}).}

In contrast to most other homogenisation approaches \cite[e.g., the review by ][]{Fronk2023}, many commonly made assumptions are not needed here.
There are no guessed fast/slow variables, no small~\(\epsilon\)s, no limits, no need to assume a variational formulation, nor energy functional, nor presume specific sub-cell modes, nor need oversampling regions, nor impose assumed boundary conditions for \rve{}s, nor assume multiple times scaled by powers of a small~\(\epsilon\).
Consequently, I contend the \im\ methodology developed herein best fits the criteria of \cite{Auriault2009} [p.59] that ``the ideal [homogenisation] process should be independent of any assumption on the physics of the model on \text{the macroscopic scale''.}

\subsubsection{Functionally graded materials}

The general theory of \cref{Sgentheory} encompasses the homogenisation of fine-scale structures in materials which also have macroscale variations in their properties (e.g., see \cref{EGmfg}).  
Two practical examples are axially functional graded beams \cite[e.g.,][]{Gantayat2022} and piezoelectric composite materials in curved shapes \cite[]{Guinovart2024}.
Generalised multi-continuum homogenisations of such curved materials would combine the theoretical framework of \cref{Sgentheory} with practical algebraic techniques developed previously for general curvilinear coordinates \cite[e.g.,][]{Roberts99b}.

\subsubsection{Computer algebra handles complicated details}

\cite{Rizzi2021} [p.2237] began by highlighting that a ``basic problem of all these theories \ldots\ is the huge number of newly appearing constitutive coefficients which need to be determined.''  
The systematic nature of our approach, together with theoretical support \cite[e.g.,][]{Potzsche2006}, empowers computer algebra to routinely handle the potentially ``huge number'' of constitutive coefficients via a robust and flexible algorithm introduced by \cite{Roberts96a} and documented for many applications in a subsequent book \cite[]{Roberts2014a}.  
\ifJ Supplementary Material \cite[Appendices~A--C]{Roberts2024a} \else\cref{cas,cashc,cas2d} \fi provide more details and list adaptable code for the three major examples of \cref{SoneDintro,Shceq,Selastic2d}---the latter two examples integrate fine-scale numerics with the macroscale algebra.
These codes use the computer algebra system Reduce\footnote{%
\url{http://www.reduce-algebra.com}} as it is free,
flexible, and fast~\cite[]{Fateman2002}.

\cite{Alavi2023} [p.2164] commented that ``micromorphic homogenization raises specific difficulties compared to higher gradient and higher-order theories'' 
In contrast, here we apply the same theoretical and practical framework for all cases with no difficulty, indeed with the same code via just a couple of parameter changes to change from standard to multi-continuum micromorphic homogenisations, and from leading order to any arbitrary higher-order in gradients \text{and/or nonlinearity.}

Let's proceed with the basics of the approach to homogenisation in 1-D, \cref{SoneDintro,Shceq}, before addressing the complexities of general nonlinear systems, potentially non-autonomous and potentially with quasi-periodic microscale, in multi-D space, \cref{Sgentheory,Selastic2d}.


\section{Systematic multi-continuum homogenisation for 1-D spatial systems}
\label{SoneDintro}

This section introduces the proven construction of multi-continuum or micromorphic homogenised models in 1-D space.  
Now a ``multi-continuum'' or ``micromorphic'' model (\cref{Dmulticont}) is phrased in terms of microscale structures\slash modes that have macroscale variations---in essence pursuing models with multiple mode-shapes in the microscale \cite[e.g.,][]{Rokos2019, Alavi2023, Sarhil2024}.
We develop multi-continuum micromorphic homogenisation by the novel combination of analysing the ensemble of phase-shifts \cite[]{Roberts2013a, Smyshlyaev2000} via innovations to the multi-modal dynamical systems modelling of \cite{Watt94b}.
Using ensembles to frame homogenisation dates back to at least \cite{Willis1985}, but a key difference here is that we do \emph{not assume} ensemble averages.

The distinguishing feature here, in contrast to most other multi-continuum approaches, is that we let the microscale physics of the system inform and determine almost all decisions.  
We do not assume any particular weighted averages are appropriate for variables \cite[e.g., as done by][p.870]{Milton2007}, nor do we assume any particular weighted averages are appropriate for any upscaling dimensional reduction.  
The \emph{main} subjective decision made herein are how many microscale modes one wishes to resolve in the multiscale modelling---in this approach we may use any number of modes.  
This decision is strongly guided by the physics-determined spectrum of the microscale (\cref{secsee,Simmmh}), and by the timescales a user needs to resolve in the model \text{(\cref{Simmmh} gives details).}

\cite{Fish2021} [p.774] identified that a challenge for ``a multiscale approach involves a trade-off between increased model fidelity with the added complexity, and corresponding reduction in precision and increase in uncertainty''. 
Although there is a trade-off with complexity, in contrast, here there is \emph{no} trade-off with increased fidelity: both here and in \cref{Sgentheory} we establish a framework with proven controllable \text{precision and certainty.}

The theoretical support for multi-continuum models such as~\cref{Epde3} depends upon the time evolution operator~\alphaD.
For the attributes of existence and construction, for general~\alphaD\ it appears best to appeal to a version of \emph{backwards theory} \cite[]{Roberts2018a, Hochs2019}, namely that the constructed invariant manifold~(\im), and homogenised evolution~\cref{Epde3} thereon, is \emph{exact} for a system \emph{close} to the specified multiscale \pde~\cref{Ehdifpde}.
The major practical difference among the various~\alphaD\ is whether the invariant manifold, multi-modal, homogenisations such as~\cref{Epde3} are \emph{emergent} in time (\cref{Dim}):  
for diffusion, \(\alpha=1\), the homogenised models are generally emergent;  
for elastic waves, \(\alpha=2\), the homogenised models are best viewed as a \emph{guiding centre} for the dynamics about the constructed manifold \cite[e.g.,][]{vanKampen85};
whereas for other cases the relevance of the homogenisation depends upon how the spectrum of the right-hand side operator of~\cref{Ehdifpde} maps into dynamics of the corresponding modes.
\cref{Simmmh} discusses these cases in detail. 

Of course, if one only addresses forced equilibrium problems, or Helmholtz-like equations for a specified frequency, then the issue of whether the modelling is relevant under time evolution need not be considered.

\subsection{Phase-shift embedding}
\label{secpse}

Our powerful innovative alternative approach to homogenisation is to firstly embed the specific given physical \pde~\eqref{Ehdifpde} into a family of \pde\ problems formed by the ensemble of all phase-shifts of the periodic microscale \cite[]{Roberts2013a, Smyshlyaev2000}.
This embedding is a novel and rigorous twist to the concept of a Representative Volume Element.
Secondly, \cref{secsmh} innovatively applies \im\ theory to then prove multi-continuum homogenisations.

Let phase~\(\phi\in[0,\ell)\) parametrise the phase-shifts of the microscale heterogeneity.   
Then our aim is to solve the ensemble of problems\footnote{The time evolution operator~\(\Dn t\alpha{}\), usually written~\alphaD\ herein, covers many cases: 
\(\alpha=1\) is a prototypical diffusion problem; 
\(\alpha=2\) is a prototypical wave problem;
fractional~\(\alpha\) could represent the fractional calculus operator \protect\cite[e.g.,][]{HongGuangSun2018};
and potentially~\alphaD\ could represent any in a wide variety of time evolution operators that commute with spatial derivatives, such as the time-step operator \(\alphaD u\mapsto (u|_{t}-u|_{t-\tau})/\tau\).  
}
\begin{equation}
\Dn t\alpha{u_\phi}=\D x{}\left\{ \kappa(x+\phi)\D x{u_\phi} \right\},
\quad x\in\XX\,,
\label{Eshdifpde}
\end{equation}
for corresponding solution fields~\(u_\phi(t,x)\).
We do this by notionally wrapping each solution~\(u_\phi\) along a wrapped diagonal in a 2D `spatial' strip, of width~\(\ell\), as illustrated by the solid blue line of \cref{Fpattensemble}.  
Then for every solution~\(u_\phi\) the value of coefficient~\(\kappa\) is only a function of~\(\theta\), as indicated by \cref{Fpattensemble}'s coloured background.
Upon doing this for all phases~\(\phi\) a field~\(\fu(t,\x,\theta)\) is defined inside the strip~\([0,L]\times[0,\ell)\) (the fraktur~\fu\ is distinct from the math-italic~\(u\)).
The field~\fu\ must be \(\ell\)-periodic in~\(\theta\) as~\fu\ arises from wrapping the continuous fields~\(u_\phi(t,x)\).

\begin{figure}\centering
\caption{\label{Fpattensemble}cylindrical domain of the embedding \pde~\eqref{Eemdifpde} for field~\(\fu(t,\x,\theta)\): the background colour represents the microscale variation in~\(\kappa(\theta)\).  
Obtain solutions of the heterogeneous \pde~\cref{Eshdifpde} on the blue line as \(u_\phi(t,x) := \fu(t,x,x+\phi)\) for any constant phase~\(\phi\).}
\setlength{\unitlength}{0.01\linewidth}
\begin{picture}(96,26)
\put(5,5){\includegraphics[width=89\unitlength ,height=16\unitlength ] {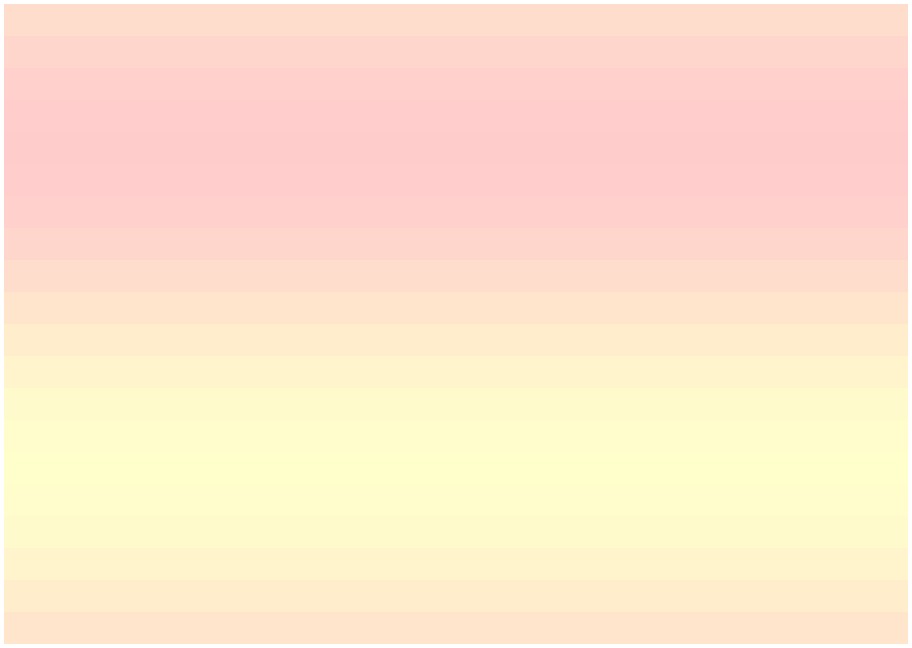}}
\put(0,5){\vector(1,0){96}}
\put(40,2){\(x,\x\)}
\put(5,0){\vector(0,1){26}}
\put(2,24){\(\theta\)}
\put(2,20){\(\ell\)}\put(4,21){\line(1,0){2}}
\put(3,3){\(0\)}
\put(93,2){\(L\)}
\put(35,16){domain \(\dom:=\XX\times[0,\ell)\)}
\put(37,10){\(\fu (t,\x,\theta)\)}
\thicklines
\put(5,5){\framebox(88,16){}}
\color{blue}
\multiput(13,5)(16,0)5{
  \put(0,0){\line(1,1){16}}
  \multiput(0,0)(0,2)8{\line(0,1){1}}
  }
  \put(13,21){\line(-1,-1){8}}
  \put(2,12){\(\phi\)}
\color{magenta}
\put(51,23){\vector(1,0){35}}
\put(47,23){\vector(-1,0){35}}
\put(48,22){$\XX$}
\put(90,13){\oval(8,20)}
\put(87,7){\rotatebox{90}{\parbox{5em}{\raggedright boundary layer}}}
\put(8,13){\oval(8,20)}
\put(5,7){\rotatebox{90}{\parbox{5em}{\raggedright boundary layer}}}
\end{picture}
\end{figure}

Physical boundary conditions at \(x=0,L\) typically require solutions with boundary layers near \(x=0,L\), of microscale size as schematically shown by \cref{Fpattensemble}, in which the solutions vary markedly different to the solution structure in the interior. 
Our modelling is restricted to a spatial domain~\XX\ in an interior that excludes such boundary layers (but see discussion in \cref{Sbcmm}).

Consequently, we create the desired embedding for the \pde~\eqref{Ehdifpde} by considering a field~\(\fu(t,\x,\theta)\) satisfying the \pde
\begin{equation}
\Dn t\alpha\fu=\left(\D \x{}+\D \theta{}\right)
\left\{\kappa(\theta)\left(\D \x\fu+\D \theta\fu \right)\right\},
\quad\fu\text{ \(\ell\)-periodic in }\theta,
\label{Eemdifpde}
\end{equation}
in the `cylindrical' domain \(\dom:=\{(\x,\theta): \x\in\XX\C
\theta\in\Theta:=[0,\ell)\}\) as illustrated in \cref{Fpattensemble}. 
We assume the lower bound \(\kappa(\theta) \geq \kappa_{\min} >0\), and that the \(\ell\)-periodic heterogeneity~\(\kappa(\theta)\) is regular enough to satisfy \cref{assKappa}.
Then we assume that general solutions~\fu\ of \pde~\cref{Eemdifpde} are in~\(\HH^N_\dom:=W^{N+1,2}(\XX)\times\HH\) for some chosen order~\(N\), for the usual Sobolev space, and for \HH\ defined by \cref{assKappa}. 
For example, choose \(N=2\) for the classic~\cref{Eedifpde} or tri-continuum~\cref{Epde3} homogenisations, or choose \(N\geq4\) for higher-gradient homogenisations (e.g., \cref{secsrh,Shcegbmm,SSShoh}).
I emphasise that, in contrast to most extant analytic homogenisations, we do \emph{not} take the usual scale separation limit \(\ell/L\to0\), but instead keep to the real physics case of finite~\(\ell/L\), and so the domain~\(\dom\) of \pde~\cref{Eemdifpde} has finite aspect ratio~\(\ell/L\) in \(\x\theta\)-space as shown by \cref{Fpattensemble}.

\begin{lemma}\label{lemeqv1}
For every solution~\(\fu(t,\x,\theta)\in\HH^N_\dom\) of the embedding
\pde~\eqref{Eemdifpde}, and for every phase~\(\phi\), the defined field \(u_\phi(t,x):=\fu(t,x,x+\phi)\) satisfies the heterogeneous diffusion \pde~\cref{Eshdifpde}.
Specifically, for the case of phase~\(\phi=0\), the field \(u(t,x):=u_0(t,x)=\fu(t,x,x)\) satisfies the given heterogeneous \pde~\eqref{Ehdifpde}.
\end{lemma}

\begin{proof}
The proofs for \cref{lemeqv1,lemcon1} are straightforward, and are encompassed by the proofs given for \cref{lemeqv,lemcon}.
\hfill\end{proof}

Note the implicit distinction among \(x\)-derivatives: \(\D \x{}\) of~\fu\ is done keeping~\(\theta\) constant; whereas \(\D x{}\) of~\(u\) and~\(u_\phi\) is done keeping phase~\(\phi\) constant.
Further, although the embedding field~\fu\ is \(\ell\)-periodic in~\(\theta\), the physical fields~\(u,u_\phi\) are generally \emph{not} periodic in~\(x\).
The microscale \(\ell\)-periodic boundary conditions are \emph{not} assumed but arise naturally due to the ensemble of phase-shifts. That is, what previously had to be assumed, here arises naturally.

\begin{lemma}[converse]\label{lemcon1}
Suppose we have a set of solutions~\(u_\phi(t,x)\) of the phase-shifted \pde~\eqref{Eshdifpde}---a set parametrised by the phase vector~\(\phi\)---and the set depends smoothly enough upon~\(t,x,\phi\) that the following \(\fu\in\HH^N_\dom\). 
Then the field \(\fu(t,\x,\theta):=u_{\theta-\x}(t,\x)\) satisfies the embedding \pde~\eqref{Eemdifpde}.
\end{lemma}

Consequently, \pde{}s~\eqref{Eemdifpde,Eshdifpde} are equivalent, and they may provide us with a set of solutions for an ensemble of materials all with the same heterogeneity structure, but with the structural phase of the material shifted through all possibilities. The critical difference between the \pde{}s~\cref{Ehdifpde,Eshdifpde} and the embedding \pde~\cref{Eemdifpde}  is that although \pde{}s~\eqref{Ehdifpde,Eshdifpde} are heterogeneous in space~\(x\), the embedding
\pde~\eqref{Eemdifpde} is \emph{homogeneous} in~\(\x\).

After establishing this embedding, the distinction between space location~\(x\) and~\x\ is largely irrelevant, and as \(x=\x\) hereafter, hereafter we just use~\(x\).

\subsection{Invariant manifolds of multi-modal any-order homogenisation}
\label{secsmh}

We now analyse the embedding \pde~\eqref{Eemdifpde} for useful \im\  models via an approach proven for systems homogeneous in the long-space \(x\)-direction. 
Recall Manfred Eigen's famous comment ``A model \ldots\ might be right but irrelevant'': the theory of \im{}s  commonly assures us that the models we construct are indeed relevant.
We use the notion of \im{}s  from nonlinear dynamical systems because for a lot of dissipative systems we can prove that \emph{all solutions} in a domain of interest are quickly attracted, on a microscale time, to solutions of the \im\  model.%
\footnote{It is not enough for all trajectories to just approach an invariant manifold:  examples exist \protect\cite[e.g.,][Example~4.11]{Roberts2014a} where the model on such an attractive manifold makes predictions with~\Ord{1} errors.}
Also, for a lot of non-dissipative, wave-like systems we can reasonably establish that \emph{all solutions} near the \im\ evolve close to solutions of the \im\  model.  
Let's define these notions more explicitly.

\begin{definition}[invariant manifold (\im)] \label{Dim}
Let (\dag) denote a given physical (non-autonomous) dynamical system, say \(\alphaD u=f(t,u)\), for a field~\(u\), such that the trajectory from initial condition \(u=u_0\) at time \(t=0\) is the flow \(u=\Phi(t;u_0)\).
Let~\(\UU(t)\)\ be a family, continuous in~\(t\), of open connected subsets of physical interest in the state space of~\(u\) over an open connected set of times~\TT\ of physical interest (\(0\in\TT\subseteq\RR\)).
\needspace{3\baselineskip}
\begin{itemize}
\item Define an \emph{invariant set}~\(\cM\), of~(\dag), to be a union of trajectories of~(\dag) in~\(\UU(t)\):  that is, for some chosen open set~\(\UU_0\subseteq\UU(0)\), \(\cM(t) := \UU(t)\cap \{\Phi(t;u_0): u_0\in\UU_0\}\) for every \(t\in\TT\) \cite[e.g.,][p.121]{Roberts2014a}.

\needspace{3\baselineskip}
\item An \emph{invariant manifold}~\cM, of~(\dag), is an invariant set of~(\dag) such that, for some natural number~\(M\), at every \(u\in\cM(t)\) and \(t\in\TT\) there exists an \(M\)-dimensional tangent space to~\(\cM(t)\)\ (the tangent spaces must be at least continuous in~\(u\in\cM(t)\) and \(t\in\TT\)).

\needspace{4\baselineskip}
\item An \emph{emergent} invariant manifold~\cM\ of~(\dag) is one where there exists a `microscale' timescale~\(t_\ell>0\) and constant~\(C\), such that for every initial condition~\(u_0\in\UU_0\) there exists an initial condition~\(u_0^\cM\in\cM(0)\) for which \(\|\Phi(t;u_0)-\Phi(t;u_0^\cM)\|\leq C\e^{-t/t_\ell}\|u_0-u_0^\cM\|\).   
That is, for every trajectory of~(\dag) in~\UU(t), there exists a trajectory in~\cM(t)\ which is approached exponentially quickly in time.\footnote{%
Of course, this requirement only makes useful practical sense if the `microscale' time~\(t_\ell\) is usefully much smaller than the timescale~\(|\TT|\) of physical interest.  
In some applications it can make sense to replace the exponential approach to~\cM\ by an algebraic-in-time approach, but then the theoretical support is much more delicate.} 

\needspace{5\baselineskip}
\item A \emph{guiding-centre} invariant manifold~\cM\ \cite[e.g.,][\S11]{vanKampen85}, of~(\dag),  is one where there exists a `microscale' timescale~\(t_\ell>0\) and constants~\(C,r_0\), such that for every~\(0<r\leq r_0\), and for every initial condition~\(u_0\in\{u_0\in\UU_0: \operatorname{dist}(u_0,\cM)<r\}\) there exists an initial condition~\(u_0^\cM\in\cM(0)\) for which \(\de t{}\|\bar\Phi(t;u_0)-\Phi(t;u_0^\cM)\|\leq Cr^2\) upon defining the microscale averaged-trajectory \(\bar\Phi(t;u_0):=\frac1{t_\ell}\int_t^{t+t_\ell}\Phi(t;u_0)\d t\)\,.

\end{itemize}
For the common case of autonomous systems, as in the rest of \cref{SoneDintro,Shceq}, neglect the \(t\)-dependence in~\UU\ and~\cM.
\end{definition}

\begin{example}[emergent invariant manifold]\label{EGeim}
\begin{subequations}\label{EEeim}%
In variables~\(u(t)\) and~\(v(t)\) consider the following toy dynamical system \cite[\S1.2]{Roberts2018a} for which, linearly, \(v(t)\)~decays like~\(\e^{-t}\):
\begin{equation}
	\de t u=-uv \qtq{and} 
	\de t v=-v+u^2-2v^2.
	\label{Ecmect}
\end{equation}
Firstly, straightforward algebra shows that both \(u=0\) and \(v=u^2\) are invariant under the \ode{}s~\cref{Ecmect}, and hence both this line and curve are invariant sets (as is, incidentally, their union and intersection).
Secondly, these two separate sets have smoothly varying tangent spaces, and hence they are both invariant manifolds~(\im{}s), here of~1D.
Thirdly, the marvellous coordinate transform
\begin{equation}
u=\frac{U}{\sqrt{1-2V/(1+2U^2)}}
\quad\text{and}\quad
v=U^2+\frac{V}{{1-2V/(1+2U^2)}}\,,
	\label{EctECT}
\end{equation}
converts the \ode{}s~\cref{Ecmect} into the system
\begin{equation}
\de t U=-U^3
\quad\text{and}\quad
\de t V=-V\left[\frac1{1+2U^2}+4U^2\right].
	\label{EcmECT}
\end{equation}
The form of~\cref{EctECT} shows the two \im{}s are mapped to the subspaces \(U=0\) and \(V=0\) respectively.  
Further, the form of the \(V\)-\ode\ in~\cref{EcmECT} leads to the bound that \(|V|\leq|V_0|\e^{-t}\) for \(t\geq0\) and every~\((U_0,V_0)\).
That is, \(V=0\) is exponentially quickly attractive, and hence so is \(v=u^2\) under~\cref{Ecmect}.
Moreover, since the \(U\)-\ode\ is independent of~\(V\), for every initial state~\(U_0,V_0\) the subsequent trajectory has precisely the same \(U\)-evolution as the trajectory from the initial state~\((U_0,0)\) on the \im~\(V=0\), and hence approaches this \im's trajectory exponentially quickly.  
Such emergence of the \im~\(V=0\), via the transform~\cref{EctECT}, establishes the emergence of \text{the \im~\(v=u^2\) under~\cref{Ecmect}.}
\end{subequations}
\qed\end{example}

The existence and emergence of \im{}s likewise express and support the relevance of a precise multi-modal homogenisation of the original heterogeneous \pde~\eqref{Ehdifpde}. 
Since, except for \cref{SShhn}, the example \pde{}s herein are linear the \im{}s  are mostly more specifically invariant \emph{subspaces}, but let's use the term \emph{manifold} as the same framework and theory immediately generalises to related nonlinear systems as discussed in \cref{Sgentheory} and the example of \cref{SShhn}.

Theory \citep{Roberts2013a} inspired by earlier more formal arguments \citep{Roberts88a, Roberts96a} establishes how to support and construct \pde\ models for the macroscale spatial structure of \pde\ solutions in cylindrical domains such as the strip~\dom. 
The practical technique is to base analysis on the foundational case where variations in~\(x\) are negligible, and then treat gradual, macroscale, variations in~\(x\) as a regular perturbation \cite[Part~III]{Roberts2014a}.%
\footnote{Alternatively, in linear problems one could justify the analysis via a spatial Fourier transform in~\(x\) \cite[\S2 and Ch.7 resp.]{Roberts88a, Roberts2014a}, as is also done in some alternative homogenisation methods \cite[e.g.,][]{Willis1985, Dohnal2015}.
However, for nonlinear problems, and also for macroscale varying heterogeneity, it is better to analyse in physical space, so we do so herein.}
Consequently, to establish the foundation of an \im~model, we consider the embedding \pde~\eqref{Eemdifpde} with \(\D{x}{}\) neglected:
\begin{equation}
\Dn t\alpha\fu=\D \theta{}\left\{\kappa(\theta)\left(\D \theta\fu \right)\right\},
\quad\fu\text{ is \(\ell\)-periodic in }\theta.
\label{Eemdifpde0}
\end{equation}
This \pde\ is a classic cell-problem in homogenisation.  
That we also use this cell-problem as a foundation means that the approach here agrees with well-established homogenisations---but the approach here extends homogenisations to much more general scenarios and to enhanced space-time scales.
The foundational cell-problem~\cref{Eemdifpde0} applies in a locale around each and every~\(x\in\XX\) \cite[\S2]{Roberts2013a}.
In general, as in some functionally graded materials and as addressed in \cref{Sgentheory}, the details of such a basis varies with locale \(x\in\XX\) and such variation is encompassed by the approach here.
But in this section, because the \pde{}s~\cref{Eemdifpde,Eemdifpde0} are translationally invariant in~\(x\), the following basis is \text{independent of~\(x\in\XX\)\,.}

\paragraph{Equilibria}
Invariant manifolds are most easily detected and constructed in a region near an equilibrium\footcite[e.g.,][]{Carr81, Chicone97, Haragus2011}, as elaborated further in \cref{SSSlbim}.  
For example, suppose a given dynamical system \(\alphaD u=f(u)\), such as one in the class~\cref{Egenpde}, has equilibrium at say \(u=0\), and its linearisation near the equilibrium is \(\alphaD u=\cL u\) for a linear operator~\cL.   
Then let \(\{v_m\}\) denotes a complete set of (generalised) eigenvectors of~\cL\ that span the state space of~\(u\).
For every chosen subset~\MM\ of indices the subspace \(\cE_\MM:=\Span\{v_j:j\in\MM\}\), often called a \emph{spectral subspace}, is an \im\ of the linearised system \(\alphaD u=\cL u\)\,.  
Nonlinear dynamical systems theory proves that, under suitable conditions, nonlinearity and/or parameter perturbations in~\(f(u)\) lead to smooth regular perturbations of the chosen subspace~\(\cE_\MM\) in a finite domain about the equilibrium.
Thus, an \im\ of \(\alphaD u=f(u)\), tangent to~\(\cE_\MM\), exists and may be constructed.

Consequently we start with equilibria of the cell problem~\cref{Eemdifpde0}. 
A family of equilibria of \pde~\cref{Eemdifpde0} is, for every constant~\(C\), \(\fu(t,\theta)=C\). 
Because the \pde~\cref{Eemdifpde0} is linear, to encompass the entire family it is sufficient to consider just the case of equilibrium \(C=0\), which we do henceforth.

\paragraph{Spectrum at each equilibrium}
Invariant manifold models \citep{Roberts2013a} are decided based upon the spectrum of the cell problem~\cref{Eemdifpde0}. In general the spectrum depends upon the microscale details of~\(\kappa(\theta)\).  

\begin{assumption}\label{assKappa}
From~\eqref{Eemdifpde0}, consider the `cell' eigen-problem for~\(\fu\)
\begin{equation}
\lambda\fu=\D \theta{}\left\{\kappa(\theta)\D \theta\fu \right\},
\quad\fu\text{ is \(\ell\)-periodic in }\theta.
\label{Eemdifeig}
\end{equation}
We assume that \(\kappa(\theta) \geq \kappa_{\min}>0\) is regular enough that the eigenvalues~\(\lambda\) are countable, real and non-positive, and also that a set of corresponding eigenvectors are complete, orthogonal, and span a weighted Sobolev space, denoted~\(\HH\), that is physically meaningful.
For generality, all derivatives may be interpreted in a weak sense, and \(\HH := \{\fu\in L^2([0,\ell]): \|\fu\|^2+\|\kappa^{1/2}\partial_\theta \fu\|^2<\infty \C \fu\text{ is }\ell\text{-periodic} \}\).%
\footnote{Physically, requiring \(\kappa_{\min}>0\) is not necessarily necessary.  An example would be a variation of \protect\cref{EGsd} (\protect\cref{Sgentheory}) to the homogenisation of shear dispersion in \emph{meandering} 'turbulent' channel flow with, in a 2D channel \(|z|<1\), a mean-flow stream-function say \(\psi=(1-z^2)[1+az\sin(2\pi x/\ell)]\) and an eddy diffusivity \(\kappa\propto 1-z^2\) that \emph{goes to zero} at the edges of the channel and hence goes to zero \text{at some cell edges.}}
\end{assumption}

For the non-dimensional \cref{EG1D}, the heterogeneity~\cref{EGkappa} is parametrised by~\(a\) (\(|a|<1\)).
In such cases we may base analysis on any convenient parameter value, here we choose to base analysis on~\(a=0\).
For this base case of constant \(\kappa=1\)\,, the spectrum of the cell eigen-problem~\cref{Eemdifpde0} is the set of eigenvalues \(\{0\C -1\C -1\C -4\C -4\C -9\C \ldots\}\) (non-zero eigenvalues have multiplicity two). 
The corresponding set of complete and orthogonal eigenvectors are~\(\{1\C \sin\theta\C \cos\theta\C \sin2\theta\C \cos2\theta\C \sin3\theta\C \ldots \}\).  
By continuity in the self-adjoint cell eigen-problem~\cref{Eemdifeig}, for at least a finite range of heterogeneity \(a\neq 0\)\,, the spectrum and eigenvectors are similar.

In general, let's order the eigenvalues such that \(0=\lambda_0 >\lambda_1 \geq\lambda_2 \geq\lambda_3 \geq\cdots\) (including repeats to account for multiplicity).
Let \(v_m(\theta)\)~denote an eigenvector corresponding to the eigenvalue~\(\lambda_m\) (suitably orthogonormalised in the case of eigenvalues of multiplicity two or more).

From the family of equilibria, \(\fu(t,\theta)=C\), we know the leading eigenvalue \(\lambda_0=0\) and the corresponding~\(v_0(\theta)\) is constant, say normalised to \(v_0=1\).
One may construct a \emph{slow invariant manifold} model based upon this eigenvalue zero.   
Such slow \im\ modelling encompasses the classic homogenised \pde\ such as~\cref{Eedifpde}, as well as its higher order generalisations.
The reason for this connection to classic homogenisation is that the constant eigenvector~\(v_0\) both matches the classic \emph{assumption} that the macroscale solutions varies little over a cell, but also matches the classic \emph{assumption} that unweighted cell averages give familiar macroscale quantities.  
Here, both such properties instead \emph{follow} from the physics-informed nature of the leading microscale eigenvector~\(v_0(\theta)\). 

In problems with more complicated physics, the correct corresponding properties follow from the physical nature of the leading eigenvector: an example is modelling the macroscale advection-diffusion in field flow fractionation channels where the leading microscale eigenvector is an exponential structure;\footcite[e.g.,][]{Suslov00a, Suslov98c} another example is the analysis of Fokker--Planck \pde{}s to construct a model for its marginal distribution where the leading eigenvector is a Gaussian structure.\footcite[e.g.,][]{Konno1994, Arnold95, Kuehn2025a}

\subsubsection{Multi-modal, multi-continuum, homogenisations exist}
\label{SSSmmmcme}

A rational multi-modal, invariant manifold, homogenisation is justified and constructed based upon the \(M\)~leading eigenvalues and eigenvectors of the cell-problem~\cref{Eemdifeig}, for any chosen~\(M\).  
This is seen most clearly in the dissipative case \(\alpha=1\) that we now discuss.\footnote{For the wave case of \(\alpha=2\), supporting theory comes from the sub-centre manifolds of \protect\cite{Sijbrand85} [\S7] or the spectral submanifolds of \protect\cite{Cabre2003}, or the backwards theory by \protect\cite{Hochs2019}.}

\begin{subequations}\label{EE1stApp1}%
\begin{corollary}[first approximation]\label{C1stApp1}
The first approximation to an \(M\)-mode, multi-continuum, invariant manifold, homogenised model is the following general linear combination of the corresponding sub-cell eigenvectors
\begin{equation}
\fu(t,x,\theta)\approx U_0v_0(\theta)+U_1v_1(\theta)+\cdots+U_{M-1}v_{M-1}(\theta),
\label{EbaseIM}
\end{equation}
for some coefficients \(U_0,\ldots,U_{M-1}\).   
The corresponding evolution is that
\begin{equation}
\alphaD U_0=\lambda_0U_0  \C\ldots\C 
\alphaD U_{M-1}=\lambda_{M-1}U_{M-1}\,.
\label{E1stApp1}
\end{equation} 
\end{corollary}
\end{subequations}

\begin{proof} 
Sections~2.1--2.3 by \cite{Roberts2013a} establishes that an \im\ homogenisation \emph{exists} based upon the spectral properties of the linear cell problem~\cref{Eemdifeig}.
Then, for example, \cite{Haragus2011} Hypothesis~2.4 identifies the set of centre eigenvalues which \cite{Aulbach2000} generalised, in their \S4 Hypothesis~B1, to include  non-zero eigenvalues as in the case here.\footnote{Here the spectral bounds of \protect\cite{Aulbach2000}  would be \(\alpha_1:=\lambda_M\), \(\beta_1:=\lambda_{M-1}\), \(\alpha_2:=\lambda_0=0\)\,.}
Then (2.2) of \cite{Haragus2011} leads to defining the corresponding \(M\)-D (generalised) centre subspace, called~\(\cE_0\).
Subsequently, their Theorem~2.9 proves that there exists a \im\  \(\cM_0:=\{u_0+\Psi(u_0): u_0\in\cE_0\}\) for some strictly nonlinear~\(\Psi\).
That is, the first approximation to the \im\ is the linear subspace~\(\cE_0\) as it is the tangent space to~\(\cM_0\) at~\(u_0\).
A general expression for any point in the tangent space~\(\cE_0\) is then the generalised eigenspace~\cref{EbaseIM}, which physically is then a tangent space approximation to macroscale varying sub-cell structures.
The first approximation to the evolution on~\(\cM_0\) is then~\cref{E1stApp1} from the corresponding eigenvalues \cite[]{Aulbach2000}, which then serves as a first approximation to an \(M\)-mode homogenisation of the embedding \pde~\cref{Eemdifpde}, and thence \text{to the original \pde~\cref{Ehdifpde}.}
\hfill\end{proof}

Due to the existence of an in-principle exact closure, that a \im\ field \(u=u_0+\Psi(u_0)\) in the notation of \cite{Haragus2011}, this \im\ approach also satisfies the second criteria of \cite{Auriault2009} [p.59] that ``the ideal [homogenisation] procedure must also permit \ldots\ the determination of the local fields of physical quantities starting from the values of macroscopic physical quantities.''  
The construction of the \im~\(u_0+\Psi(u_0)\), such as~\cref{EG1Dv1}, providesthe local fields  \text{(for all heterogeneity phase-shifts).}

Physically, \(U_0,\ldots,U_{M-1}\) are macroscale `variables', `amplitudes' or `order parameters' that vary acceptably gradually over the macroscale~\(x\)  \cite[]{Roberts2013a}, that is, they are functions~\(U_0(t,x)\C \ldots\C U_{M-1}(t,x)\).
For eigenvalues of multiplicity\({}>1\) it is most appealing to choose corresponding eigenvectors~\(v_j(\theta)\) with distinct physical meaning as then the corresponding variables~\(U_j\) have a clear \text{physical meaning.}\footnote{In principle, \cref{EbaseIM} only needs to span the chosen 'centre' subspace~\(\cE_0\), so any linearly independent set of \(M\)~linear combinations of the eigenvectors could instead be chosen for basis functions in~\cref{EbaseIM}.  
However, solving the homological equation~\cref{E1dHomologic} for the necessary corrections to~\cref{EbaseIM} then becomes significantly more complicated.}

For \cref{EG1D}, choosing to resolve modes corresponding to the three eigenvalues \(\lambda_0=0\) and \(\lambda_1=\lambda_2=-1\) leads to the tri-mode, tri-continuum, homogenised model~\cref{Epde3}. 
In the terminology of \cite{Haken83, Chouksey2023}, modes with eigenvalues\({}\leq \lambda_3= -4\) are slaved.
From nonlinear dynamical system theory, the spectral gap~\((-1,-4)\) between~\(\lambda_2\) and~\(\lambda_3\) (\cref{Dgap}) caters for the perturbing heterogeneity and perturbing \(x\)-gradients of~\fu.\footnote{In nonlinear problems, this spectral gap also caters for the perturbing nonlinearity.}  
This tri-mode choice leads to the (tangent space) approximation~\cref{EbaseIM} being here \begin{equation}
\fu\approx U_0 +U_1\sin\theta +U_2\cos\theta
\C\quad \alphaD U_0=0
\C\  \alphaD U_1=-U_1
\C\  \alphaD U_2=-U_2
\,.
\label{Eegim3}
\end{equation}
Hence the \im\ parameters have the following physical interpretation: 
\(U_0\)~measures the mean\slash average of field~\(u\) in each cell;
\(U_1\)~measures the leading sub-cell structures in~\(u\) that are \emph{out-of-phase} with the material heterogeneity~\eqref{EGkappa}; and 
\(U_2\)~measures leading structures \emph{in-phase} with the heterogeneity.
Because the leading neglected mode has corresponding eigenvalue \(\lambda_3=-4\), this tri-mode homogenisation will resolve dynamics on (non-dimensional) time scales significantly longer than \(1/|\lambda_3|^{1/\alpha} =1/4^{1/\alpha}\).
In contrast, many other homogenisation approaches do not explicitly identify such a lower bound on \text{the timescale resolution.}

The task of \cref{SSScmmm} is to establish that this dynamical systems framework systematically derives physics-informed refinements to~\cref{EbaseIM} (the function~\(\Psi(u_0)\) in the proof of \cref{C1stApp1}) and simultaneously derives accurate \pde{}s, such as~\cref{Eedifpde,Epde3} for \cref{EG1D}, governing the evolution of the macroscale variables~\(U_0,\ldots,U_{M-1}\).

Importantly, there are only \emph{two} subjective decisions made in this approach.
The first subjective decision is where to divide the spectrum into sub-cell modes whose dynamics we model explicitly, namely~\(v_0,\ldots,v_{M-1}\), and sub-cell modes which are accounted for implicitly, which are `slaved', namely~\(v_m\) for \(m\geq M\)\,. 
\cref{SSSsd,SSSwls} discuss that an \(M\)-mode model should resolve time scales significantly longer than~\(1/|\lambda_M|^{1/\alpha}\).
That is, this first decision is largely about resolved time scales.
The second subjective decision is to choose an order~\(N\) of accuracy for the constructed \(M\)-continuum model (e.g., \cref{SSScmmm}) which can improve the spatial resolution.
For \cref{EG1D}, the given classic homogenisation~\cref{Eedifpde} neglects third-order spatial derivatives and so corresponds to choosing order \(N=2\).
\cref{SSSisr} discusses subtleties in choosing an appropriate order~\(N\).
A rough approach is to construct a few more orders than you expect, then in use truncate to an order~\(N\), and estimate its errors by the magnitude of the leading neglected term, whether of order~\(N+1\) or~\(N+2\).
This would serve as a rough estimation of the quantitative \text{remainder error~(23) of \cite{Roberts2013a}.}


\subsubsection{Multi-modal, multi-continuum, homogenisations are relevant}

The justification for the relevance of an homogenisation founded on~\cref{EbaseIM} is the following. 
However, the argument depends upon the nature of the time evolution operator~\alphaD\ (\cref{Simmmh} discusses more broadly with detailed justification).
\begin{itemize}
\item In the diffusive case, \(\alpha=1\), the homogenisation is relevant because the slow \im\ tangent to~\cref{EbaseIM} exponentially quickly attracts solutions from all initial conditions\footcite[e.g.,][]{Aulbach2000, Prizzi02, Roberts2013a}.   
This emergence (\cref{Dim}) is because all the slaved sub-cell modes decay roughly like~\(\e^{\lambda_mt}\) for \(m=M,M+1,\ldots\)\,.
The slowest of these is~\(\e^{\lambda_Mt}\) and so we expect, and can often prove (e.g., \cref{EGeim}), that solutions from all initial conditions approach a corresponding \(M\)-mode \im\ on fast times, roughly~\(1/|\lambda_M|\).

\item In the wave case, \(\alpha=2\), all sub-cell eigenvector modes are oscillatory with frequency \(\omega_m:=\sqrt{-\lambda_m}\)\,.
Hence any \im\  founded on~\cref{EbaseIM} appears not to be emergent in time.
Instead, one may argue its relevance via one of at least three sub-cases:
\begin{itemize}
\item commonly one views the model as a a \emph{guiding centre} (\cref{Dim}) for the dynamics on timescales longer than~\(1/\omega_M\) about the constructed manifold \cite[e.g.,][]{vanKampen85}, albeit despite controversies about the existence of such slow \im\ models \cite[e.g.,][\S13.5.3]{Lorenz87, Roberts2014a}, controversies that may be resolved via backwards theory \cite[e.g.,][]{Hochs2019};

\item or there may be physical processes not represented in the given sub-cell dynamics~\cref{Eemdifpde0} that physically cause sufficient attraction to the chosen \(M\)-mode subspace~\cref{EbaseIM};
\item or one is only interested in predicting equilibria or periodic orbits in which case emergence is largely irrelevant.
\end{itemize}

\end{itemize}

\subsubsection{Construct multi-modal, multi-continuum, homogenisations}
\label{SSScmmm}

The construction of multi-continuum homogenisation relies on theory proven by \cite{Roberts2013a}, which in turn rests on general theory by \cite{Aulbach2000, Potzsche2006, Hochs2019}.
Corollary~13 of \cite{Roberts2013a} proves that an established procedure \citep{Roberts88a, Roberts96a} is indeed a rigorous method to construct \im\  homogenised models such as~\eqref{Epde3}.
The crucial result is that if a derived approximation satisfies the embedding \pde~\cref{Eemdifpde} to a residual of~\Ord{\partial_x^{N+1}}, then the corresponding homogenisation is correct to an error~\Ord{\partial_x^{N+1}}. 
For a class of linear advection-diffusion systems, \cite{Watt94b} developed practical procedures to derive multi-mode models to any chosen order~\(N\) of residual.
These procedures are extended herein to \cref{Pconstruct} that encompasses a much wider class of homogenisation problems, \text{including nonlinear problems.}

The procedure is based upon the following.
For a defined vector of local amplitudes~\(\Uv(t,x) := (U_0,\ldots,U_{M-1})\),
we seek an \im\  of the embedding \pde~\cref{Eemdifpde} in the form \(\fu(t,x,\theta)=v(\Uv,\theta)\) such that \(\alphaD \Uv=G(\Uv)\) where the right-hand side dependence upon~\Uv\ in both of these implicitly involves its gradients~\(\Uv_x\C \Uv_{xx}\C\)etc, as in~\eqref{Epde3}.
For any given approximations~\(\tv,\tG\) to~\(v,G\), define~\(\Res(\tv,\tG)\) to be the residual of \text{the embedding \pde~\cref{Eemdifpde}.}

\begin{lemma} \label{Lupdates}
For the embedding \pde~\cref{Eemdifpde}, compute corrections~\(v',G'\) to an \im\ approximation~\(\tv,\tG\) by solving a variant of the usual linear cell problem forced by the residual, namely
\begin{equation}
\D\theta{} \left\{ \kappa(\theta) \D\theta{v'} \right\}
-\sum_{m=0}^{M-1}\lambda_m\D{U_m}{v'}U_m
-\sum_{m=0}^{M-1}v_mG'_m
=\Res(\tv,\tG),
\label{E1dHomologic}
\end{equation} 
often called the \emph{homological equation}.%
\footcite[e.g.,][]{Potzsche2006, Roberts2014a, Siettos2021, Martin2022} 
Interpret the factor \((\D{U_m}{v'})U_m\) in the Fre\'chet derivative, Calculus of Variations, sense that it represents \(v'_{U_m}U_m +v'_{U_{mx}}U_{mx} +v'_{U_{mxx}}U_{mxx} +\cdots\) where these subscript-derivatives of~\(v'\) are done with respect to the subscript symbol \cite[]{Roberts88a}.
\end{lemma}

\begin{proof} 
Obtain the homological equation~\cref{E1dHomologic} straightforwardly by the following \cite[e.g.,][\S2.1]{Roberts96a}.
Substitute correcting approximations \(\fu=\tv+v'\) and \(\alphaD \Uv =\tG+G'\) into the governing \pde~\cref{Eemdifpde}.
Linearise in small corrections (primed).
The terms independent of primed corrections form the residual~\(\Res(\tv,\tG)\).
Approximate the terms linear in primed corrections by evaluating their coefficients at the leading (tangent space) approximation of \cref{C1stApp1}.
Because the \im\ invariance equation is not singular \cite[e.g.,][]{Potzsche2006}, then solving~\cref{E1dHomologic} for~\(v',G'\) ensures \(\Res(\tv+v',\tG+G')=\ord{\Res(\tv,\tG)}\), and so \(\fu=\tv+v'\) and \(G=\tG+G'\) gives an improved approximation to the \text{\im\ and evolution thereon.}
\hfill\end{proof}

The distinction between the homological equation~\cref{E1dHomologic} and the usual cell-problems for homogenisation arises because this \im~framework systematically accounts for all physical out-of-equilibrium effects (because an \im~is composed of actual trajectories of a given physical system, \cref{Dim}).

Nonetheless, Corollary~12 by \cite{Roberts2013a} establishes that a traditional multiple scale homogenisation \cite[e.g., \S3.3.2.1 of][]{Auriault2009} when compared to the \im~homogenisation here will both have the same leading order terms. 
This agreement arises because the algebraic details of the two approaches start out as much the same: for example, the \im\ embedding \pde~\cref{Eemdifpde} has much the same algebraic form as the corresponding \pde\ arising in multiple scales.
The key difference is the framework in which the algebraic details are realised: the \im\ approach is more flexible and powerful, with better supporting theory.
Differences may arise in non-leading order terms depending upon particular assumptions made in order to take the method of multiple scales to higher-order.
Differences are more significant for micromorphic-modes of non-zero~\(\lambda_m\) as the leading order term is then~\(\lambda_mU_m\), and for multiple scales to correctly obtain the interesting spatial gradient terms one needs to extend its heuristics to multiple time-scales and more space-scales that only then properly account for the \text{chain rule within~\(\alphaD u\).}

To construct \im\ homogenisations, the procedure is to thus iteratively evaluate the residual~\(\Res(\tv,\tG)\), and solve~\cref{E1dHomologic} for corrections, until the residual is of the desired \text{order of error (\cref{Pconstruct}).}

\begin{subequations}%
For \cref{EG1D} let's detail the first iteration.
The leading, tangent space, approximation for~\(\tv\C \tilde\Gv\) is~\cref{Eegim3}.
For this approximation, computing the residual of the embedding \pde~\cref{Eemdifpde} gives
\begin{align}&
\Res(\tv,\tG)
=-U_1a\sin2\theta
-U_2a\cos2\theta
-U_{0x}a\sin\theta
\nonumber\\&\quad{}
+U_{2x}(2\sin\theta-\tfrac32a\sin2\theta)
+U_{1x}(\tfrac12a -2\cos\theta +\tfrac32a\cos2\theta)
+\Ord{a^2}.
\label{EG1Drhs}
\end{align}
The homological equation~\cref{E1dHomologic} to solve for (dashed) corrections is then
\begin{equation}
\DD\theta{v'} 
+\D{U_1}{v'}U_1 +\D{U_2}{v'}U_2
-G'_0 -G'_1\sin\theta -G'_2\cos\theta
=\Res(\tv,\tG),
\label{EG1Dhomo}
\end{equation}
where here the left-hand side is approximated by setting parameter \(a=0\).
To solve~\cref{EG1Dhomo} we first, the so-called \emph{solvability condition}, choose~\(G'_0,G'_1,G'_2\) to, respectively, eliminate constant-in-\(\theta\), \(\sin\theta\), and~\(\cos\theta\) components in the residual.
This leads to the physics-informed correction \(G'=(-\tfrac12aU_{1x}\C aU_{0x}-2U_{2x}\C 2U_{1x})+\Ord{a^2}\) that give all the first-order gradient terms listed in the tri-continuum model~\cref{Epde3}.
The second step is to solve the rest of the homological~\cref{EG1Dhomo} for the correction~\(v'\) \text{to the sub-cell field:}
\begin{align}
\DD\theta{v'} 
+\D{U_1}{v'}U_1 +\D{U_2}{v'}U_2
&=-U_1a\sin2\theta
-U_2a\cos2\theta
\notJbreak{}
-U_{2x}\tfrac32a\sin2\theta
+U_{1x}\tfrac32a\cos2\theta
+\Ord{a^2}.
\label{EG1Dvd}
\end{align}
Straightforwardly solving this for \(\ell\)-periodic solutions, say via undetermined coefficients, gives that the sub-cell physics updates the \im\ field from~\cref{Eegim3} to
\begin{align}
\fu&\approx U_0 
+U_1(\sin\theta +\tfrac13a\sin2\theta) 
+U_2(\cos\theta +\tfrac13a\cos2\theta)
\notJbreak{}
+U_{2x}\tfrac12a\sin2\theta
-U_{1x}\tfrac12a\cos2\theta
\,.
\label{EG1Dv1}
\end{align}
This second approximation~\cref{EG1Dv1} to the \im\ illustrates how the `slaved' sub-cell modes, here the \(\sin2\theta\) and \(\cos2\theta\) components) are determined via the homological equations. 
Other multi-continuum approaches would derive similar first updates, but one key difference here is that we systematically account for the out-of-equilibrium physical effects encoded by~\(+\D{U_1}{v'}\,U_1 +\D{U_2}{v'}\,U_2\) in the homological equation~\cref{EG1Dvd}.
In contrast, many other homogenisation approaches omit such out-of-equilibrium physics because they either only address statics, or by invoking traditional average-based heuristics: for example, the comment by \cite{Auriault2009} [p.58] that ``to take the mean of a system of partial differential equations \ldots\ is valid'' overlooks subtle effects of the chain rule for \emph{dynamics} in \(\partial_tv=\sum_m\D {U_m}v\cdot\D t{U_m}\).%
\footnote{For another example, the method of \cite{Milton2007} cannot be exact because they rely on [pp.868, 873 e.g.] ``ensemble averaging the equation of motion'' and ``averaging (6.2) gives''. 
Hence they also implicitly omit nonlinear subtleties of the chain rule for \(\partial_tv=\sum_m\D {U_m}v\cdot\D t{U_m}\).}
\end{subequations}


For \cref{EG1D} further similar iterations give corrections of higher-order in both heterogeneity~\(a\) and spatial gradients~\(\partial_x\).
For example, the computer algebra of \ifJ Supplementary Code \cite[Appendix~A,][]{Roberts2024a} \else\cref{cas} \fi systematically constructs the tri-continuum model to be
\begin{subequations}\label{EE3Mevol}%
\begin{align}
\alphaD U_{0}&=
-\tfrac{1}{2}a U_{1x} 
+\left(1+\tfrac{1}{2} a^{2}\right)U_{0xx} 
-\tfrac{1}{2}a U_{2xx} 
+\Ord{\partial_x^3,a^3},
\label{EE3Mevol0}
\\
\alphaD U_{1}&=
-\left(1+\tfrac{5}{12} a^{2}\right)U_{1} 
+aU_{0x} 
-\left(2+\tfrac{2}{9} a^{2}\right)U_{2x} 
\notJbreak{}
+\left(1-\tfrac{17}{54} a^{2}\right)U_{1xx}
+\Ord{\partial_x^3,a^3},
\label{EE3Mevol1}
\\
\alphaD U_{2}&=
-\left(1-\tfrac{1}{12} a^{2}\right)U_{2} 
+\left(2+\tfrac{2}{9} a^{2}\right)U_{1x} 
\notJbreak{}
-a U_{0xx} 
+\left(1+\tfrac{5}{27} a^{2}\right)U_{2xx} 
+\Ord{\partial_x^3,a^3}.
\end{align}
These \pde{}s form a rigorous second-order tri-continuum, three-mode,
homogenised model for the heterogeneous system~\cref{Ehdifpde} with heterogeneity~\cref{EGkappa}.
Physically, \cref{EE3Mevol1} shows that gradients~\(U_{0x}\) of the mean mode predominantly create out-of-phase microscale structures, measured by~\(U_1\), that then affect the macroscale effective diffusivity of the mean mode via the \(U_{1x}\)-term in~\cref{EE3Mevol0}. 
\end{subequations}

In the above and other `multi-physics' asymptotic statements I invoke the following definition.
\begin{definition}[asymptotic order] \label{Dao}
Define that \Ord{\rho,\psi} means~\(\Ord{\rho}+\Ord{\psi}\), and 
define that~\Ord{\partial_x^p} or~\Ord{\grad^p} is to mean terms with \(p\)~or more spatial derivatives.
More precisely, \Ord{\partial_x^p} and \Ord{\grad^p} are to mean of the order of the corresponding remainder term, respectively~(23) by \cite{Roberts2013a} and~(52) by \cite{Roberts2016a}, for the case \(N=p-1\)\,.
\end{definition}

\subsubsection{Flexible microscale parametrisation}
\label{Sfmp}

Many researchers express concerns that solutions of homological equations, such as~\cref{E1dHomologic}, are not unique---there are more variables than equations.
To resolve such non-uniqueness one \emph{must}, explicitly or implicitly, define a physical meaning for each of the amplitude variables~\(U_0,\ldots,U_{M-1}\).
Such definition makes the solution of~\cref{E1dHomologic} unique \cite[e.g.,][\S2.1 and \S5.3 resp.]{Roberts96a, Roberts2014a}.

For example, this \cref{SoneDintro} \emph{implicitly} chooses to parametrise multi-continuum homogenisations directly in the amplitudes~\(U_0,\ldots,U_{M-1}\) of the microscale sub-cell eigenvectors~\(v_0,\ldots,v_{M-1}\) via the ``elimination'' mentioned following the homological equation~\cref{EG1Dhomo}.  
This common parametrisation is straightforward to do because these are the (orthogonal) microscale eigenvectors, but in principle we should be \emph{explicit} and define the amplitudes as \(U_m:=\tfrac1\pi\int_0^{2\pi}v_m\fu\d\theta\)\,.  
However, the \im\ framework potentially empowers us to parametrise almost arbitrarily the manifold of the homogenisation \cite[\S5.3,e.g.]{Roberts2014a}---akin to general change of basis for a linear operator---and so in your homogenisations you may choose almost any physical meaning you like for your macroscale variables~\(U_0,\ldots,U_{M-1}\): they just have to be able to smoothly parametrise the \im.
For example, in shear dispersion (\cref{EGsd} introduces shear dispersion in the general theory of \cref{Sgentheory}) \cite{Strunin01a} [\S2 and~\S3 resp.] showed how to transform an \im\ \emph{from a two-mode model to a two-zone model} either via transforming the derived modal equations, or equivalently via defining two coupled zones at the outset and deriving the interaction between and within the two zones in terms of their chosen defined \text{averages over each zone.}

\cite{Auriault2009} [p.58] comment that a ``general characteristic of all the [homogenisation] methods is that they use mean values to define macroscopic quantities''.
For micromorphic models, this is only appropriate when interpreted as a \emph{weighted mean}, typically weighted by a corresponding eigenvector (adjoint if not self-adjoint).
But, in contrast, our placing homogenisation in an invariant manifold framework proves that one can parametrise the macroscale by very flexibly defined quantities.  
It is largely up to a user to \emph{choose} what quantities are physically appropriate on the macroscale---of course, most users will choose to invoke \text{(weighted) mean values.}

Alternatively, one could instead \emph{adaptively modify} the definition of the amplitudes to simplify the algebraic form~\cref{EE3Mevol} of the evolution on the invariant manifold---a \emph{normal form} of the homogenisation \cite[e.g.,][]{Arneodo85b}. 
But such adaptive modification is often unphysical, often quite tedious, and it can be hard to achieve desired simplifications.\footnote{%
For example, \cite{Roberts2014a} [p.178] showed that, given two different parametrisations of an \im, if the transformation between the parametrisations commutes with the evolution operator, then the two corresponding evolution equations on the \im\ are algebraically identical.  
That is, in these cases algebraic simplification is impossible.}
The crucial point throughout is that although the definition of amplitudes may differ, one uses linear combinations of the same physical sub-cell modes~\(v_0,\ldots,v_{M-1}\) as a basis to detail the microscale structures (the multiscale lifting).
\cite{Alavi2023} [p.2166] commented that ``proper elaboration of the macroscopic kinematic and static quantities that pertain to the micromorphic continuum is a problematic issue''.
In contrast, in our dynamical system framework there is no problematic issue: this paragraph indicates how the precise physical meaning of the variables used to parametrise a multi-modal, multi-continuum, homogenisation need be only mildly constrained by the physics of the problem and so is largely a \text{subjective aesthetic decision.} 

The discussion and theory so far is summarised in the following procedure to construct multi-continuum models.  
Variations of this procedure adequately encompass the more general cases of \cref{Shceq,Sgentheory,Selastic2d}.
This procedure may be implemented completely algebraically (e.g., this section), or mixed algebra-numerics (e.g., \cref{Shceq,Selastic2d}), as completely detailed with the aid of code listed in \ifJ Supplementary Material \cite[Appendices]{Roberts2024a}\else \cref{cas,cashc,cas2d}\fi.

\needspace{2\baselineskip}
\begin{procedure}[construct multi-continuum \im\ homogenised models]\label{Pconstruct}
Given a field~\(u(t,x)\) satisfies the dynamical system \(\alphaD u=f(u,x)\) with \(\ell\)-periodic microscale heterogeneity in space~\(x\).
\needspace{3\baselineskip}
\begin{enumerate}
\item Embed the given system in the ensemble of phase-shifts of the heterogeneity to form the embedded system \(\alphaD\fu = f(\fu,\theta)\), name it~(\ddag), where microscale heterogeneity is encoded via cell variable~\(\theta\) (e.g., \pde~\cref{Eemdifpde,Eempde}).  
Solutions~\fu\ are to be \(\ell\)-periodic in~\(\theta\).

\needspace{3\baselineskip}
\item Find a useful equilibrium (often \(\fu=0\), sometimes a family of equilibria), linearise~(\ddag) about the equilibrium, neglect macroscale gradients~\(\D x{},\gradx\), to obtain the basic sub-cell dynamics in a form such as \(\alphaD\fu = \cL_\theta\fu\) (e.g., \pde~\cref{Eemdifpde0,Eemeig}).

\needspace{2\baselineskip}
\item Solve the cell eigen-problem, \(\lambda v=\cL_\theta v\) with \(v\) is \(\ell\)-periodic in~\(\theta\), to obtain cell eigenvalues~\(\lambda_m\) and corresponding (generalised) eigenvectors~\(v_m(\theta)\).

\needspace{3\baselineskip}
\item Select \(M\) eigenvalues and eigenvectors of interest in the physical scenario, say numbered~\(\lambda_0,\ldots,\lambda_{M-1}\).   
Commonly, to construct a model to resolve dynamics on macroscale timescales longer than some time~\(T\), choose all modes for which \(|\lambda_m|^{1/\alpha}<1/T\).

\needspace{3\baselineskip}
\item In view of the subspace spanned by the corresponding cell eigenvectors~\(v_0,\ldots,v_{M-1}\), \emph{explicitly define} what you choose the multi-continuum amplitudes~\(U_0,\ldots,U_{M-1}\) to measure within cells (often averages or weighted averages of sub-cell properties).  
Normalise~\(v_0,\ldots,v_{M-1}\) to suit.

\needspace{3\baselineskip}
\item Set the first approximation of the corresponding \im~model for macroscale variables~\(U_m(t,x)\) according to a tangent space model, such as~\cref{EE1stApp1,Emode}: \(\fu\approx\sum_{m=0}^{M-1} U_mv_m(\theta)\) such that \(\alphaD U_m \approx \lambda_mU_m\)\,.

\needspace{2\baselineskip}
\item Select the orders of asymptotic truncation for the construction:  the order of macroscale spatial gradients;  the order of nonlinearity; and the order in any other `perturbative' physical effect(s).

\needspace{2\baselineskip}
\item Iterate to successively improve the \im~approximation.
\begin{enumerate}
\item Evaluate the residual \Res(\ddag) for the current \im~approximation.
\needspace{3\baselineskip}
\item Solve the appropriate homological equation (e.g.,~\cref{E1dHomologic,EEgenHomologic}) together with equations for amplitudes~\(U_m\), to make unique corrections to the \im~approximation (this updates the `slaved' modes).
\needspace{2\baselineskip}
\item Terminate the iteration when \Res(\ddag) is as small as desired---smallness is expressed in terms of order of macroscale spatial gradients~\(\D x{},\gradx\), and order of nonlinearity and/or other perturbation.
\end{enumerate}

\needspace{1\baselineskip}
\item Optionally regularise the resultant model \pde\ (\cref{SSSnrm,SSSisr}).
\end{enumerate}
\end{procedure}

\subsubsection{Accuracy compared to other homogenisation methods}

\im\ models are in principle made of trajectories of the original system (\cref{Dim}), and hence \im~models are in principle exact.  
This exactness holds throughout the domain of validity of the \im~model, namely the set~\(\UU(t)\) of \cref{Dim}. 

In practice we almost always construct approximations to a chosen \im\ and the evolution thereon---so in practice there is some error.   
For the broad range of systems to which they apply, \im~approximation theorems guarantees that the order of error in an approximation is the same as the order of error of the residual of the governing equations evaluated at a given \im~approximation, leading to the following \cref{CordAcc} that ensures a successful termination of \cref{Pconstruct} constructs an \im~model to the specified order of error.

\begin{corollary}[invariant manifold model error]\label{CordAcc}
Upon termination of the iterative loop in \cref{Pconstruct}, and provided suitable conditions on~(\ddag) hold for the theorems to apply, the following order-of-errors hold:
\begin{itemize}
\item for the case of an emergent \im\ (\cref{Dim}), the constructed manifold~\cM\ and the evolution thereon have errors~\Ord{\Res(\ddagger)};
\item for both emergent and guiding centre \im{}s (\cref{Dim}), the constructed manifold~\cM\ and the evolution thereon are exact for a system~\Ord{\Res(\ddagger)} close to the given system~(\ddag);
\end{itemize}
\end{corollary}

\begin{proof} 
In both cases the theory of \cite{Roberts2013a, Bunder2018a}, under suitable preconditions, transforms the analysis of the spatial structure across a spatial domain~\XX\ into a family of non-autonomously forced local dynamics at each station~\(X\in\XX\).
Then in both cases Thm.~2.18 of \cite{Hochs2019}, under suitable preconditions, guarantees the constructed~\cM\ and evolution thereon are exact for a system \Ord{\Res(\ddagger)}~close to the local dynamics of~(\ddag). 
In the emergent case, Prop.~3.6 of \cite{Potzsche2006}, under suitable preconditions, guarantees the local constructed~\cM\ and evolution thereon approximate an \im~of~(\ddag) to errors~\Ord{\Res(\ddagger)}.
\hfill\end{proof}

Consequently, if an \im~model and another model are compared in a regime where they are both valid, then they must agree to their order of error. 
If a discrepancy occurs between the two that is outside their order of error, then there must be some unaccounted for approximation made in one or the other.

\subsection{Spatial resolution of three-mode homogenisation}
\label{secsrh}

With the computer algebra of \ifJ Supplementary Code \cite[Appendix~A,][]{Roberts2024a} \else\cref{cas} \fi we easily construct multi-modal homogenisations, such as~\cref{EE3Mevol}, to any chosen order in gradients, and for a wide range of periodic heterogeneity.
In this linear class of problems we use high-order to quantitatively estimate limits of approximate homogenisations \text{such as~\cref{EE3Mevol}.}

\subsubsection{\protect\cref{EG1D}: convergence in heterogeneity~\(a\)}

Recall that in this example the heterogeneity~\cref{EGkappa}, \(\kappa=1/(1+a\cos \theta)\), has strength parametrised by~\(a\).
Let's first explore the dependence in~\(a\) of the tri-continuum model~\cref{EE3Mevol}.

\begin{figure}
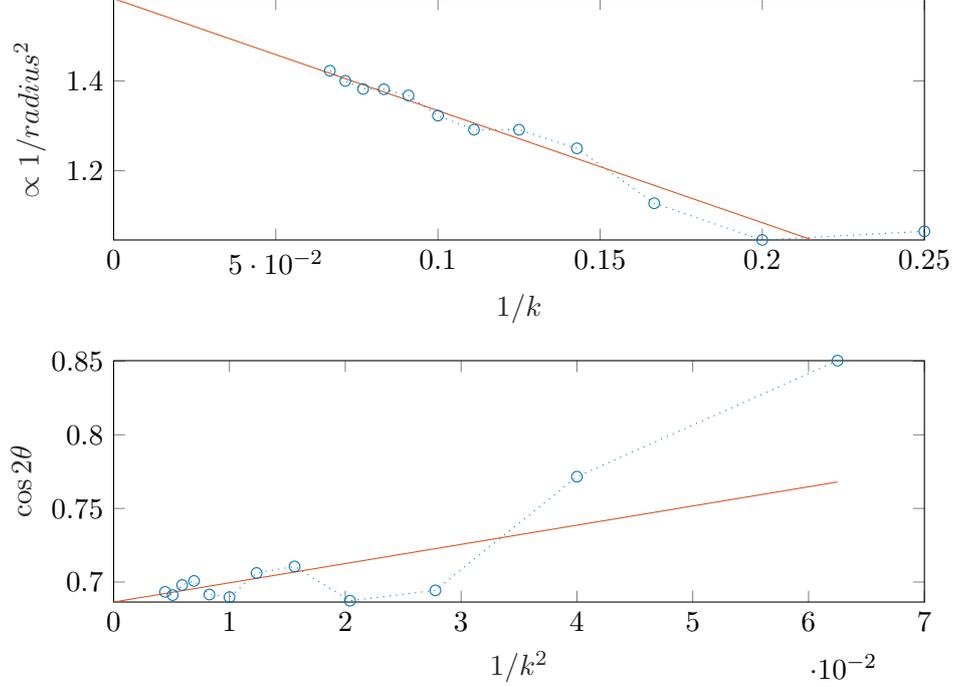
\centering
\caption{\label{Fcas1dCs2}Mercer--Roberts plot for the series in heterogeneity~\(a\) of the coefficient of~\(U_1\) in a high-order extension of~\cref{EE3Mevol1} for~\(\alphaD U_{1}\).  
The extrapolated intercepts to \(1/n=0\) predict the location of convergence limiting singularities in the complex \(a\)-plane.}
\inPlot{Figs/cas1dCs2}
\end{figure}

Let's construct the tri-mode homogenisation to low order in spatial gradient, errors~\Ord{\partial_x^3}, but here to high-order error in heterogeneity, namely~\Ord{a^{31}}.  
The computer algebra code takes less than three minutes to execute---the results are independent of \(\alpha\in\{1,2\}\).  
I chose three important coefficients in the extension of the homogenisation~\cref{EE3Mevol}:
\begin{itemize}
\item coefficient of~\(U_{0xx}\) in~\(\alphaD U_{0}\) that starts \(1+\tfrac12a^2+\tfrac5{24}a^4\cdots\);
\item  coefficient of~\(U_1\) in~\(\alphaD U_{1}\) that starts
\(-1-\tfrac5{12}a^2-\tfrac{437}{3456}a^4+\cdots\);
\item  coefficient of~\(U_2\) in~\(\alphaD U_{2}\) that starts
\(-1+\tfrac1{12}a^2-\tfrac{53}{3456}a^4+\cdots\).
\end{itemize}
These series are all of the form \(\sum_{n=0}^\infty c_nz^{n}\) for \(z=a^2\): 
for every four-tuple of consecutive computed coefficients,~\(c_{n-1}\C c_n\C c_{n+1}\C c_{n+2}\), \cite{Mercer90} [Appendix] give formulas to roughly estimate the location of the nearest complex-conjugate pair of singularities in the complex \(z\)-plane; 
a plot of these estimates versus~\(1/n\) is then used to extrapolate to the ideal estimate as \(n\to\infty\).
\cref{Fcas1dCs2} shows a Mercer--Roberts plot for the coefficients of~\(U_1\) in~\(\alphaD U_{1}\), the plot for the coefficients of~\(U_{0xx}\) in~\(\alphaD U_{0}\) is almost the same.  
For the third series, a simpler Domb--Sykes plot suffices \cite[e.g.,][]{Domb57, Hunter87}.
Such plots predict the radius of convergence limiting singularity in heterogeneity~\(a\) as \(1.21\C 1.21\C 1.57\), respectively, due to singularities in the complex
\(a\)-plane at respective angles~\(23^\circ\C 23^\circ\C 90^\circ\) to the real-\(a\) axis.
Remarkably, the predicted radius of convergence~\(1.21\) indicates that we may use the three-mode homogenisation model even up to (and~past!) the extreme contrast \text{heterogeneity of \(a\to1\).}

Practically, the radius of convergence indicates that via expansion to errors~\Ord{a^{11}} one would compute coefficients to four decimal places over the range \(|a|\lesssim1/2\).  
Further exploration indicates that the \([6,6]\) Pad\'e approximations in~\(a\) appear to be similarly accurate over the range \(|a|\leq 1\)\,.

\subsubsection{Convergence in spatial wavenumber}
\label{SSScsw}

Recall that traditional mathematical proofs of homogenisation require the scale separation limit that the length-scale ratio \(\ell/L\to0\).
In practice, engineers and scientists presume that \(\ell/L<0.1\) or~\(0.01\) is sufficient \cite[e.g.,][p.57]{Auriault2009}.  
For example, \cite{Somnic2022} comment~[p.4] ``For a periodic network of lattices to be considered as material, the characteristic length of its cells needs to be at least one or two orders of magnitude below the medium's overall length scale.''
In contrast, and surprisingly, numerical studies by \cite{Ameen2018} suggested~\(\ell/L\) could be as large as one.
By exploring the tri-modal homogenised evolution~\cref{EE3Mevol} to high-order in spatial gradients~\(\partial_x\), albeit  to low-order in heterogeneity~\(a\), we here quantify the range of valid scale ratios~\(\ell/L\), and reasonably agree with Ameen et al.   

Choosing errors~\Ord{a^4}, amazingly we find the evolution effectively truncates.
The computer algebra derives the \im\  homogenised evolution to~\Ord{\partial_x^{31} ,a^4} in just a few seconds. 
After a spatial Fourier transform  to wavenumber~\(k\), some Domb--Sykes plots in powers of wavenumber~\(k\) of the \(a^2,a^3\)~terms then shows that all are limited by simple pole singularities at \text{wavenumber \(k=\pm3/2\)\,.}

A first consequence is that these Domb--Sykes plots show that for small heterogeneity the tri-continuum modelling resolves all wavenumbers \(|k|<3/2\)\,.
This bound corresponds to all wavelengths bigger than~\(4\pi/3=2\ell/3\).  
That is, potentially the resolved macroscales are all wavelengths bigger than just \(2/3\)~of the microscale periodicity!  
However, a practical lower bound may be about twice this.
Moreover, be aware that higher orders in heterogeneity~\(a\) appear to be more restrictive \cite[]{Roberts2022a}.

A second consequence is that the tri-continuum, three-mode, homogenisation algebraically simplifies using the nonlocal operator
\(\cD:=(1+4/9\,\partial_x^2)^{-1}\).  
The \ifJ Supplementary Code \cite[Appendix~A,][]{Roberts2024a} \else computer algebra of \cref{cas} \fi finds the following to arbitrarily \text{high order in~\(\partial_x\):}
\begin{subequations}\label{eqimphpde}%
\begin{align}
\alphaD U_{0}&= (1+\tfrac12a^2)U_{0xx} -\tfrac12aU_{1x} -\tfrac12aU_{2xx}
-\tfrac1{72}a^3\cD(15\partial_x+11\partial^3+3\partial_x^5)U_1 
\nonumber\\&\quad{}
-\tfrac1{72}a^3\cD(11\partial_x^2+5\partial_x^4)U_2
+\Ord{a^4},
\\
\alphaD U_{1}&= -U_1-2U_{2x}+aU_{0x}+U_{1xx}
-\tfrac1{36}a^2\cD(15+18\partial_x^2+5\partial_x^4)U_1
\nonumber\\&\quad{}
-\tfrac1{18}a^2\cD(4\partial_x+4\partial_x^3+\partial_x^5)U_2
+\tfrac1{144}a^3\cD(60\partial_x+56\partial_x^3+13\partial_x^5)U_0
\notJbreak{}
+\Ord{a^4},
\\
\alphaD U_{2}&= -U_2+2U_{1x}-aU_{0xx}+U_{2xx}
+\tfrac1{36}a^2\cD(3+8\partial_x^2+3\partial_x^4)U_2
\nonumber\\&\quad{}
+\tfrac1{18}a^2\cD(4\partial_x+4\partial_x^3+\partial_x^5)U_1
-\tfrac1{72}a^3\cD(22\partial_x^2+12\partial_x^4+\partial_x^6)U_0
\notJbreak{}
+\Ord{a^4},
\end{align}
\end{subequations}
Since the operator~\cD\ is nonlocal, the model~\cref{eqimphpde} is an example of a nonlocal homogenisation \cite[e.g.,][]{Bazant2002}.

Effects of higher order than cubic in~\(a\) have yet to be explored in detail.
However, plots like \cref{Fcas1dCs2} for terms~\(a^4\) in heterogeneity indicate convergence limiting singularities at wavenumber \(k=\pm 1/2\).
That is, although for small~\(a\) we can get the `exact' nonlocal model~\cref{eqimphpde}, at larger heterogeneity~\(a\) singularities limit the homogenisation to macroscales longer than twice the length~\(\ell\) of the microscale period, as also found by \cite{Roberts2022a} for one-mode homogenisation.
Thus, more generally, the high-order homogenisation appears valid for macroscale \(L>2\ell\)\,.
Equivalently, the homogenisation is valid for scale ratios \(\ell/L<0.5\) (although a practical bound might be a half of this), which, in contrast, is significantly better than the ``one or two orders of magnitude'' usually assumed.

\subsection{\protect\cref{EG1D}: homogenise heterogeneous nonlinearity}
\label{SShhn}

Because the underlying theory of invariant manifold models is that of nonlinear dynamical systems theory, as discussed generally by \cref{Sgentheory},  homogenising nonlinear heterogeneous systems requires just a few straightforward modifications. 

Here consider constructing a one-mode homogenisation of the heterogeneous diffusion~\cref{Ehdifpde} with the addition of nonlinear heterogeneous advection, namely
\begin{equation}
\alphaD u = \D x{}\left\{\kappa(x)\D xu -\gamma\eta(x) u^2/2 \right\},
\quad 0<x<L\,.
\label{EhHetNon}
\end{equation}
For the \(\alpha=1\) case, this is akin to a Burgers' \pde\ with nonlinearity strength parametrised by~\(\gamma\), and heterogeneous coefficient in both the diffusion and nonlinear advection term, \(\kappa(x)\) and~\(\eta(x)\) respectively.
Such systems have previously been explored for propagating fronts \cite[e.g.,][]{JackXin2000}.
As before, let's non-dimensionalise on the microscale length~\(\ell\) so that the heterogeneities~\(\kappa(x)\) and~\(\eta(x)\) are \(2\pi\)-periodic: specifically
\begin{align}&
\kappa(\theta):=1/(1+a\cos \theta),
&&
\eta(\theta):=c_1\cos\theta+c_2\sin2\theta\,.
\label{egkappaeta}
\end{align}
Some exploration indicated this particular choice for~\(\eta(\theta)\) has interesting interactions with this~\(\kappa(\theta)\).

The corresponding phase-shift embedding modifies \pde~\cref{Eemdifpde} to the nonlinear
\begin{equation}
\alphaD \fu=\left(\D\x{}+\D \theta{}\right)
\left\{\kappa(\theta)\left(\D\x\fu+\D \theta\fu \right)
-\gamma\eta(\theta) \fu^2/2\right\},
\label{EemHetNon}
\end{equation}
for fields~\(\fu(t,\x,\theta)\) being \(2\pi\)-periodic in~\(\theta\).
As an example, we construct the one-mode \im\ homogenisation (\(M=1\)) of this embedding \pde\ (choose \(M>1\) for multi-continuum homogenisations).
The accessible class of nonlinear \im{}s is to construct homogenisations as a regular perturbation in nonlinearity parameter~\(\gamma\).  
There are just three necessary changes in the \ifJ Supplementary Code \cite[Appendix~A,][]{Roberts2024a}\else computer algebra of\cref{cas}\fi: 
firstly, truncating to some specified order of error in~\(\gamma\), here choose errors~\Ord{\gamma^2} to just report leading order effects of the nonlinearity;  
secondly, setting the extra heterogeneity~\(\eta(\theta)\) as~\cref{egkappaeta} specifies; and 
thirdly, modifying the computation of the residual by including the additional term~\(-\gamma\eta(\theta) \fu^2/2\) \text{in the flux.}

For the dissipative Burgers' case of \(\alpha=1\), executing the code constructs the one-mode homogenisation and finds that the \im\  ensemble field (to one order lower)\begin{subequations}\label{EEhetnonIM}%
\begin{align}
\fu&= U_0 +a\sin\theta\,U_{0x}
+\gamma c_1( \tfrac12\sin\theta +\tfrac18a\sin2\theta )U_0^2
\nonumber\\&\quad{}
-\gamma c_2( \tfrac14a\cos\theta +\tfrac14\cos2\theta +\tfrac1{12}a\cos3\theta )U_0^2
+\gamma c_1( \cos\theta -\tfrac18a\cos2\theta )U_0U_{0x}
\nonumber\\&\quad{}
+\gamma c_2( a\sin\theta +\tfrac14\sin2\theta -\tfrac19a\sin3\theta )U_0U_{0x}
+\Ord{\gamma^2,\partial_x^2,a^2}.
\label{EhetnonIM}
\end{align}
Such expressions for the \im\ field systematically account for all the sub-cell physics effects of the interaction of diffusion, nonlinearity, micro-structures and macro-gradients.
Simultaneously the code constructs that the evolution on the \im\ obeys the homogenised \pde
\begin{equation}
U_{0t}=U_{0xx}-\tfrac14\gamma\partial_x\big( c_1aU_0^2 
+c_2a^2U_0U_{0x} \big)
+\Ord{\gamma^2,\partial_x^3,a^3}.
\label{EhetnonIMe}
\end{equation}
\end{subequations}
Although the microscale nonlinear advection coefficient has zero-mean, \(\overline{\eta(\theta)}=0\), nonetheless the interaction of the two heterogeneities~\eqref{egkappaeta} generates a non-zero effective nonlinear advection\({}\approx-\frac12\gamma c_1aU_0U_{0x}\), as well as an effectively nonlinear macroscale  diffusion\({}\approx\partial_x\big[ (1-\tfrac14\gamma c_2a^2U_0)U_{0x}\big]\).

Wave systems, \(\alpha=2\), have more complicated models. 
But fortunately the models can be expressed as modifications to~\cref{EEhetnonIM}, modifications involving the time derivative \(V_0:=\D t{U_0}\).
The computer algebra of \ifJ Supplementary Code \cite[Appendix~A,][]{Roberts2024a}\else\cref{cas}\fi\ constructs that the \im\  ensemble field is
\begin{subequations}\label{EEhetnonIM2}%
\begin{align}
\fu&=(\text{right-hand side \cref{EhetnonIM}})
-\gamma c_1( \tfrac12\sin\theta +\tfrac5{16}a\sin2\theta )V_0^2
\nonumber\\&\quad{}
-\gamma c_2( \tfrac58a\cos\theta +\tfrac18\cos2\theta +\tfrac{13}{216}a\cos3\theta )V_0^2
-\gamma c_1( 6\cos\theta +\tfrac{21}{16}a\cos2\theta )V_0V_{0x}
\nonumber\\&\quad{}
-\gamma c_2( \tfrac92a\sin\theta +\tfrac38\sin2\theta +\tfrac19a\sin3\theta )V_0V_{0x}
\,.
\label{EhetnonIM2}
\end{align}
Simultaneously the code constructs that the evolution on the \im\ obeys the homogenised \pde
\begin{equation}
U_{0tt}=
(\text{right-hand side \cref{EhetnonIMe}})
+\tfrac12\gamma\partial_x\big( c_1aV_0^2 +3c_2a^2V_0V_{0x} \big).
\label{EhetnonIM2e}
\end{equation}
\end{subequations}
As in linear analysis, the main part of the model is the classic wave \pde, \(U_{0tt}=U_{0xx}\)\,.
However, the nonlinearity interacts with the heterogeneity to generate \Ord{a\gamma}~nonlinear effects on the macroscale.

In contrast, many other homogenisation methods do not apply to such nonlinear systems.


\section{Example: high-contrast multi-continuum homogenisation of a laminate}
\label{Shceq}

The modelling of materials with so-called \emph{high contrast} is of interest.%
\footcite[e.g.,][]{WingTatLeung2024, Efendiev2023, YifanChen2023}
This section extends the example of the multiscale \pde~\cref{Ehdifpde} to a \emph{high-contrast}, \(\ell\)-periodic, heterogeneous coefficient~\(\kappa(x)\) of a laminate in two spatial dimensions \cite[analogous to the example of stratified composites by][\S4]{Boutin1996}, namely
\begin{equation}
\Dn t\alpha u=\D x{}\left\{\kappa(x)\D xu\right\} +\kappa(x)\DD yu\,,
\quad (x,y)\in\Omega\,,
\label{Ehdifpde2}
\end{equation} 
where \(\Omega\) is a connected 2-D spatial domain that is of typical diameter~\(L\).
Specifically, in each microscale period, the coefficient \(\kappa(x)=\kappa_1\) constant for most~\(x\), except in a thin near-insulating layer of width~\(\eta\ll\ell\) where the coefficient \(\kappa(x) = \kappa_0\ll\kappa_1\).
We focus on modelling interesting dynamics in an interior~\(\XX\subset\Omega\), separated from~\(\partial\Omega\) by boundary layers of \Ord{\ell}~thickness, for domains~\(\Omega\) of relatively large size~\(L\).

This section shows how the novel and powerful invariant manifold~(\im) framework of \cref{SoneDintro}  establishes rigorous multi-mode multi-continuum homogenised models of the high-contrast laminate material.
For example, \cref{Shcegbmm} addresses the bi-continuum case and derives the following homogenised model in terms of two physics-informed macroscale quantities~\(U_0,U_1\) that evolve according to two coupled macroscale \pde{}s of the form, non-dimensionalised,
\begin{equation}
\alphaD U_{0}=0.81U_{0xx} +0.95U_{0yy} +0.36U_{1x} 
+\cdots,\quad
\alphaD U_{1}= -0.46U_1 -0.36U_{0x} 
+\cdots,\quad
\label{Epde2}
\end{equation}
where here the numerical coefficients are for a specific high-contrast case in the class defined above (\(\ell=2\pi\C \eta=0.377\C \kappa_1=1\C \kappa_0=0.06\)), and for both diffusion and waves, \(\alpha=1,2\)\,.
Physically, \(U_0\)~is the cell-mean of the field~\(u\), whereas the micromorphic variable~\(U_1\) measures the \(x\)-gradient in the field~\(u\) between consecutive insulating layers (see \cref{FhcegEvecs3}). 
The aim of this section is to justify and construct homogenisations such as~\cref{Epde2}.

Because our \im\ modelling is transitive (\cref{Dmulticont}),
the corresponding classic homogenised one-mode macroscale \pde\ for the laminate may be recovered by the adiabatic quasi-static approximation of the second mode in~\cref{Epde2} that gives \(U_1\approx -0.77U_{0x}\)\,.  
Thence the first \pde\ of~\cref{Epde2} reduces quantitatively to the usual homogenised macroscale \pde\ \(\alphaD U_{0}\approx 0.54U_{0xx} +0.95U_{0yy}\) for the cell-mean field of~\(u\) (which agrees with the directly constructed~\cref{EhcegManifold1U}).
All sound homogenisation methods would derive this particular model: one distinction here is that we place its derivation in a unified framework that extends straightforwardly to arbitrarily high-order gradients (\cref{Shcegsmm}), and also extends straightforwardly to arbitrarily many micromorphic modes (\cref{Shcegbmm}).
Indeed the number~\(M\) of modes, and the number~\(N\) of gradients, are coded as arbitrary parameters in the computer algebra of \ifJ Supplementary Code \cite[Appendix~B,][]{Roberts2024a}\else\cref{cashc}\fi.
One advantage of such multi-continuum models, such as the bi-continuum~\eqref{Epde2}, is that they usually resolve dynamics on shorter \text{identifiable space-time scales.}

Recall that we make rigorous progress through considering the ensemble of all phase-shifts.
For this laminate we only consider phase shifts in the direction~\(x\) of the heterogeneous variations. 
It is thus a straightforward variation of \cref{SoneDintro} to see that here we need to solve the embedding \pde
\begin{equation}
\Dn t\alpha\fu=\left(\D \x{}+\D \theta{}\right)
\left\{\kappa(\theta)\left(\D \x\fu+\D \theta\fu \right)\right\}
+\kappa(\theta)\DD y\fu\,,
\quad\fu\text{ \(\ell\)-periodic in }\theta,
\label{Eemdifpde2}
\end{equation}
in the `cylindrical' domain \(\dom:=\{(\x,y,\theta): (\x,y)\in\XX\C
\theta\in[-\ell/2,\ell/2]\}\).
A straightforward variation of \cref{secpse} would establish that solutions of such an embedding \pde~\cref{Eemdifpde2} provide us with solutions to the original heterogeneous \pde~\cref{Ehdifpde2}, for every phase-shift of the heterogeneity.

To establish the foundation of an \im~homogenisation, we consider the modified embedding \pde~\eqref{Eemdifpde2} with \(\D{\x}{}\) and~\(\D y{}\) neglected, and solved with \(\ell\)-periodic boundary conditions in~\(\theta\). 
This neglect gives the basic physical sub-cell problem explored in \cref{secsee}.
When the \pde\ is linear, as here, it is sufficient to consider the dynamics about the zero equilibrium, which we do henceforth.

\subsection{Spectrum at each equilibrium}
\label{secsee}

The approach is to choose invariant manifold models \citep{Roberts2013a} based upon the sub-cell physics encoded in the spectrum of the cell problem~\cref{Eemdifpde0}, and to choose options depending upon desired macroscale attributes. 

\begin{SCfigure}
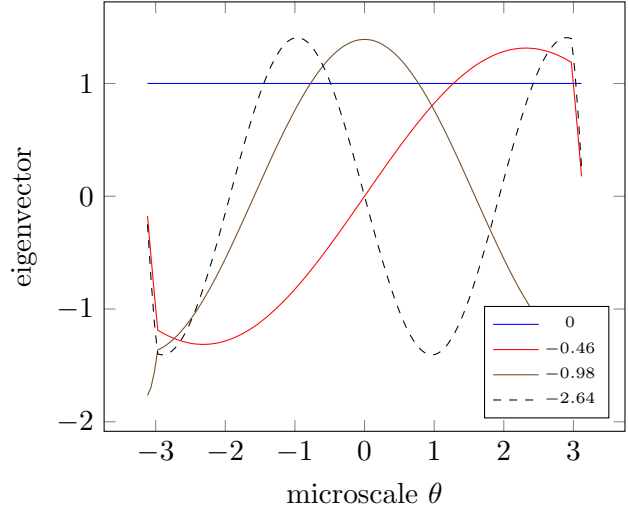
\centering
\caption{\label{FhcegEvecs3}an example of the leading four eigenvectors~\(v_m(\theta)\) of the cell-problem~\cref{Eemdifpde0}, and corresponding eigenvalues~\(\lambda_m\) in the legend, for a high-contrast heterogeneity.
For the non-dimensional case of \(\ell=2\pi\), \(\kappa_1=1\), and a thin insulating layer \(\eta/\ell=0.06\) located at \(\theta=\pm\pi\) of insulation, \(\kappa_0=0.06\) (\(\chi=1\)).}
\inPlot{Figs/hcegEvecs3}
\end{SCfigure}

The spectrum of the cell problem is in turn determined from the eigenvalues of the operator on the right-hand side of the cell-problem~\cref{Eemdifpde0} with its  boundary conditions of \(\ell\)-periodicity in~\(\theta\).
Say the eigenvalues are \(\lambda_0=0\geq\lambda_1\geq\lambda_2\geq\cdots\) with corresponding eigenvectors~\(v_0(\theta)\C v_1(\theta)\C v_2(\theta)\C\ldots\)\,, such as the example eigenvectors drawn in \cref{FhcegEvecs3}.
Generally we choose our modelling to focus on the sub-cell modes associated with the small magnitude eigenvalues because these are either the \emph{emergent} sub-cell modes, \(\alpha=1\), or the \emph{guiding centre} sub-cell modes, \(\alpha=2\) \cite[e.g.,][]{vanKampen85}.

For the specific example leading to the model~\cref{Epde2}, the non-dimensional spectrum of the cell eigen-problem~\cref{Eemdifeig} is \(\lambda_m\kappa_14\pi^2/\ell^2 \in \{0\C {-0.46}\C {-0.98}\C \ldots\}\).
In the wave case of \(\alpha=2\), these correspond to non-dimensional frequencies \( \sqrt{-\lambda_m\kappa_1}2\pi/\ell \in \{0\C {\pm0.68}\C {\pm0.99}\C \ldots\}\) (the zero frequency has multiplicity two).
The homogenised bi-continuum \pde~\cref{Epde2} is constructed by forming the approximate invariant manifold based upon the two smallest eigenvalues and corresponding eigenvectors (blue and red in the example of \cref{FhcegEvecs3}): that is, eigenvalues~\(\{0\C {-0.46}\}\) (or wave frequencies~\(\{0\C {\pm0.68}\}\)).
Since \(\lambda_2\approx-1\), recall from \cref{SSSmmmcme} that this two-mode choice would follow from a user needing the resultant model to resolve timescales longer than \(t_\ell\approx (\ell/2\pi\sqrt{\kappa_1})^{2/\alpha}\).

For these spatiotemporal systems a crucial issue is the practical spatiotemporal resolution of a model.
\cref{SSSsd,SSSwls,SSSfde} discuss the expected temporal resolution, whereas \cref{SSSisr} discusses general spatial resolution, and the links between the two.

\subsubsection{High-contrast thin layer}
\label{SSShctl}

To soundly homogenise the high contrast problem we need to determine the spectrum of the eigen-problem~\cref{Eemdifeig}.
We first analytically approximate the eigenvalue spectrum in cases when a layer of near `insulator' is so thin that we can replace it by a `jump' condition (as suggested by the example eigenvectors of \cref{FhcegEvecs3}).  
These analytic approximations guide subsequent numerical-algebraic construction and interpretation.
We deduces that eigenvalues are \(\lambda_m\approx -m^2\kappa_1\pi^2/\ell^2\) corresponding to eigenvectors, over the cell \(\theta\in(-\ell/2,+\ell/2)\), of \(v_m\approx \cos(m\pi \theta/\ell)\) for even~\(m=0,2,4,\ldots\)\,, and, with a `jump' across the layer at \(\theta=\pm\ell/2\), of \(v_m\approx \sin(m\pi \theta/\ell)\) for odd~\(m=1,3,5,\ldots\)\,.

We seek solutions to the eigen-problem~\cref{Eemdifeig} for the right-hand side operator. 
That is, we find \(\ell\)-periodic solutions to
\begin{equation*}
\kappa_1\DD\theta v =\lambda v\text{ on }(-\ell/2,\ell/2),\quad
\text{except in a layer where }\kappa_0\DD\theta v =\lambda v\,,
\end{equation*}
for an `insulating' layer of small thickness~\(\eta\) and where \(\kappa_0=\Ord{\eta}\).
We know all eigenvalues \(\lambda\leq0\).

\paragraph{Thin insulating layer}
Here derive jump conditions across the thin layer.
For algebraic simplicity, \emph{temporarily} set the origin of~\(\theta\) at the centre of the thin layer so the layer is the interval~\((-\eta/2,+\eta/2)\).

Within the thin layer an eigenvector
 is of the form
\def\ko{\sqrt{-\lambda/\kappa_0}}
\(v=C\cos\big(\ko \theta\big)+D\sin\big(\ko \theta\big)\).
Define
\begin{align*}&
[v]:=v_{\eta/2}-v_{-\eta/2}
=2D\sin\big(\ko\eta/2\big)
=D\ko\eta+\Ord{\eta^{3/2}};
\\&
\overline v:=\tfrac12(v_{\eta/2}+v_{-\eta/2})
=C\cos\big(\ko \eta/2\big)=C+\Ord{\eta}.
\end{align*}
Hence \(C\approx\overline v =\Ord{1}\) and \(D\approx [v]/\big(\ko \eta\big) =\Ord{\eta^{-1/2}}\).
Then within the layer the derivative \begin{align*}
v_\theta &= -\ko C\sin\big(\ko \theta\big) +\ko D\cos\big(\ko \theta\big)
\\&\approx -\ko \overline v \sin\big(\ko \theta\big) +([v]/\eta)\cos\big(\ko \theta\big).
\end{align*}
So the jump and the mean of the derivative are
\begin{align*}&
[v_\theta]\approx-2\ko \overline v \sin\big(\ko \eta/2\big)
\approx (\lambda\eta/\kappa_0)\overline v\,,
\\&
\overline{v_\theta}\approx([v]/\eta)\cos\big(\ko \eta/2\big)
\approx [v]/\eta\,.
\end{align*}

\paragraph{Outside the layer}
The eigenvectors~\(v(\theta)\) are to be continuous so the jump~\([v]\) and mean~\(\overline v\) are the same inside and outside the layer.  
And the flux has to be continuous across the layer boundary, that is \(\kappa_1v_\theta^\text{outside}=\kappa_0v_\theta^\text{inside}\).
Hence outside the layer we have the conditions
\begin{subequations}\label{eqs:}%
\begin{align}&
[v_\theta]=\big[(\kappa_0/\kappa_1)v_\theta^\text{inside}\big]
\approx (\lambda\eta/\kappa_1)\overline v\,,
\label{Eqhcjumpy}
\\&
\overline{v_\theta} =\overline{(\kappa_0/\kappa_1)v_\theta^\text{inside}}
\approx \kappa_0/(\kappa_1\eta)[v] =1/(\chi\ell)[v] \,,
\label{Eqhcmeany}
\end{align}
\end{subequations}
as we choose to scale~\(\kappa_0\) so that \(\kappa_0/\kappa_1=\eta/(\chi\ell)\) for some insulating parameter~\(\chi\).
That is, the layer diffusivity\slash elasticity \(\kappa_0\propto\eta\kappa_1\) decreases with the relative layer thickness~\(\eta/\ell\).  
Thus small~\(\eta\) characterises a high-contrast material.
Parameter~\(\chi\) characterises the strength of the `insulation' in the thin layer: larger is more insulating, whereas smaller is less so.
 
For algebraic simplicity we now reset the origin of~\(\theta\) so that the thin layer is at \(\theta=\pm\ell/2\), and hence a jump across the thin layer is hereafter \([v]=v_{-\ell/2}-v_{+\ell/2}\).

There are two families of eigenvectors and eigenvalues.
\begin{itemize}
\item The symmetric family is eigenvectors 
\(v=\cos k\theta\) for eigenvalue \(\lambda=-\kappa_1k^2\) for some wavenumbers~\(k\) to be determined. 
For this eigenvector
 \begin{align*}&
\overline v=\cos(k\ell/2),
&& [v]=0,
&& \overline{v_\theta}=0,
&& [v_\theta]=2k\sin(k\ell/2).\lambda
\end{align*}
Hence \cref{Eqhcmeany} is satisfied, whereas \cref{Eqhcjumpy} requires that \(2k\sin(k\ell/2)=-k^2\eta\cos(k\ell/2)\), that is, \(2k\tan(k\ell/2)=-k^2\eta\to 0\) as \(\eta\to0\)\,.
Hence these eigenvectors occur for wavenumber \(k=m\pi/\ell\) for even integer~\(m\).
That is, \(v_m=\cos(m\pi \theta/\ell)\) and corresponding eigenvalues \(\lambda_m=-\kappa_1\pi^2m^2/\ell^2\) for \(m=0,2,4,\ldots\).
\cref{Tspec} list the first two of these eigenvalues, \(\lambda_0,\lambda_2\), for four selected parameters~\(\chi\) of thin insulation layer width~\(\eta\).
\begin{SCtable}
\caption{\label{Tspec}first two \(\fK\)-values that solve \(\tan\fK=-\chi\fK\) for four values of insulation strength~\(\chi\).
Below are the leading four eigenvalues for the non-dimensional case of \(\kappa_1=1\) and \(\ell=2\pi\)\,: these approximate the eigenvalues for small insulation layer width~\(\eta\). }
\(\ \begin{array}{l|llll}
\chi&1/3&1&3&9\\
\hline
\fK_1&2.46&2.03&1.74&1.63\\
\fK_3&5.23&4.91&4.78&4.74\\
\hline
\lambda_0&0&0&0&0\\
\lambda_1&-0.61&-0.42&-0.31&-0.27\\
\lambda_2&-1&-1&-1&-1\\
\lambda_3&-2.77&-2.45&-2.32&-2.27\\
\hline
\end{array}\)
\end{SCtable}

\item The asymmetric family is eigenvectors of the form \(v=\sin k\theta\) for eigenvalue \(\lambda=-\kappa_1k^2\) for some wavenumbers~\(k\) to be determined.
For this eigenvector
 \begin{align*}&
\overline v=0,
&& [v]=-2\sin(k\ell/2),
&& \overline{v_\theta}=k\cos(k\ell/2),
&& [v_\theta]=0.
\end{align*}
Hence \cref{Eqhcjumpy} is satisfied, whereas \cref{Eqhcmeany} requires that \(k\cos(k\ell/2)=-2/(\chi\ell)\sin(k\ell/2)\), 
\begin{equation}
\text{that is, }\tan(k\ell/2)=-\chi(k\ell/2).
\label{Eqhcaspec}
\end{equation}
For \emph{odd} integer~\(m\), let \(\fK_m\) be the solutions of \(\tan\fK=-\chi\fK\) in sequence so that \(m\pi/2<\fK_m\leq (m+1)\pi/2\).  
Then wavenumber \(k=2\fK_m/\ell\) satisfies~\cref{Eqhcaspec} and so asymmetric eigenvectors are \(v_m=\sin(2\fK_m\theta/\ell)\) corresponding to eigenvalues \(\lambda_m=-\kappa_14\fK_m^2/\ell^2\) for \(m=1,3,\ldots\)\,.%
\footnote{Using \(\tan\fK\approx 1/(\pi/2-\fK)\), gives \(1/(\pi/2-\fK_1)\approx -\chi\fK_1\) which leads to \(\fK_1 \approx \tfrac\pi2+\tfrac2{\pi\chi}\).
Similarly, \(\fK_3\approx \tfrac{3\pi}2+\tfrac2{3\pi\chi}\).
These reproduce \cref{Tspec} within errors~\(0.005\)--\(0.2\) over \(\chi\geq1\)\,.}
\cref{Tspec} lists the first two of these eigenvalues, \(\lambda_1,\lambda_3\), for four selected parameters~\(\chi\) of thin insulation \text{layer width~\(\eta\). }
\end{itemize}

\subsection{One-mode slow manifold homogenisation}
\label{Shcegsmm}

One may construct a \emph{slow}~\im\ homogenised model for the laminate based upon the eigenvalue zero, here corresponding to the one sub-cell mode~\( v_0=\cos 0\theta=1\).
Since~\(v_0\) is constant, such slow~\im\ modelling gives the classic homogenised \pde, but generalised to higher-order gradients at finite scale separation \cite[e.g.,][]{Roberts2022a}.
In the diffusion case, \(\alpha=1\), the argument for its emergence as a valid model is that all other sub-cell modes decay exponentially quickly in time, the slowest of which is~\(\exp({-\kappa_14\fK_1^2t/\ell^2})\).

Let's explore the non-dimensional case of \(\kappa_1=1\) and \(\ell=2\pi\) (\cref{Tspec}) with the specific insulating thin layer \(\eta/\ell=\kappa_0=0.06\) (i.e., \(\chi=1\)).
All the cell-problems are solved numerical on a sub-cell grid with 128~points per cell.
The numerically obtained leading non-zero eigenvalue is \(\lambda_1=-0.4597\), so in the diffusion case the decay to the slow~\im\  homogenisation is roughly like~\(\e^{-0.46\,t}\), from any given initial condition.
The computer algebra of \ifJ Supplementary Code \cite[Appendix~B,][]{Roberts2024a}\else\cref{cashc}\fi\ uses \cref{Pconstruct}, with modes \(M=1\) and via the residual of the embedding \pde~\cref{Eemdifpde2} (discretised in~\(\theta\)), to construct a slow~\im\ to any specified order in~\(\grad{}=(\partial_x,\partial_y)\).
The result is that the detailed slow~\im\  field
\begin{subequations}\label{EEhcegManifold1}%
\begin{equation}
\fu(t,x,\theta)=U_0 +u_{10}(\theta)U_{0x} +u_{20}(\theta)U_{0xx} +u_{02}(\theta)U_{0yy} +\Ord{\grad^3},
\label{EhcegManifold1u}
\end{equation}
in terms of the coefficient functions plotted in \cref{FhcegManifold1}.
\begin{SCfigure}
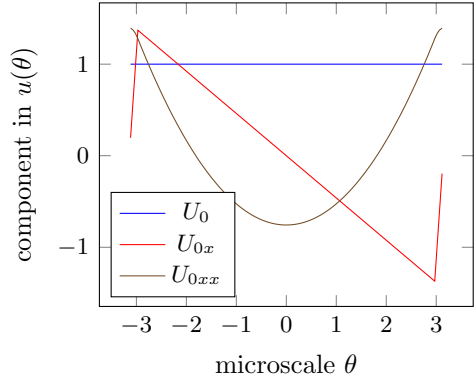
\centering
\caption{\label{FhcegManifold1}the sub-cell structure of the one-mode slow~\im\  field~\cref{EhcegManifold1u} in the high-contrast thin-layer problem.  
Specifically, this is the non-dimensional case of \(\ell=2\pi\), \(\kappa_1=1\), and a thin insulating layer \(\eta/\ell=0.06\) located at \(\theta=\pm\pi\) of insulation parameter \(\chi=1\) so \(\kappa_0=0.06\). }
\inPlot{Figs/hceg2Manifold1}
\end{SCfigure}%
For this high-contrast thin layer, the \text{\(U_{0x}\)-component} of \cref{FhcegManifold1} shows that \(x\)-gradients of the field lead to a sub-cell field where rapid spatial variation takes place in the thin insulating layer, as expected physically.
Corresponding to~\cref{EhcegManifold1u}, but to higher-order gradients, the homogenised evolution is
\begin{align}
\alphaD{U_0}&= .5386\,U_{0xx}  +.9486\,U_{0yy}
+(.3379\,U_{0xxxx} -.1381\,U_{0xxyy} +.0145\,U_{0yyyy})
\nonumber\\&{}
+(-.9019\,U_{0xxxxxx} +.2824\,U_{0xxxxyy} -.0021\,U_{0xxyyyy} 
-.0038\,U_{0yyyyyy})
\notJbreak{}
+\Ord{\grad^7}.
\label{EhcegManifold1U}
\end{align}
The leading-order homogenised \pde, \begin{equation}
\Dn t\alpha{U_0} = \divv(K\grad U_0) +\Ord{\grad^3},
\quad\text{with }K:=\operatorname{diag}(0.5386,0.9486),
\label{Ehceg1classic}
\end{equation} is the classic homogenisation with the anisotropic diffusion\slash elasticity tensor~\(K\) (as also deduced quasi-adiabatically following~\cref{Epde2}).
In such approximations, the error term denoted by~\(\Ord{\grad^p}\) represents a complicated remainder expression (\cref{Dao}).

Classic homogenisation approaches would construct parts of the model~\cref{EhcegManifold1U} because both classic methods and the approach here solve cognate cell problems obtained faithfully from the given physical equations.
For example, the formulas by \cite{Boutin1996} for fourth-order homogenisation of a static laminate would here reproduce the right-hand side of the first line of~\cref{EhcegManifold1U}.
However, in contrast, our dynamical systems framework, empowered by computer algebra, supports arbitrarily high-order construction (e.g., \cref{secsrh}), and explicitly encompasses generalisations to dynamics, to non-autonomous systems, to nonlinearity, and to multi-continuum `micromorphic' models \text{(e.g., \cref{Shcegbmm,SSpemn}).}

\subsubsection{Nonlocality regularises models}
\label{SSSnrm}

Higher-gradient models such as~\cref{EhcegManifold1U} often need regularisation \cite[often essential in practical wave modelling:][]{Benjamin1972, Bona76}. 
For example, upon neglecting the sixth-order derivatives, the homogenisation~\cref{EhcegManifold1U} may be regularised to \begin{equation}
(1-0.63\partial_x^2-0.02\partial_y^2)\,\alphaD U_0= 0.54\,U_{0xx}+0.95\,U_{0yy} -0.74U_{0xxyy} +\Ord{\grad^5}.
\label{EhcegManifold1Uc}
\end{equation} (to two decimal places). 
Whereas upon retaining the sixth-order derivatives, the homogenisation may be regularised to
\begin{align}&
(1 -0.63\partial_x^2 -0.02\partial_y^2 +2.07\partial_x^4 +0.007\partial_x^2\partial_y^2  +0.004\partial_y^4)\,\alphaD U_0
\nonumber\\&{}
= 0.54\,U_{0xx} +0.95\,U_{0yy} -0.74U_{0xxyy} +2.33U_{0xxxxyy}
 +\Ord{\grad^7}.
\label{EhcegManifold1Ud}
\end{align}
\end{subequations}
An equivalent form of these two regularised \pde{}s shows they are effectively nonlocal in space.
For example, we may rewrite~\cref{EhcegManifold1Uc} as the explicit nonlocal homogenisation \(\alphaD U_0\approx \cK\star (0.54\,U_{0xx}+0.95\,U_{0yy} -0.74U_{0xxyy})\) in terms of the convolution kernel \(\cK(x,y) \propto K_0\big( \sqrt{x^2/0.63+y^2/0.02} \big)\) and the modified Bessel function~\(K_0\).
This kernel decays to zero on the heterogeneity length~\(\ell=2\pi\) since \(K_0(r)\sim\sqrt{2/(\pi r)}\e^{-r}\) as \(r\to+\infty\) (similarly for higher-order regularised homogenisations).  
\cite{Bazant2002} discussed how such nonlocal models may improve spatial resolution to desirably capture smaller-scale effects, produce convergent numerical solutions, and capture size effects seen in experiments.
\cref{SSSisr} links such improvements to the characteristics of general \text{perturbed eigenvalue problems.}

But what these high-gradient regularised homogenisations \emph{usually} cannot do is to make a significant difference to the time resolution of a model---the reason is that they all neglect, except for the first `mean' mode, the independent dynamics of all the sub-cell modes.
In principle, one could involve extra time derivatives in regularising model \pde{}s.
But such higher-order in time-derivatives are equivalent to a first-order system in more variables, so one might as well instead construct a properly physics-informed micromorphic model such as \text{that of \cref{Shcegbmm}.}

\subsubsection{Boundary conditions for macroscale models}
\label{Sbcmm}

Traditional homogenisation theory provides a leading order macroscale model, such as~\cref{Ehceg1classic}, with boundary conditions at \((x,y)\in\partial\Omega\), often Dirichlet conditions on the macroscale field~\(U_0\).
Such macroscale boundary conditions are rigorous \emph{only} in the extreme scale separation limit \(\ell/L\to0\)\,.
At the finite~\(\ell/L\) of real applications, and especially for higher gradient models such as~\cref{EhcegManifold1U},  corrections are needed.
Further, new arguments are needed for multi-continua models \text{such as~\cref{Epde2}.}

\cite{Boutin1996} [\S3.2]  discussed such boundary conditions and observed that ``boundary layers have to be introduced in order to match real and homogenized conditions. 
The analysis of these layers, up to third-order, requires specific developments \ldots\ [the] layers have a `thickness' of about one period''.
\cite{Auriault2009} [pp.59,72] likewise commented that ``the introduction of matching boundary layers then makes it possible to complete the solution to the problem.''
Recently, \cite{Fergoug2022a} asserted [p.3] that the ``construction of a solution near the vicinity of the boundaries remains beyond capabilities of the classical homogenization''.
They went on to propose [p.4] ``a general boundary layer correction methodology for asymptotic homogenization in order to approximate real microscale fields near the boundaries.''
These methods envisage constructing such boundary-layer structures explicitly.

Alternative suitable basic theory and approach had already been developed for general \pde\ systems in 1D~space \cite[]{Roberts92c}, again using nonlinear dynamical systems theory and techniques \cite[see also][]{Chen2016}.
In contrast to the direct methods, this approach solves a dual problem to efficiently and accurately project through the boundary layers without explicitly constructing physical \text{boundary-layer structures.}

However, here we restrict attention to modellinginside a (large) domain~\XX\ separated from the physical boundaries at~\(\partial\Omega\), as indicated in \cref{Fpattensemble},  leaving boundary layer analysis to \text{further research.}

In the case of nonlocal models, such as~\cref{EhcegManifold1Uc}, \cite{Maugin2010a} [p.9] asked ``What about boundary conditions that are in essence foreign to this representation of matter-matter interaction?''
An answer here is that the nonlocality is on a length scale that is surely comparable to or smaller than the necessary boundary layers excluded from~\XX.  
For example, the convolution kernel~\cK\ for~\cref{EhcegManifold1Uc} decays to a negligible~\(10^{-3}\) over the microscale cell-width in that laminate.
I contend that for points in~\XX, the nonlocal convolutions generally do not reach the physical boundaries to any practical extent, and so do not need to be addressed herein.
Nonetheless, in principle, a nonlocal kernel like~\cK\ should itself be provided with boundary conditions that would change the kernel near the boundaries to a form that respects the boundaries and their associated boundary layers.   
\text{This case also needs research.}

\subsection{Two-mode, bi-continuum, homogenisations exist and emerge}
\label{Shcegbmm}

A two-mode \im~homogenised model for the laminate, such as~\cref{Epde2}, may be constructed based upon the leading two eigenvalues, generally resulting in improved space-time resolution.
For definiteness we continue to non-dimensionalise space-time so that cell-length \(\ell=2\pi\) and the coefficient \(\kappa_1=1\), and also focus on the case of thin layer parameter \(\chi=1\), the second \text{column of \cref{Tspec}.}

In this case the leading two eigenvalues are \(\lambda_0=0\) and \(\lambda_1\approx-0.42\) corresponding to the two
sub-cell modes \(v_0=1\) and \(v_1\approx\sin(2.03\,\theta/\pi)\) (the two blue curves in the two plots of \cref{FhcegManifold2} are more precise). 
In the diffusion case, \(\alpha=1\), and since the next eigenvalue \(\lambda_2=-1\), such a two-mode model is emergent with the slowest transient decaying roughly like~\(\e^{-t}\).   
The eigenvalue gap of~\((-0.42,-1)\) caters for
perturbing macroscale \(x,y\)-gradients.  

\cref{Pconstruct} constructs the \im\ homogenisation corresponding to these two modes.
The iterative construction starts from the initial approximation (\cref{C1stApp1}) that \begin{equation}
\fu\approx U_0 +U_1v_1(\theta)
\quad\text{such that}
\quad\alphaD U_0\approx0\,,
\quad \alphaD U_1\approx -0.42\,U_1\,.
\label{EhcegInit}
\end{equation}
The computer algebra of \ifJ Supplementary Code \cite[Appendix~B,][]{Roberts2024a}\else\cref{cashc}\fi\ then iteratively corrects approximations until the governing embedding \pde~\cref{Eemdifpde2} has residual smaller than a chosen specified order of error. 

In this laminate example we modify the procedure in two ways.
Firstly, recall that \cref{EG1D} implemented the procedure by doing all steps exactly in algebra.   
But for this high-contrast media we do not have exact algebraic expressions for the eigenvectors (\cref{SSShctl}).
Consequently, we adopt a simple sub-cell \(\theta\)-space discretisation of the cell eigen-problem~\cref{Eemdifeig} and the corresponding homological equation~\cref{E1dHomologic}.
The computer algebra sets \(n=128\) points per period of sub-cell variable~\(\theta\), and uses centred differences in~\(\theta\), which should be fine enough to be faithful to the microscale differentials to about four significant digits.
The macroscale variations in~\(x,\x\) and~\(y\) are still represented explicitly in algebra---there is no numerical approximation on the macroscale.
\cite{MacKenzie05} first discussed such fine-grid numerics for constructing invariant manifolds of the macroscale dynamics of the 1-D and 2-D Kuramoto--Sivashinky \pde\ and the need for numerics in 2-D \citep[see also ][]{Roberts2011a}.
\cref{Selastic2d} also uses such mixed numerics-algebra for multi-continuum homogenisation \text{of 2-D heterogeneous elasticity.}

The second modification is that for generalised multi-continua homogenisation it is often awkward to code solutions to the homological equation~\cref{E1dHomologic} when it involves modes with non-zero eigenvalue\slash frequency.  
Instead we may more quickly code, and do so here, simpler updates which just take more computer iterations to be accurate.
The simplification is to omit the tricky \(\sum_{m=0}^{M-1}\lambda_m(\D{U_m}{v'})U_m\) on the left-hand side of the homological equation~\cref{E1dHomologic}.
Then the left-hand side operator is a straightforward constant matrix which is efficiently inverted or LU-factored just once \ifJ\cite[Appendix~B.3,][]{Roberts2024a}\else(\cref{SSfifhu})\fi\ and used repeatedly.  
The updates~\(v'\) and~\(G'\) are then not precise, but they are accurate enough to make systematic progress.
The numerical error in the various coefficients of~\(v\) and~\(G\) here decreases each iteration by the ratio of the largest magnitude eigenvalue in the model to the smallest magnitude eigenvalue neglected by the model, that is, the ratio~\(|\lambda_{M-1}|/|\lambda_M|\).  
That is, the numerical convergence is quicker for a larger spectral gap between~\(\lambda_{M-1}\) and~\(\lambda_M\) (\cref{Dgap}).
We simply let the computer do more iterations until the numerical error is small enough: \ifJ the Supplementary Code \cite[Appendix~B,][]{Roberts2024a}\else\cref{cashc}\fi\ sets a maximum \text{relative error of~\(10^{-8}\).}

We now explore the bi-continuum homogenisation up to fourth-order in macroscale gradients.
Upon executing the code, 25~iterations are sufficient to give the following detailed physics-informed sub-cell field, to relative errors\({}\approx10^{-8}\), and in terms of the eight coefficient functions plotted in \cref{FhcegManifold2}:
\begin{subequations}\label{EEhcegManifold2}\needspace{2\baselineskip}%
\begin{align}
\fu(t,x,\theta)& =\phantom{v_1(\theta)}U_0 +u_{01}(\theta)U_{0x} +u_{02}(\theta)U_{0xx}
 +u_{03}(\theta)U_{0yy}
\nonumber\\&{}
+v_1(\theta)U_1+u_{11}(\theta)U_{1x}+u_{12}(\theta)U_{1xx}
 +u_{13}(\theta)U_{1yy}
+\Ord{\grad^3}.
\label{EhcegManifold2u}
\end{align}
\begin{figure}\centering
\caption{\label{FhcegManifold2}the sub-cell structure of the bi-continuum, \im\ field~\cref{EhcegManifold2u} in the high-contrast laminate problem.  
Specifically, the non-dimensional case of \(\ell=2\pi\), \(\kappa_1=1\), and a thin insulating layer \(\eta/\ell=0.06\) located at \(\theta=\pm\pi\) of insulation parameter \(\chi=1\) so \(\kappa_0=0.06\). }
\begin{tabular}{@{}c@{}c@{}}
\inPlot{Figs/hceg2Manifold20}
&
\inPlot{Figs/hceg2Manifold21}
\end{tabular}
\end{figure}

The corresponding homogenised evolution, but to higher order in gradients, for the macroscale variables~\(U_0,U_1\) is constructed to be, for both \(\alpha=1,2\),
\begin{align}
\alphaD{U_0}&=\phantom{- 0.4597\,U_{1}} +0.3552\,U_{1x} +0.8130\,U_{0xx} 
+0.9486\,U_{0yy} +0.0087\,U_{1yy}
\nonumber\\&\quad{} 
+1.063\,U_{1xxx} +0.2159\,U_{1xyy}
 - 0.0026\,U_{1xxyy} - 0.0024\,U_{1yyyy}
 \nonumber\\&\quad{}
 + 0.6146\,U_{0xxxx} + 0.1323\,U_{0xxyy} + 0.0143\,U_{0yyyy}
 +\Ord{\grad^5},
\label{EhcegManifold2U0}
\\
\alphaD{U_1}&= - 0.4597\,U_{1} - 0.3552\,U_{0x} - 2.620\,U_{1xx} 
+0.9742\,U_{1yy} +0.0087\,U_{0yy} 
\nonumber\\&\quad{} 
- 1.736\,U_{0xxx} - 0.1143\,U_{0xyy}
+ 0.0284\,U_{0xxyy} - 0.0014\,U_{0yyyy}\nonumber\\&\quad{}
- 17.11\,U_{1xxxx} + 0.7877\,U_{1xxyy}+ 0.0018\,U_{1yyyy}
 +\Ord{\grad^5}
\label{EhcegManifold2U1}
\end{align}
\end{subequations}
Various truncations and regularisations of these two \pde{}s form bi-continuum, homogenised models for this high-contrast laminate.
As discussed in \cref{SoneDintro,Sgentheory} this homogenisation is supported by extant rigorous dynamical systems theory.
In application, one truncates the \pde{}s~\cref{EEhcegManifold2} to an order of error suitable for the purposes at hand (and possibly with some suitable regularisation).
In solutions obtained using a truncated~\cref{EEhcegManifold2}, one could quantitatively estimate the modelling error via the remainder expression~(52) of \cite{Roberts2016a}, or more approximately by evaluating the neglected term(s) of the \text{next-higher order.}


\subsection{Potentially extend to more modes and to nonlinearity}
\label{SSpemn}

The computer algebra of \ifJ Supplementary Code \cite[Appendix~B,][]{Roberts2024a}\else\cref{cashc}\fi\ constructs corresponding \(M\)-mode invariant manifold homogenisations for any specified number~\(M\)\,.
For example, a three-mode, tri-continuum, invariant manifold, homogenisation of the laminate may be constructed based upon the leading three eigenvalues: for example, \(\lambda=0,-0.42,-1\) in the non-dimensional case \(\chi=1\), \(\ell=2\pi\), and \(\kappa_1=1\) (\cref{Tspec}).
In the diffusive case, \(\alpha=1\), such a three-mode
homogenisation emerges with the slowest
transient decaying roughly like~\(\e^{-2.5\,t}\) (\cref{Tspec}). 
That is, such a three-mode homogenisation is valid over shorter times than the bi-continuum, homogenisation~\cref{EEhcegManifold2}, which in turn is valid over shorter times than the one-mode classic homogenisation~\cref{EhcegManifold1U}.
The eigenvalue gap~\((-1,-2.5)\) caters for the perturbing macroscale \(x,y\)-gradients.
Similarly for any other chosen~\(M\).

This methodology also readily adapts to homogenising \emph{nonlinear} heterogeneous systems, as in \cref{SShhn}, with theoretical support established in the next \cref{Sgentheory}. 
\ifJ The Supplementary Code \cite[Appendix~B,][]{Roberts2024a}\else\cref{cashc}\fi\ 
includes an optional nonlinear term as an example.
With nonlinearity a multi-mode construction generally takes many more iterations.
Theory for constructing nonlinear invariant manifolds \cite[e.g.,][]{Potzsche2006} indicates the bound that the spectral gap ratio (\cref{Dgap}) \emph{must} be larger than the order of the nonlinearity.
For the example laminate here, the spectral gap ratio (\cref{Tspec}) for bi-continuum or tri-continuum models is only about~\(2.4\) so a quadratically nonlinear model is straightforwardly constructed, but a cubic model is problematic.
In contrast a one-mode classic homogenisation has an infinite spectral gap ratio and so models may be constructed for arbitrarily high-orders of nonlinearity.  
\emph{The modelling or homogenisation of nonlinear systems generally requires a bigger spectral gap than that needed \text{for linear systems.}}


\section{General multi-continuum, multi-mode, homogenisation of heterogeneity}
\label{Sgentheory}

Generalising the previous \cref{SoneDintro,Shceq}, this section develops this innovative approach to the rigorous multi-continuum, multi-mode, homogenisation of the dynamics of \emph{nonlinear, non-autonomous, multi-physics problems in multiple large space dimensions with periodic or quasi-periodic heterogeneity}.
Encompassing multi-physics scenarios, as is done here, is an outstanding challenge according to the review by \cite{Fronk2023} [\S1.4.3].
Our approach does \emph{not} invoke any variational principle and so applies to a much wider variety of systems than many homogenisation methods.
Instead, this general approach is supported by the rigorous dynamical system framework and theory of invariant manifolds.\footcite[e.g.,][]{Carr81, Muncaster83b, Bates98, Aulbach2000, Prizzi02, Haragus2011, Roberts2014a, Chekroun2015a, Hochs2019}

Consider quite general multiscale materials with complicated microstructure.  
Suppose that the spatial domain has \(d\)~dimensions of large extent, the macroscale, and possibly some thin spatial dimensions: examples include elastic beams and plates, or thin fluid films and shallow water, but also include in scope extensive 3-D materials with no apparent thin physical dimension.  
Consider times in a physically relevant interval~\(\TT\subseteq\RR\)\,. 
Let position in the large dimensions be denoted by~\xv, and when relevant let position in the thin dimensions be denoted by~\zv.%
\footnote{The \(z\) could be a discrete finite index, as in the single scalar field cases of \cref{SoneDintro,Shceq}, or it could be a continuous domain such as the cross-section of a beam, shell, or channel (e.g., \cref{EGsd}).   
When a continuum, the \zv-dimensions need not necessarily be physically thin.  
Instead we just need the dynamics in the \zv-directions to be like those we associate with `thin' domains: e.g.,  quasi-stationary distributions of multiscale Fokker--Planck \pde{}s \cite[e.g.,][\S18 and \S21.2 resp.]{vanKampen85, Roberts2014a}.}
Here we homogenise the dynamics away from the boundary~\(\partial\Omega\) of the macroscale spatial domain~\(\Omega\subseteq\RR^d\) so we consider \(\xv\in\XX\subset\Omega\) for some spatial domain~\XX\ of interest that does not include boundary layers.
Let the thin domain of~\zv\ be denoted by~\ZZ.
Let the field(s) of interest be a function of~\(t,\xv\) such that \(u(t,\xv)\in\HH_\ZZ\) for some Hilbert space~\(\HH_\ZZ\subset L^2(\ZZ)\)\ that contains the \zv-dependence.
For most of the following, the \zv-structure is implicit via this Hilbert space of~\(u\): this implicit dependence empowers us to focus on the multiscale character of the \xv-dependence in the large domain~\XX.  
The class of heterogeneous problems we address is then of the general form
\begin{equation}
\alphaD u=\fL(\xv,\thetav) u 
-\divv \fv(\xv,\thetav,u,u_\xv)
+\gamma g(t,\xv,\thetav,u,u_\xv)
\quad\text{for }\thetav:=\ellm^+\xv\,,
\label{Egenpde}
\end{equation}
where \(\alphaD{}\) denotes a time evolution operator as introduced by \cref{SoneDintro} for \pde~\cref{Ehdifpde}, and where the right-hand side is \(1\)-periodic in~\thetav.
The linear operator \(\fL(\xv,\thetav):\HH_\ZZ\to\HH_\ZZ\) encapsulates many purely \zv-direction processes, and may parametrically depend upon~\xv,\thetav\ as indicated.
The gradient operator \(\grad:=(\D{x_1}{},\ldots,\D{x_d}{})\), whereas \(\divv\)~denotes the corresponding divergence.
The `flux' function~\fv\ and the `forcing' function~\(\gamma g\) \emph{may} both be nonlinear functions of the field~\(u\) and its gradient \(u_\xv:=\grad u\) (for generality, derivatives may be interpreted in the weak sense).
We assume that the form of~\(\fL\C \fv\C g\) are such that there exist general solutions~\(u(t)\) of~\cref{Egenpde} in a weighted Sobolev \text{space~\(W^{1,2}(\XX)\times\HH_\ZZ\) for every \(t\in\TT\)\,.}

\cite{Fish2021} [p.775] commented that the ``engineering counterpart [homogenisation] based on the so-called Hill--Mandel macrohomogeneity condition assumes equivalency between the internal virtual work at an \rve\ level and that of the overall coarse-scale fields.''  
In contrast, our approach here makes \emph{no} such assumption and so applies to a much wider range of systems, such as \text{the class~\cref{Egenpde}.}

\begin{example}[shear dispersion]\label{EGsd}
In addition to the example heterogeneous \pde{}s~\cref{Ehdifpde,EhHetNon,Ehdifpde2},
a straightforward example of~\cref{Egenpde} is the shear dispersion in a 2D~channel, long in the \(x\)-direction and narrow in the \(z\)-direction (say \(|z|<1\)), and with heterogeneous advection-diffusion. 
The concentration~\(u(t,x,z)\) of the material is governed by the following (non-dimensional) \pde\ in the form of~\cref{Egenpde}:
\begin{equation*}
    \partial_t u = \underbrace{\partial_z[\kappa(z)u_z]}_{\fL u}
-\partial_x[\underbrace{v(z) u-\kappa(z)u_x}_{\fv}]
+\underbrace{0}_{\gamma g}\,,
\end{equation*}
where the diffusive mixing~\(\kappa(z)\) varies in~\(z\), and the advection velocity~\(v(z)\)\ has shear in~\(z\), such as the classic parabolic profiles \(\kappa,v\propto 1-z^2\).
For this shear dispersion, \cite{Watt94b} showed how to develop multi-mode, multi-continuum models for the emergent macroscale dynamics.
A derived low-order bi-continuum model was found to be that the concentration \(u\approx U_0(t,x)+(3z^2-1)U_2(t,x)\) for the homogenised \pde{}s  
\begin{equation*}
    \partial_tU_0\approx -\bar vU_{0x}+\tfrac25\bar vU_{2x}\,,
\qquad \partial_tU_2\approx -6U_2+\tfrac12\bar vU_{0x}\,,
\end{equation*} 
in terms of the average advection~\(\bar v:=\overline{v(z)}\) (their~(2.22)--(2.24)).  
Further, \cite{Strunin01a} discussed interpreting such two-mode bi-continuum models as physical zonal models---realised by parametrising the same \im\ through choosing a different definition of amplitudes~\(U0,U_2\).
\qed\end{example}

A second example application in the class~\eqref{Egenpde} is the one-mode modelling of Taylor dispersion in a channel with wavy walls, see Fig.~2.1 by \cite{Rosencrans93}.
A nonlinear example of~\eqref{Egenpde} with forcing is the accurate two-mode bi-continuum modelling of the inertial dynamics of thin fluid flow over a substrate surface which is arbitrarily curved over macroscale space \cite[]{Roberts99b}.
A fourth example is that \cite{Guinovart2024} developed an asymptotic homogenisation of piezoelectric composite materials in generalised curvilinear coordinates that would also fit in the class of~\eqref{Egenpde}: although they only considered equilibrium problems whereas we model dynamics; and they only addressed homogenisations obtained by cell-averaging which implicitly can only encompass cell modes that are essentially constant within a cell, whereas we encompass any selection of physics-informed sub-cell structures.

\subsubsection{Multiscale nature}
The appearance of the repeated dependence upon space~\xv\  in \pde~\cref{Egenpde}, both directly via \xv, and indirectly via \(\thetav=\ellm^+\xv \pmod1\), is a consequence of the multiscale spatial structure of the material. 
To reflect multiscale structure, the \pde~\cref{Egenpde} poses that the spatial variations of the coefficients~\(\fL,\fv,\gamma g\) may occur on both a macroscale directly via~\xv\ and on microscales using~\(\thetav=\ellm^+\xv\), where \(\ellm^+\)~is defined via~\cref{Eellm}. 
The macroscale spatial variations cater for functionally graded materials\footcite[e.g.,][\S6.1]{DaChen2024, Anthoine2010, Roberts2022a} or in the nonlinear modulation of spatial patterns\footcite[e.g.,][]{Cross93, Newell93, Roberts2013a}. 
Whereas the microscale heterogeneity to be homogenised has spatial variations represented via the vector of phase variables~\(\thetav:=\ellm^+\xv \pmod 1\) corresponding to the phase variable in \cref{SoneDintro,Shceq}.
We aim to prove the existence and construction of \emph{closed and accurate} models of the macroscale dynamics of \pde~\cref{Egenpde} via a purely-macroscale varying, system-level, multi-modal, field~\(\Uv(t,\xv)\in\RR^M\) satisfying a homogenised macroscale \pde\ system of the form 
\begin{equation}
\alphaD \Uv =  \Gv(t,\xv,\Uv,\Uv_\xv,\Uv_{\xv\xv}, \ldots),
\label{Eepde}
\end{equation}
for some purely-macroscale functional~\(\Gv\) (which implicitly depends upon parameters such as~\(\gamma\)).

\begin{example}[macroscale functional graduations]\label{EGmfg}
A very simple illustration of homogenisation with macroscale variations in material properties is to consider the 1-D heterogeneous \cref{EG1D} with coefficient \(\kappa(x):=1/[1+a(x)\cos x]\) where the strength~\(a(x)\), instead of being constant, now varies gradually in~\(x\), that is, over lengths longer than the microscale period~\(2\pi\) of~\(\cos x\)\,.  
The computer algebra of \ifJ Supplementary Code \cite[Appendix~A,][]{Roberts2024a} \else\cref{cas} \fi constructs the resultant functionally graded homogenisation by the inclusion of one extra command that simply asserts that parameter~\(a\) depends upon~\(x\).
The resulting one-mode homogenised model is 
\begin{equation*}
\alphaD U_0 = U_{0xx} +\tfrac12\left[a(x)^2U_{0xxx}\right]_x +\Ord{\partial_x^5,a^5}.
\end{equation*}
This homogenisation is particularly simple since for this form of~\(\kappa(x)\) we know the leading coefficient on the \textsc{rhs} is one for all strengths~\(a\), and so it is no surprise that it remains independent of macroscale varying~\(a\).
Here the macroscale variations in heterogeneity strength only affects higher-order gradient terms, such as the shown \text{fourth-order term.}  
\qed\end{example}

\begin{assumption}[smoothness]\label{assFandG}
The operator~\fL\ and functions~\fv\ and~\(g\) on the right-hand side of \pde~\cref{Egenpde} are to be smooth functions of their arguments~\(\xv,u,u_\xv\).
We define \emph{smooth} to mean continuously differentiable to an order~\(q\) sufficient for the purposes at hand, uniformly~\(C^q\) for some~\(q\), or possibly infinitely \text{differentiable,~\(C^\infty\).}
\end{assumption}

The \thetav-dependence in \pde~\cref{Egenpde} need not be so smooth.  
An example being the piecewise constant coefficient~\(\kappa(x)\) in the high-contrast example of \cref{Shceq}.
The crucial constraint on the \(\thetav\)-dependence is that \cref{Aeig} on a general eigen-decomposition needs to be met.
Nonautonomous effects via the \(t\)-dependence in~\(g\) need not necessarily be smooth: we need certain integrals to be bounded, which is often the case provided the \(t\)-dependence is measurable \text{\cite[e.g.,][]{Aulbach96, Roberts06k}.}

\subsubsection{Microscale heterogeneity}
We suppose that the microscale heterogeneity in~\(\xv\in\RR^d\), represented via the variable~\(\thetav=\ellm^+\xv\), is possibly quasi-periodic with some number~\(P\) of incommensurable vector periods~\(\ellv_p\in\RR^d\) for \(p=1,\ldots,P\). 
In contrast to \cite{AbdoulAnziz2018} we do \emph{not} assume \(P\leq d\) but cater for much more general heterogeneity, such as quasi-periodic \cite[e.g.,][]{Roberts2022a}.
For example, a 3-D bulk material with microscale heterogeneity varying \(\ell\)-periodically in each direction (a cubic cell) has periods \(\ellv_1:=(\ell,0,0)\), \(\ellv_2:=(0,\ell,0)\), and \(\ellv_3:=(0,0,\ell)\); whereas if the \(x_1\)-direction was instead quasi-periodic with periods~\(\ell\) and~\(\ell/\sqrt2\), then include the additional fourth vector period \(\ellv_4:=(\ell/\sqrt2,0,0)\). 
Define both the \(d\times P\) matrix%
\footnote{For the given four-vector quasi-periodic example, the \(3\times4\) matrix~\ellm\ and its Moore--Penrose pseudo-inverse are
\(\ellm=\ell\begin{bmat} 1&0&0&1/\sqrt2\\ 0&1&0&0\\ 0&0&1&0 \end{bmat}\) and 
\(\ellm^+=\dfrac1\ell\begin{bmat} 2/3&0&0\\ 0&1&0&\\ 0&0&1\\ \sqrt2/3&0&0 \end{bmat}\), for which \(\ellm\ellm^+=I_3\).}
\begin{equation}
\ellm:=\begin{bmatrix} \ellv_1&\cdots&\ellv_P \end{bmatrix},
\quad\text{and }\ellm^+:={}\text{its Moore--Penrose pseudo-inverse}
\label{Eellm}
\end{equation}
\cite[e.g.,][]{Golub1996}. 
This pseudo-inverse~\(\ellm^+\) appears in the general system~\cref{Egenpde}.
In the approach here, the appearance of~\(\ellm\) often parallels that of the asymptotically small parameter~\(\epsilon\) in asymptotic homogenisation methods.
The pseudo-inverse~\(\ellm^+\) correspondingly parallels that of~\(1/\epsilon\).
However, in contrast to other methods, we do \emph{not} invoke limits \(|\ellv_p|
\to0\)\,: all the vector periods~\(\ellv_p\) are some fixed physical microscale displacements in~\(\RR^d\) that happen to be relatively small compared to the length~\(L\) of the macroscales of interest.
Our approach and results here apply to the physically relevant cases of finite scale separation ratios~\(|\ellv_p|
/L\).
The results are \emph{not} restricted to the \text{mathematical limits \(|\ellv_p|
/L\to0\)\,.}

\subsection{Phase-shift embedding}
\label{Spse}

Generalising \cref{secpse}, and in a novel, rigorous and efficient twist to the concept of a Representative Volume Element, let's embed any specific given physical \pde~\eqref{Egenpde} into a family of \pde\ problems formed by the ensemble of all phase-shifts of the (quasi-)periodic microscale.
Such embedding is cognate to that used for quasicrystals and quasi-periodic systems in multi-D space by \cite{KaiJiang2024, KaiJiang2014}.  
Their approach and methods are in a global Fourier space and so do not appear to cater for macroscale spatial modulation of the microscale heterogeneity, nor in the solution, nor for the general class of dynamic problems~\cref{Egenpde} considered here.
\cite{Rokos2019} used a cognate family of phase-shifts of the material shown in \cref{FegDeform} in order to compute its deformed equilibria.
\cite{Smyshlyaev2000} noted the interpretation of cell-averaging as an ensemble average, but again in just addressing materials in equilibrium.
\cite{Milton2007} [p.868, eqn.~(4.3)] also commented upon ensembles of phase-shifts, but limited analysis to ensemble averages.
Indeed the algebra consequent to the embedding invoked here parallels that of many other two-scale analyses, such as that by \cite{Guinovart2024}, but complementing such previous approaches, the phase-shift embedding provides a clear physical interpretation of the algebra---an interpretation that is valid at the finite scale separation of real applications, instead of relying upon the mathematical artifice of the \text{limit of infinite scale separation.}

\begin{figure}
\centering
\caption{\label{Fpattensemble}%
schematic domain of the multiscale embedding \pde~\eqref{Eempde} for a field \(\fu(t,\x,\thetav)\), for \(\x\in\XX\subset\RR\) and for \(\thetav\in\Theta:=[0,1]^2\). 
Here the periodicities \(\ell_1=1.62\) and \(\ell_2=0.72\) so \(\ellm=\begin{bmatrix} 1.62&0.72 \end{bmatrix}\).  
We obtain solutions of the heterogeneous \pde~\eqref{Egenpde} on such blue lines as \(u_\phiv(t,x):=\fu(t,x,\phiv+\ellm^+x)\) for every constant phase~\(\phiv\in\RR^2\), here \(\phiv=(0.82,0.32)\), and where the third argument of~\fu\ has components modulo~\(1\).}
\inPlot{Figs/phaseShiftEmbedding}
\end{figure}%
As indicated by the schematic case illustrated in \cref{Fpattensemble}, let's create the desired phase-shift embedding by considering a field~\(\fu(t,\fx,\thetav)\), implicitly depending on~\zv, and satisfying the \pde%
\footnote{Many systems not in the quite general form~\cref{Egenpde} may be similarly embedded.  
The rule for derivatives is that a gradient \(\grad\mapsto \gradx+\ellm^{+T}\gradq\) whereas a divergence \(\divv\mapsto \divx+\divq\ellm^+\).
The fraktur~\fL\ introduced in~\cref{Eempde} is distinct from the mathcal~\cL\ defined for the cell-problem~\cref{Eemeig}.
}
\begin{align}
\alphaD\fu &=\fL(\fx,\thetav)\fu 
-(\divx+\divq\ellm^+) \fv\big(\fx,\thetav,\fu,\fu_\fx+\ellm^{+T}\fu_\thetav\big)
\notJbreak{}
+\gamma g\big(t,\fx,\thetav,\fu,\fu_\fx+\ellm^{+T}\fu_\thetav\big),
\label{Eempde}
\end{align}
in the domain \(\dom:=\XX\times\Theta\times\ZZ\) for the unit \(P\)-cube \(\Theta:=[0,1]^P\), and with boundary conditions of \(1\)-periodicity in~\(\theta_p\).
The subscripts~\fx\ and~\thetav\ denote the respective gradient operator, that is, \(\fu_\fx:=\gradx\fu\) and \(\fu_\thetav:=\gradq\fu\)\,.
We assume that the heterogeneous explicit dependence upon~\(\fx,\thetav\) in~\(\fL\C \fv\C g\) are regular enough that general solutions~\fu\ of \pde~\cref{Eempde} are in \(\HH^N_\dom :=W^{N+1,2} (\XX)\times\HH_\Xv\) for some chosen order~\(N\), and for a space \(\HH_\Xv\subset L^2(\Theta\times\ZZ)\) satisfying \cref{Aeig}. 
The domain~\dom\ (\cref{Fpattensemble}) is multiscale as it is large in~\fx, and relatively `thin' in both~\zv\ and~\thetav.
I emphasise that this domain has finite aspect ratio: we do \emph{not} take any limit involving an aspect ratio tending to zero \text{nor to infinity.}

\begin{figure}\centering
\caption{\label{Fcells}%
(left)~material with periodic holes distributed \(\ell\)-periodic in 2-D space; (right)~is embedded in a 4-D space, here schematic, with the $x_1x_2$-plane represented by one axis, and the holes occurring across the $\ell\times\ell$ $\theta_1\theta_2$-cross-section.}
\begin{tabular}{@{}cc@{}}
\setlength{\unitlength}{0.7ex}
\begin{picture}(40,39)
\put(4,3){
\put(17,-3){$x_1$}
\put(-4,17){$x_2$}
\put(0,1){\colorbox{blue!30!white}{\makebox(34,34){}}}
\color{blue}
\multiput(0,0)(6,0)7{\line(0,1){36}}
\multiput(0,0)(0,6)7{\line(1,0){36}}
\color{white}
\multiput(3,3)(6,0)6{\multiput(0,0)(0,6)6{\circle*3}}
\color{red!80!black}
\multiput(3,3)(6,0)6{\multiput(0,0)(0,6)6{\circle3}}
}\end{picture}
&
\inPlot{Figs/cellxsec2}
\end{tabular}
\end{figure}
Material with exclusion holes, or rigid inclusions, are also encompassed in this phase-shift embedding.
For example, \cref{Fcells} schematically shows how a 2-D material with periodic circular holes may be embedded in 4-D space where the holes manifest themselves as a cylindrical exclusion from the 4-D~domain.
The sub-cell problems are then in the \(\theta_1\theta_2\)-cross-section with the hole excluded.

\cref{Fpattensemble} indicates that we regard \(\fx=\xv\)\,.
The distinction between \fx~and~\xv\ is that partial derivatives in~\fx\ are done keeping~\thetav\ constant (e.g., parallel to the \(x\)-axis in \cref{Fpattensemble}), whereas partial derivatives in~\xv\ are done keeping the phase-shift~\phiv\ constant (e.g., along the (blue) diagonal lines in \cref{Fpattensemble}).

\begin{lemma}\label{lemeqv}
For every solution~\(\fu(t,\fx,\thetav)\in\HH^N_\dom\) of the embedding \pde~\eqref{Eempde}, and for every vector of phases~\phiv, the field \(u_\phiv(t,\xv):=\fu(t,\xv,\phiv+\ellm^+\xv)\) (for example, the field~\(\fu\) evaluated on the solid-blue lines in \cref{Fpattensemble}) satisfies the heterogeneous, phase-shifted, \pde
\begin{align}
\alphaD{u_\phiv} &=\fL(\xv,\phiv+\ellm^+\xv) u_\phiv 
-\divv \fv(\xv,\phiv+\ellm^+\xv,u_\phiv,\grad u_\phiv)
\notJbreak{}
+\gamma g(t,\xv,\phiv+\ellm^+\xv,u_\phiv,\grad u_\phiv).
\label{Eshpde}
\end{align}
Hence \(u_\ov(t,\xv):=\fu(t,\xv,\ellm^+\xv)\) satisfies the given heterogeneous \pde~\eqref{Egenpde}.
\end{lemma}

Recall that the most common boundary conditions \emph{assumed} for microscale cells (\rve{}s) are periodic, although in the usual homogenisation arguments other boundary conditions appear equally as valid despite giving slightly different results~\cite[e.g.,][]{Mercer2015}. 
For example, \cite{Liupekevicius2024} [pp.9,10] ``impose periodic conditions on the micro-fluctuation for the solid phase \ldots\ and zero uniform conditions on the micro-fluctuation for the fluid boundary'', without any apparent good justification for such varied imposition.
In contrast, here the boundary conditions of \(1\)-periodicity in microscale~\thetav\ are \emph{not} assumed but instead arise naturally due to the ensemble of phase-shifts.
That is, what in other approaches has to be assumed, in this approach \text{arises naturally.}

\begin{proof}
Start by considering the left-hand side of \pde~\eqref{Eshpde}, namely the time evolution operator
\begin{align*}
\alphaD{u_\phiv}&
=\alphaD{}\fu(t,\xv,\phiv+\ellm^+\xv)
\\&
=\left[\alphaD\fu\right]_{(t,\xv,\phiv+\ellm^+\xv)}
\qquad\text{(which by \pde~\eqref{Eempde} becomes)}
\\&
=\left[\fL\fu
-(\divx+\divq\ellm^+) \fv(\fx,\thetav,\fu,\fu_\fx+\ellm^{+T}\fu_\thetav)
\right.\notJbreak{}\left.{}
+\gamma g(t,\fx,\thetav,\fu,\fu_\fx+\ellm^{+T}\fu_\thetav)
\right]_{(t,\xv,\phiv+\ellm^+\xv)}
\\&
=\left[\fL(\fx,\thetav)\fu\right]_{(t,\xv,\phiv+\ellm^+\xv)}
-\divv\left\{\left[ \fv(\fx,\thetav,\fu,\fu_\fx+\ellm^{+T}\fu_\thetav)
\right]_{(t,\xv,\phiv+\ellm^+\xv)}\right\}
\nonumber\\&\qquad{}
+\gamma g(t,\xv,\xv+\phiv,u_\phiv,\grad u_\phiv)
\\&
=\fL(\xv,\phiv+\ellm^+\xv)u_\phiv
-\divv \fv(\xv,\phiv+\ellm^+\xv,u_\phiv,\grad u_\phiv)
\notJbreak{}
+\gamma g(t,\xv,\phiv+\ellm^+\xv,u_\phiv,\grad u_\phiv),
\end{align*}
namely the right-hand side of~\eqref{Eshpde}.
Hence, provided \pde~\eqref{Eempde} has boundary conditions of \(1\)-periodicity in~\(\theta_p\), every solution of the embedding \pde~\eqref{Eempde} gives a solution of \pde~\eqref{Eshpde}, namely the  \pde~\eqref{Egenpde} for any \(P\)-dimensional phase-shift~\phiv\ of the heterogeneity.

In particular, \(u_\ov(t,\xv)\), of phase-shift \(\phiv=\ov\), satisfies the given heterogeneous \pde~\eqref{Egenpde}.
\hfill\end{proof}

\begin{lemma}[converse]\label{lemcon}
Suppose we have a family of solutions~\(u_\phiv(t,\xv)\) of the phase-shifted \pde~\eqref{Eshpde}---a family parametrised by the phase vector~\(\phiv\in\RR^p\)---and the family depends smoothly enough upon~\(t,\xv,\phiv\) such that the following \(\fu\in\HH^N_\dom\).
Then the field \(\fu(t,\fx,\thetav):=u_{\thetav-\ellm^+\fx}(t,\fx)\) satisfies the \text{embedding \pde~\eqref{Eempde}.}
\end{lemma}

\begin{proof} 
First, from the \pde~\eqref{Eempde}, consider
\begin{align*}
\fu_\fx+\ellm^{+T}\fu_\thetav
&=\left[(-\ellm^+)^T\D{\phiv}{u_\phiv}+\D \xv{u_\phiv}+\ellm^{+T}\D{\phiv}{u_\phiv}
\right]_{\phiv=\thetav-\ellm^+\fx,\xv=\fx}
\notJbreak{}
=\big[\grad{u_\phiv}\big]_{\phiv=\thetav-\ellm^+\fx,\xv=\fx}\,.
\end{align*}
Second, since \(\phiv=\thetav-\ellm^+\fx=\thetav-\ellm^+\xv\)\,, that is \(\thetav=\phiv+\ellm^+\xv\)\,, then for every smooth~\(\fv(\fx,\thetav)\),
\((\divx+\divq\ellm^+)\fv
=\divv\big\{f|_{\fx=\xv,\thetav=\phiv+\ellm^+\xv}\big\}\).
Thirdly, hence the right-hand-side of \pde~\eqref{Eempde} becomes
\begin{align*}
&
\fL(\fx,\thetav) \fu
-(\divx+\divq\ellm^+) \fv(\fx,\thetav,\fu,\fu_\fx+\ellm^{+T}\fu_\thetav)
\notJbreak{}
+\gamma g(t,\fx,\thetav,\fu,\fu_\fx+\ellm^{+T}\fu_\thetav)
\\&
=\big[\fL(\fx,\thetav) \fu\big]|_{\fx=\xv,\thetav=\phiv+\ellm^+\xv}
-\divv\big\{ \fv(\fx,\thetav,\fu,\fu_\fx+\ellm^{+T}\fu_\thetav)|_{\fx=\xv,\thetav=\phiv+\ellm^+\xv}\big\}
\nonumber\\&\quad{}
+\gamma g(t,\fx,\thetav,\fu,\fu_\fx+\ellm^{+T}\fu_\thetav)|_{\fx=\xv,\thetav=\phiv+\ellm^+\xv}
\\&
=\fL(\xv,\phiv+\ellm^+\xv) u_\phiv
-\divv \fv(\xv,\phiv+\ellm^+\xv,u_\phiv,\grad u_\phiv)
\notJbreak{}
+\gamma g(t,\xv,\phiv+\ellm^+\xv,u_\phiv,\grad u_\phiv)
\,,
\end{align*}
namely the right-hand side of \pde~\eqref{Eshpde}.
Lastly, since \(\alphaD\fu(t,\fx,\thetav)=\alphaD{u_{\thetav-\ellm^+\fx}}\) it follows that \(\fu:=u_{\thetav-\ellm^+\fx}(t,\fx)\) satisfies the embedding \pde~\eqref{Eempde}.
\hfill\end{proof}

Consequently, \pde{}s~\eqref{Eempde,Egenpde} are equivalent, and they may provide us with an ensemble of solutions for an ensemble of materials all with the same heterogeneity structure, but with the structural phase of the material shifted through all possible phases.
The key difference between \pde{}s~\cref{Egenpde,Eempde} is that although \pde~\eqref{Egenpde} is heterogeneous in space~\xv, the embedding \pde~\eqref{Eempde} is \emph{homogeneous} in space~\fx.
Because of this homogeneity, \cref{Simmmh} is empowered to apply an existing rigorous theory for gradual variations in space that leads to desired multi-continuum homogenised models of the \pde~\cref{Eempde}, and \text{thence to that of~\eqref{Egenpde}.}

\subsection{Invariant manifolds of multi-continuum, micromorphic, any-order homogenisation}
\label{Simmmh}

Generalising \cref{secsmh}, let's analyse the embedding \pde~\eqref{Eempde} for useful `homogenised' \im{}s.  
Such \im{}s are to express and support the relevance of a  potential hierarchy of accurate homogenised models for the original heterogeneous \pde~\eqref{Egenpde}.  
The systematic approach developed simplifies considerably much of the ``difficulty to choose a priori an appropriate model for a given microstructure'' discussed by \cite{Alavi2023} [p.2164].

Developments in dynamical systems theory inspired by earlier more formal arguments\footcite{Bunder2018a, Roberts2016a, Roberts2013a, Roberts88a, Roberts96a} establishes how to construct a \pde\ model for the macroscale spatial structure of \pde\ solutions in multiscale domains~\dom\ such as \cref{Fpattensemble}.
The technique is to base analysis on the case where variations in~\fx\ are over a large enough scale that they are approximately negligible---the variations are both directly in the parametric~\fx\ dependence of~\(\fL\C \fv\C g\) and indirectly via the field~\(\fu\).
Then we treat finite, macroscale, variations in~\fx\ as a regular perturbation.
Despite the gradient~\(\grad_\fx\) being an unbounded operator, the theoretical developments justify being able to treat such derivatives as `small' \cite[e.g.,][p.987]{Roberts2016a}.
Hence the \emph{regular} perturbation analysis proceeds to any chosen order~\(N\) in the `small' gradients~\(\grad_\fx\), and with quantifiable remainder error \cite[(52)]{Roberts2016a}.

For two examples in linear elasticity systems, the usual leading order homogenisations are the case \(N=2\), and the so-called second-order homogenisations\footcite[e.g.,][]{Anthoine2010, Cornaggia2020, AewisHii2022} correspond to the higher-order case \(N=4\).
Because of the power of the established dynamical system framework, here we allow \text{arbitrary order~\(N\).}

\subsubsection{Linear basis of invariant manifolds}
\label{SSSlbim}

\begin{quoted}{The Little Prince, Antoine de Saint Exup\'ery}
it is only with the heart that one can see rightly; 
what is essential is invisible to the eye.
\end{quoted}

Invariant manifolds (\im{}s) are mostly constructed from the base of an equilibrium or a family of equilibria, as we do here to generalise \cref{secsmh}.
In the vicinity of each and every equilibria we characterise all solutions in terms of spectral properties of the system's linearisation.
It is these outwardly invisible spectral properties that lie at the heart of proven modelling.
We choose an appropriate subset of modes based upon the spectrum of eigenvalues and the intended purpose of the model.
This choice empowers us to construct corresponding \im{}s that pass through the base equilibria and extend out into the state space.
The correspondingly constructed evolution then forms a closed, accurate, in-principle exact, \text{macroscale homogenisation.}

\paragraph{Equilibria}
Following \cite{Bunder2018a}, consider the dynamics of the embedding \pde~\cref{Eempde} in a mesoscale locale of each and every `cross-section' position \(\fx=\Xv\in\XX\) of interest in the physical problem at hand.
In such a mesoscale locale the variations in macroscale variable~\fx\ are small enough so that for \emph{linearisation} purposes we consider the macroscale gradients negligible, \(\gradx\equiv0\).
Many secondary physical nonlinearities or forcing effects are gathered into \(g(t,\fx,\thetav,\fu,\fu_\fx+\fu_\thetav)\), multiplied by perturbation parameter~\(\gamma\), so we seek equilibria with coefficient \(\gamma=0\)\,.
That \im\ theory supports modelling in the vicinity of each equilibria\footcite{Aulbach2000, Prizzi02, Hochs2019, Bunder2018a} empowers us to systematically and accurately model effects with non-zero gradients~\(\gradx\), non-zero~\(\gamma\), \text{and with nonlinearities.}

With effectively zero~\(\gradx\) and zero~\(\gamma\) the dynamics of the embedding \pde~\cref{Eempde} in the locale near~\(\Xv\) reduces to the cross-sectional, \emph{cell-problem}
\begin{equation}
\alphaD\fu = \fL(\Xv,\thetav)\fu -\divq \fv(\Xv,\thetav,\fu,\ellm^{+T}\fu_\thetav),
\quad \fu\text{ is 1-periodic in }\theta_p\,,
\label{EempdeX}
\end{equation}
at every cross-section~\(\Xv\in\XX\)\ of interest.
For each~\Xv, the cross-sectional \pde~\cref{EempdeX} is a \emph{cell-problem} in that its right-hand side only contains derivatives in~\thetav\ and implicitly~\zv.
The following treatment of \pde~\cref{EempdeX} is parametrised by the macroscale locale~\Xv, and so many of the quantities identified may depend upon~\Xv, as in functionally graded materials\footcite[e.g.,][\S6.1]{DaChen2024, Anthoine2010, Roberts2022a}, although often not.
For brevity, any such \Xv-dependence is mostly implicit in the following.

We assume that the cell-problem~\cref{EempdeX} has one or more chosen equilibria \(\fu=\fu^*(\thetav)\in\HH_\Theta\times\HH_\ZZ\); sometimes these equilibria depend upon location~\Xv, but for simplicity we leave such dependence implicit.
Often these equilibria will zero the flux~\fv.
Often the equilibria~\(\fu^*\) are constant in~\thetav\ and in~\zv---but they need not be constant.
Often the chosen equilibria form a subspace such as \(\fu^*\propto v_0(\thetav)\) for some~\(v_0\).
Define \begin{equation}
\EE:= \{\text{chosen equilibria } \fu^*\} \subset \HH_\Theta\times\HH_\ZZ\,;
\label{Eequil}
\end{equation}
that is, \(\fL\fu^*-\divq\fv(\Xv,\thetav,\fu^*,\ellm^{+T}\fu^*_\thetav)=0\) \text{for every~\(\fu^*\in\EE\).}

In application, one often has useful physical intuition about a suitable base set of equilibria~\EE\ for the scenarios of interest.
One then introduces the artificial ordering parameter~\(\gamma\) into the governing equations so that \(\fL\fu+\divq\fv\) has the desired equilibria~\EE, and all other terms are gathered into the `perturbing'~\(\gamma g\) \cite[e.g.,][\S9.1 and Part~V]{Roberts2014a}.
The systematic framework here empowers arbitrary order construction in such an artificial~\(\gamma\) (\cref{SSSccimmch}) to often enable accurate prediction at the physically relevant~\(\gamma\) (usually at \(\gamma=1\)).

\paragraph{Characterise nearby dynamics via linearisation}
For each \(\fu^*\in\EE\), explore the nearby dynamics by seeking solutions to the cell-problem~\cref{EempdeX} in the form \(\fu=\fu^*(\thetav)+\fu'(t,\thetav)\) for small~\(\fu'\).
Then, invoking Fre\'chet derivatives, the flux becomes
\begin{align*}
\fv&=\fv(\Xv,\thetav,\fu^*+\fu',\ellm^{+T}[\fu^*_\thetav+\fu'_\thetav])
\\&\approx \fv(\Xv,\thetav,\fu^*,\ellm^{+T}\fu^*_\thetav)
+\underbrace{\D \fu\fv(\Xv,\thetav,\fu^*,\ellm^{+T}\fu^*_\thetav)}_{{}=:\jv^*(\Xv,y)}\fu'
+\underbrace{\D {\fu_\thetav}\fv(\Xv,\thetav,\fu^*,\ellm^{+T}\fu^*_\thetav)}_{{}=:J^*(\Xv,\thetav)}\fu'_\thetav\,.
\end{align*}
To characterise general dynamics near the equilibria~\EE, we thus address the linearised cell-problem  
\begin{align}&
\alphaD v=\cL v\,,\quad
\text{where }\cL v:=\fL v-\divq(\jv^*v+J^*v_\thetav),\quad 
1\text{-periodic in }\theta_p\,,
\nonumber\\&
\text{and its corresponding eigen-problem}\quad
\lambda v=\cL v\,,
\label{Eemeig}
\end{align}
for cell eigenvalues~\(\lambda\) and cell eigenvectors~\(v(\thetav)\).
Often the set~\EE, \cref{Eequil}, is chosen so that these cell-problems are independent of the equilibria \(\fu^*\in\EE\)\,.

\begin{assumption}[eigen-decomposition]\label{Aeig}
Assume the following.
\begin{enumerate}
\item A non-empty~\XX\ and~\EE\ exist where for every cross-section \(\Xv\in\XX\) and every equilibrium \(\fu^*\in\EE\), there exists for the cell eigen-problem~\cref{Eemeig} a complete countable set of (generalised) eigenvectors~\(v_m(\thetav)\) for corresponding eigenvalues~\(\lambda_m\) (sometimes complex valued), for index \(m=0,1,2,\ldots\) (often ordered so that \(\Re\lambda_{m+1} \leq \Re\lambda_m\), and if \(\Re\lambda_{m+1} = \Re\lambda_m\) then \(|\Im\lambda_{m}|\leq|\Im\lambda_{m+1}|\)).

\item There exists a countable set of corresponding (generalised) adjoint eigenvectors~\(w_m(\thetav)\) normalised so that \(\left<w_m,v_m\right>=\delta_{m,n}\)\,, such that the linear operator defined in~\cref{Eemeig} may be written as \(\cL=\sum_{m=0}^\infty \big[ \lambda_mv_m\left<w_m,\cdot\right> +\zeta_mv_m\left<w_{m+1},\cdot\right> \big]\) in a space denoted~\(\HH_\Xv\subset H:=L^2(\Theta\times\ZZ)\) (where \(\zeta_m=0\) when \(\lambda_m\neq\lambda_{m+1}\) so that this series for~\(\cL\fu\) is effectively a Jordan form of~\cL).
Define the space \(\HH_\Xv\)~as the closure of~\(\big\{\fu\in H:{}\)this series for~\(\cL u\) is an absolutely convergent series\(\,\big\}\).

\item These eigenvectors and eigenvalues are~\(C^{N+1}\) \text{in macroscale~\(\Xv\in\XX\).}
\end{enumerate}
\end{assumption}

Often the cell eigen-problem~\cref{Eemeig} is independent of macroscale position~\Xv.
Examples of \Xv-dependence are functionally graded materials that have graduations in a large space dimension. 
If such a graded material has a sudden\slash step change in material properties, then such a change generates local physical `boundary' layers that have to be excised from the spatial domain~\XX.  
Such physical `boundary' layers about a material change would instead be resolved in an homogenisation via an argument akin to that for domain boundary conditions (\cref{Sbcmm}).

For each in \(\Xv\in\XX\)\ and \(\fu^*\in\EE\), \cref{Eequil}, since the eigenvectors are complete in~\(\HH_\Xv\), a general solution to the \emph{linearised} cell-problem~\cref{Eemeig} is thus%
\footnote{In cases where generalised eigenvectors occur for the same eigenvalue (\(\zeta_m\neq0\) in \protect\cref{Aeig}), then the evolution of those corresponding~\(a_m\) is more complicated, namely \(\alphaD a_m = \lambda_m a_m+\zeta_m a_{m+1}\).  
Typically there are a finite number of generalised eigenvectors for a given eigenvalues, and the evolution then includes some multiplicative factors that grow algebraically in time.  
For brevity we do not detail such cases here.}
\begin{align}&
\fu=\fu^*+\sum_{m=1}^\infty a_m(t)v_m(\thetav),
&&\text{where }
\alphaD a_m = \lambda_m a_m\,,
\label{Emode}
\end{align}
and \(a_m(t)\) denotes a general solution to this \(m\)th-mode \ode.
Now the developing argument splits depending upon the nature of~\(\alphaD\), the spectrum of eigenvalues~\(\{\lambda_m:m=0,1,2,\ldots\}\), and also upon external knowledge about the physical problem and its context.
It is not feasible to address all the myriad of possibilities for the time evolution operator~\(\alphaD\). 
Instead we focus mainly on the two main cases of first- and second-order time derivatives, \(\alpha=1,2\).

\paragraph{Perhaps use cells twice the minimum size}
But before we leave the cell-problem, recall the material and deformation shown by \cref{FegDeform}.
In the un-deformed state, such as the left and right ends of the material as shown, the periodic heterogeneous cell is clearly a square, say with side length~\(\ell\),  each microscale square with one circular inclusion.
However, choosing such a cell results in homogenisations that are almost certainly unable to predict the checkerboard pattern in the middle of \cref{FegDeform}---the necessary `zig-zag' variations would be too short for good accuracy in a macroscale homogenisation.  
To encompass this checkerboard deformation one should instead aim for a square cell of size~\(2\ell\times2\ell\) containing four inclusions; that is, embed the physical problem using heterogeneity period~\(2\ell\) in both directions.
Then a local checkerboard pattern can be chosen as one of the sub-cell modes~\(v_m(\thetav)\) in a multi-continuum micromorphic homogenisation \cite[e.g.,][]{Rokos2019}. 

\cite{Combescure2022} used cognate 1-D examples in discussing selecting generalized continuum models for materials displaying microstructure instabilities.
Indeed, recall Mathieu's equation for oscillations in~\(u(t)\) with parametric forcing, \(u_{tt}+(\omega^2+\gamma\cos t)u=0\) \cite[e.g.,][\S11.4]{Bender81}.  
This system has instabilities when the natural frequency \(\omega=k/2\) for integer~\(k\).
The strongest instability occurs when \(\omega=1/2\), that is, at twice the period of the forcing.
The spatial analogue for homogenisation is that a likely candidate for a mode in a micromorphic homogenisation is one with twice the wavelength of the underlying heterogeneity, as in \cref{FegDeform}, and captured in our systematic homogenisation via an embedding with cells of twice the minimal size.

\subsubsection{Systems with significant dissipation}
\label{SSSsd}

\begin{SCfigure}[0.6]
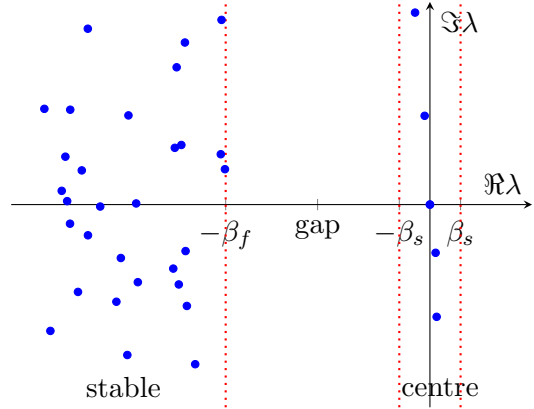
\centering
\caption{\label{Fnfseval}schematic picture of the complex plane of eigenvalues~\(\lambda_m^{1/\alpha}\) (blue dots) forming two separate sets characterised by bounding parameters~\(\beta_L\) and~\(\beta_\ell \) of the centre modes and stable modes, respectively.}
\inPlot{Figs/nfseval}
\end{SCfigure}

The case \(\alpha=1\) is the case of first-order in time \pde~\cref{Egenpde}.  
This case usually has the cleanest argument and most rigorous support. 
General solutions to the \(m\)th-mode \ode~\cref{Emode} are exponentials in time: \(a_m=A_m\e^{\lambda_mt}\).\footnote{Albeit possibly multiplied by a polynomial in~\(t\) in the case of generalised eigenvectors.}
Typically most of these cell eigenvalues~\(\lambda_m\) have large negative real-part (see the schematic example of \cref{Fnfseval}), and so these corresponding cell modes decay to zero very quickly.
It is the relatively few cell eigenvalues~\(\lambda_m\) with small real-part that determine the long-time \text{macroscale evolution.}

As indicated in \cref{Fnfseval} (and also \cref{Fnfreval0}) we suppose that there are no eigenvalues with significantly positive real-part as then the linear dynamics would predict cell mode(s) with rapid exponential growth.  
In most dynamics scenarios, such `exploding' modes ruin the usefulness of the equilibrium as a base from which to form a homogenised model.
Hence we restrict attention to equilibria in~\EE, \cref{Eequil}, whose spectrum has no eigenvalues with significantly positive real-part.

Modelling nonlinear systems requires a good gap to occur in the spectrum, although linear autonomous system are not so restrictive.  
The following definition covers the two most useful cases, as illustrated by \cref{Fnfseval,Fnfreval0}.

\begin{definition}[spectral gap] \label{Dgap}
Let \(\lambda_0,\lambda_1,\ldots\) denote the countable eigenvalues (\cref{Aeig}) of the dynamical system linearised about some equilibrium.
A \emph{spectral gap} is characterised by two bounds, say \(0\leq\beta_L<\beta_\ell \) (a large ratio~\(\beta_\ell /\beta_L\) is desirable), and commonly is identified in one of the following two ways.
\begin{enumerate}
\item \label{DIsge} If \(|\Re\lambda^{1/\alpha}_j|\leq\beta_L\) for \(j=0,\ldots,M-1\) and \(\Re\lambda^{1/\alpha}_j\leq-\beta_\ell \) otherwise, then a \emph{spectral gap} occurs between the \emph{centre modes}, indexed by \(j=0,\ldots,M-1\), and the \emph{stable modes} otherwise (\cref{Fnfseval}).

\item \label{DIsgc} If \(|\lambda^{1/\alpha}_j|\leq\beta_L\) for \(j=0,\ldots,M-1\) and \(|\lambda^{1/\alpha}_j|\geq\beta_\ell \) otherwise, then a \emph{spectral gap} occurs between the \emph{slow modes}, indexed by \(j=0,\ldots,M-1\), and the \emph{fast modes} otherwise (\cref{Fnfreval0}).
\end{enumerate}
\end{definition}

A rational multi-continuum model is formed by identifying a significant \emph{spectral gap} in the spectrum of cell-eigenvalues (e.g., \cref{Fnfseval}) that holds for all~\(\Xv\in\XX\), albeit that the bounds~\(\beta_L,\beta_\ell \) may vary with~\Xv.
In a physical application one aims to resolve macroscale time variations longer than some minimum timescale of interest~\(t_L\): in such a case one seeks a spectral gap with bounds \(\beta_L<1/t_L<\beta_\ell \) (preferably \(\beta_L\ll1/t_L\ll\beta_\ell \)).
Identifying such a gap identifies \(M\)~\emph{centre} cell-eigenvalues, as in \cref{Fnfseval}, and for convenience we index the eigenvalues so that these are \(\lambda_0,\ldots,\lambda_{M-1}\) (repeated according to multiplicity). 
Commonly there is one or two conserved cell modes (cell-eigenvalues zero) and all other eigenvalues have significantly negative real-part, whence we may choose to identify \(\beta_L=0\), and the argument here then usually leads to the common homogenised models.
But when appropriate for the physical scenarios of interest, one may alternatively choose to include more sub-cell modes in the modelling, as we encompass here. 
The corresponding \(M\)~cell-eigenvectors \(v_0,\ldots,v_{M-1}\) are \(M\)-microscale modes that form the basis of an \(M\)-D subspace~\(\MM_\Xv\) of~\(\HH_\Xv\) that is invariant to the linearised cell-dynamics~\cref{EempdeX}.
All the other cell-modes~\(v_m\), for \(m\geq M\), called \emph{stable} modes, decay faster than~\(\e^{-\beta_\ell t}\).
Hence, within the linearised dynamics the invariant subspaces~\(\MM_\Xv\) of the centre modes is a local \emph{emergent} (\cref{Dim}) subspace on the chosen macroscale times~\(t_L\) of interest.

In \emph{linear} problems with \emph{homogeneous macroscale} the  spectral gap may be small (e.g., \cref{Selastic2d}): one may choose a sharp distinction between the so-called centre and stable modes.  
The reason is that in linear, macroscale-homogeneous, problems the macroscale modes do not interact and hence do not `fillin' the gap.
But in problems with \emph{either nonlinear or macroscale-inhomogeneous} effects (e.g., \cref{SShhn}), both of which we encompass here, the macroscale modes and/or variations interact to generate effects `within' the gap.  
At high enough order these effects `cross' the gap and are then likely to cause troublesome small divisors, divisors that limit the order of approximation and the validity of the modelling.
To avoid such small divisors one needs a gap big enough for the desired \text{order of approximation.}

The nonlinear theory of \im{}s for non-autonomous systems\footcite[e.g.,][]{Aulbach2000, Potzsche2006, Haragus2011, Roberts2013a,  Hochs2019, Bunder2018a}, subject to various preconditions, asserts that under perturbation by nonlinearity and macroscale spatial modulation (\(\gradx\neq0\)) the qualitative nature of this linear picture is preserved throughout~\XX\ and~\EE.
The theoretical support is summarised in the following \cref{Pgenft,Pgenabt}.

\begin{proposition}[Forward Theory]\label{Pgenft}
Under \cref{assFandG,Aeig}, and that the cell-eigenvalues have a spectral gap in the sense of \cref{Dgap}.\ref{DIsge} (as in \cref{Fnfseval}),  and when the cell-operator~\cL, \cref{Eemeig}, generates a strongly continuous semigroup, and when~\eqref{Eempde} is non-autonomous if \cL~is bounded,   then the following holds for \text{the system~\eqref{Eempde}.}
\needspace{2\baselineskip}
\begin{enumerate}
\item  There exists an \(M\)-mode centre invariant manifold~\cM\ in some neighbourhood of the equilibria~\EE, \cref{Eequil}, such that for every~\(\fu^*\in\EE\), \cM~has tangent space~\(\MM_\Xv\).
\needspace{2\baselineskip}
\item All solutions in the neighbourhood of~\EE\ are exponentially quickly attracted to solutions on~\cM\ (approximately like~\(\e^{-\beta_\ell t}\)).
\needspace{4\baselineskip}
\item If an approximation to~\cM\ and the evolution thereon (the homogenisation~\eqref{Eepde}) satisfies~\eqref{Eempde} to a residual of order~\(N+1\) in spatial gradients, nonlinearity, and~\(\gamma\), and provided~\(N\) is constrained by the smoothness of~\(\fL,\fv,g\) and for nonlinear systems by the spectral gap \(N+1<\beta_\ell /\beta_L\) (\cref{Dgap}.\ref{DIsge}), then the approximations have errors of the \text{same order~\(N+1\).}
\end{enumerate}
\end{proposition}

\begin{proof} 
Autonomous~\eqref{Eempde} systems under the stated preconditions satisfy Assumption~3 of \cite{Bunder2018a} to establish the validity of \Xv-local analysis.
Then the existence, emergence and approximation theorems of \cite{Haragus2011, Aulbach2000, Potzsche2006} collectively apply and give a spatially-local version of the results of \cref{Pgenft} at each~\Xv, with one extra ingredient, namely that the spatial gradients give a \pde~\cite[(44)]{Bunder2018a} with a complicated remainder expression that is of the order~\(N+1\).
This existence, emergence and approximation results hold for every locale~\(\Xv\in\XX\), and so the results hold globally in the spatial domain~\XX.
For non-autonomous systems~\eqref{Eempde}, the extant theories of \cite{Haragus2011, Aulbach2000, Potzsche2006} additionally require that \cL~be bounded.
\hfill\end{proof}

The examples of \cref{SoneDintro,Shceq}  have self-adjoint linear operators for~\cL\ that generate the required strongly continuous semigroup.
Since they are also autonomous, and for the dissipative case \(\alpha=1\), the above properties hold for the constructed approximations to their \(M\)-mode, \(M\)-continuum, homogenisations. 

However, there are also many physical systems of great interest that do not satisfy the required preconditions for \cref{Pgenft} but do satisfy the less stringent preconditions of the following backward proposition.

\begin{proposition}[alternate Backward Theory]\label{Pgenabt}
Under \cref{assFandG,Aeig}, and \cL~of~\cref{Eemeig} has a spectral gap according to \cref{Dgap}.\ref{DIsge} (as in \cref{Fnfseval}), and when~\(\HH^N_\dom\) is a graded Fre\'chet space and \cL\ is a continuous linear operator on~\(\HH^N_\dom\) \cite[in the sense of][]{Hochs2019}, then the following holds for the embedding system~\eqref{Eempde}:
\needspace{2\baselineskip}
\begin{enumerate}
\item there exists a constructible smooth system in~\XX\ close to the embedding system~\eqref{Eempde}, close to within differences of order~\(N+1\), an order described  and constrained as in \cref{Pgenft}; 
\needspace{2\baselineskip}
\item which, in a finite domain of state space containing~\EE, \cref{Eequil}, \emph{exactly} possesses the constructed \(M\)-mode centre \im~\cM, with \emph{exactly} the constructed evolution on~\cM; and 
\needspace{2\baselineskip}
\item \label{PgenabtEmrg}where the centre \im~\cM\ is emergent (\cref{Dim}, and emerges approximately like~\(\e^{-\beta_\ell t}\)).
\end{enumerate}
\end{proposition}

\begin{proof} 
Under the stated preconditions, system~\eqref{Eempde} satisfies Assumption~3 of \cite{Bunder2018a} to establish the validity of \Xv-local analysis.
Under \cref{Aeig} a linear coordinate transform exists to convert~\cL\ to having orthogonal eigenvectors.
Then the theorems of \cite{Hochs2019} apply to give a spatially-local version of the results of \cref{Pgenabt} at each~\(\Xv\in\XX\), with one extra ingredient, namely that the spatial gradients give a \pde~\cite[(44)]{Bunder2018a} with a complicated remainder expression that is of the order~\(N+1\).
These existence, emergence and invariance results hold for every locale~\(\Xv\in\XX\), and so the results hold globally in the spatial domain~\XX.
\hfill\end{proof}

In these two propositions, the statements of order~\(N+1\) errors and residuals justify the iterative construction algorithm of \cref{Pconstruct}, as implemented in the computer algebra of \ifJ Supplementary Code \cite[Appendices~A to~B,][]{Roberts2024a}\else\cref{cas,cashc}\fi.

In their review, \cite{Fish2021} [p.774] commented that in ``providing a link between fine and coarse scales \ldots\ the undertaking becomes challenging for heterogeneous systems, particularly for describing large deformation and failure of materials, which often involve history-dependent mechanisms.'' 
To answer this challenge, we establish a framework and straightforward systematic construction whose proven supporting theory, such as these two propositions, encompasses large nonlinear out-of-equilibrium modelling.
Moreover, the approach often avoids the need for history-dependent mechanisms through a systematic approach to the extra kinematic variables of a \text{multi-continuum micromorphic homogenisation.}

However, in multi-D space the spectrum of eigenvalues is often `rich'---it  has \emph{many} modes with comparable eigenvalues of slowish decay.  
In such scenarios these many independent modes may combine to look like one (or a few) fractional derivative modes \cite[e.g.,][]{HongGuangSun2018}.
This combination may be especially observable when the many are primarily forced by a few dominant modes.
An argument to justify such a fewer-mode fractional calculus model of such a collective is not attempted herein.

\subsubsection{Wave-like systems} 
\label{SSSwls}

\begin{SCfigure}[2]
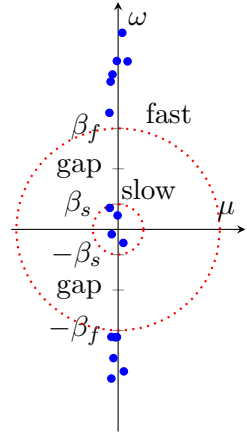
\centering
\caption{\label{Fnfreval0}schematic complex plane of cell-eigenvalues \(\mu+\i\omega\) (blue discs) in the case of wave-like dynamics when within each cell there are slow modes among fast oscillations.
This case is for small real-part of~\(\lambda_m^{1/\alpha}\).  
The cell-eigenvalues form two separate sets characterised by bounding parameters~\(\beta_L\) and~\(\beta_\ell \) of the slow and fast modes, respectively.  
}
\inPlot{Figs/nfreval0}
\end{SCfigure}

Typical wave-like systems, like the elasticity homogenisation of \cref{Selastic2d}, have  eigenvalues where~\(\lambda_m^{1/\alpha}\) are pure-imaginary, or nearly so, as illustrated schematically in \cref{Fnfreval0}.
In the second-order case, \(\alpha=2\), each eigenvalue~\(\lambda_m\) of the cell-problem~\cref{Eemeig} give rise to two linearly independent solutions of the modal \ode~\cref{Emode}, namely \(a_m(t)=A_m\e^{(\mu_m+\i\omega_m)t}+B_m\e^{-(\mu_m+\i\omega_m)t}\) for \(\mu_m+\i\omega_m:=\sqrt{\lambda_m}\).  

\paragraph{Slow manifold scenarios}
In the scenarios of \cref{Fnfreval0}, the usual physical argument is that it is the relatively few cell-eigenvalues with small frequency~\(\omega_m:=\Im\lambda_m^{1/\alpha}\) that best characterise the long-time macroscale evolution.
If so, then a rational multi-continuum model is formed by identifying a significant \emph{spectral gap} (\cref{Dgap}.\cref{DIsgc}) in the frequencies, such as that shown in \cref{Fnfreval0}.
Physically, one would usually be aiming to resolve macroscale time variations longer than some minimum timescale of interest~\(t_L\), and so seek a gap with bounds \(\beta_L<1/t_L<\beta_\ell \).
Identifying such a frequency gap identifies \(M\)~\emph{slow} cell-eigenvalues, as in \cref{Fnfreval0}, and for convenience suppose these eigenvalues correspond to \(\lambda_0,\ldots,\lambda_{M-1}\). 
But if appropriate for the physical application, then one may choose to include more sub-cell modes into the modelling as we allow here. 
The corresponding \(M\)~cell-eigenvectors \(v_0,\ldots,v_{M-1}\) are \(M\)-microscale modes that form the basis of an \(M\)-D subspace~\(\MM_\Xv\) of~\(\HH_{\Theta\times\ZZ}\).
The subspace~\(\MM_\Xv\) is invariant for the linearisation of the cell-dynamics~\eqref{EempdeX}.
Because the other cell-modes~\(v_M,v_{M+1},\ldots\), all oscillate faster than~\(\e^{\pm\i\beta_\ell t}\), within the linearised dynamics the invariant subspace~\(\MM_\Xv\), called a \emph{slow subspace}, is expected to act as a \emph{guiding centre} (\cref{Dim}) for the nearby dynamics of the \text{cell problem~\eqref{EempdeX}.}\footcite[e.g.,][]{Muncaster83c, Lorenz86, vanKampen85}  

For \emph{linear} systems~\cref{Eempde} the idea of a useful `guiding-centre' slow subspace is usually sound (because of the simplicity of linear superposition of solutions).  
The homogenisation of 2-D elasticity in \cref{Selastic2d} is an example.

However, for \emph{nonlinear} systems~\cref{Eempde} a corresponding `guiding-centre' slow invariant manifold is very delicate.
\cite{Sijbrand85} established some criteria for the existence of some such subcentre manifolds.
But there is a long running controversy in geophysical dynamics about how the non-existence of slow manifolds\footcite[e.g.,][]{Lorenz87, Lorenz92, Boyd95, Bokhove96, Vanneste2004} accords with their enormous practical usefulness.
Perhaps the most useful rigorous result\footcite[Introduced in \S13.5 of the book by][]{Roberts2014a} is the following corresponding version of the Backwards \cref{Pgenabt} with some \text{appropriate changes.}

\begin{proposition}[slow manifold]\label{Pgensm}
Under \cref{assFandG,Aeig}, and \cL~of~\cref{Eemeig} has a spectral gap according to \cref{Dgap}.\ref{DIsgc} (as in \cref{Fnfreval0}), and when~\(\HH^N_\dom\) is a graded Fre\'chet space and \cL~is a continuous linear operator on~\(\HH^N_\dom\) in the sense of \cite{Hochs2019},  
then the following holds: 
\needspace{2\baselineskip}
\begin{enumerate}
\item there exists a constructible smooth system close to the embedded system~\eqref{Eempde}, close to within errors of order~\(N+1\), an order constrained as in \cref{Pgenft}; and
\needspace{2\baselineskip}
\item which, in a finite domain of state space containing~\EE, possesses \emph{exactly} the constructed \(M\)-mode slow \im~\cM, with \emph{exactly} the constructed evolution on~\cM. 
\end{enumerate}
\end{proposition}

\begin{proof} 
Under the stated preconditions, system~\eqref{Eempde} satisfies Assumption~3 of \cite{Bunder2018a} with the change \(\Re\lambda_j\mapsto |\lambda_j^{1/\alpha}|\) which does not change any of the subsequent construction in that paper, then the analysis there establishes the validity of \Xv-local analysis.
Under \cref{Aeig} a linear coordinate transform exists to convert~\cL\ to having orthogonal eigenvectors.
Then the theorems of \cite{Hochs2019} apply to give a spatially-local version of the results of \cref{Pgenabt} at each~\(\Xv\in\XX\) with two aspects:  firstly, change \(\Re\lambda_j\mapsto |\lambda_j^{1/\alpha}|\) which does not change any of the subsequent construction (such as the inequalities in Thm.~2.18 which still apply) provided additional non-resonance criteria are met in order to ensure various time integrals exist and are bounded;   and secondly, with one extra ingredient, namely that the spatial gradients give a \pde~\cite[(44)]{Bunder2018a} with a complicated remainder expression that is of the order~\(N+1\).
Hence the stated existence and invariance results hold for every locale~\(\Xv\in\XX\), and so the results hold globally in the \text{spatial domain~\XX.}
\hfill\end{proof}

Such a slow \im\ model is expected to be a guiding centre (\cref{Dim}), but nonlinearity and/or non-autonomous effects may cause a wide variety of unexpected behaviours near~\cM\ which are not present on~\cM.

However, a ``guiding centre'' is not the only justification for an \(M\)-mode, \(M\)-continuum \im~model.
Another justification for a slow manifold model is that it may be physically observable due to physical effects not encoded in the mathematics or computation.
For example, one such often uncoded physical mechanism is weak viscoelastic effects that damp out fast elastic waves.
Another often uncoded physical mechanism is that the fast waves may radiate out of the spatial domain~\XX\ of interest to leave the slow subspace within~\XX: for two examples, sound may radiate energy from a vibrating beam; and in the atmosphere fast inertial waves propagate up into the upper atmosphere where they break and dissipate, leaving the bulk of the atmosphere in the \text{quasi-geostrophic slow manifold.}

\paragraph{Scenarios modelling modulation of small-wavelength waves}
Two further scenarios for wave-like systems may be invoked.
A first scenario is when non-autonomous forcing, in~\(\gamma g\) of~\cref{Egenpde}, forces the system near some `fast' natural frequency of the cell problem, and when the forcing overcomes (weak) dissipation or radiation out of~\XX\ of the corresponding `fast' cell modes.
In this case those fast modes may be needed in a macroscale model \cite[e.g.,][]{Touze2006, Touze2021}.  
Two significant macroscale features are: firstly, the modulation of the fast mode; and secondly, the potential for nonlinear wave-wave resonance of the fast mode to drive mean `flow' in the slow modes.
An example of the latter is the large scale Stokes drift driven by relatively small scale ocean surface waves \cite[e.g.,][]{Craik2005}.
However, the existence and smoothness of such a model for nonlinear systems is highly problematic due to the likely plethora of nonlinear wave-wave (near) resonances \text{with other cell-modes.}

A second scenario is either when the initial conditions are that of a near uniform small-wavelength wave, or when a localised physical initial condition together with wave dispersion combine to cause a near uniform small-wavelength in the spatial domain~\XX\ of interest.
Then one wishes to focus purely on the one fast mode of the particular small wavelength mode, under the convenient but questionable assumption that other modes remain small.
One example would be the so-called \emph{high-frequency homogenisation} used to understand band gaps, Bloch waves, and Brillouin zones in heterogeneous material \cite[e.g.][]{Craster2010, Touboul2024}.
In this scenario one would construct an \im\ describing the large-scale modulation of the wave from a linear base of the one mode, and via embedding the system in the ensemble of all phase-shifts of the wave \cite[\S2.5]{Roberts2013a}. 
Generally, the resulting models are variants of the so-called nonlinear Schr\"odinger \pde~\cite[e.g.,][]{Whitham74, Totz2012}. 
The existence and smoothness of such a model is highly problematic due to potential nonlinear (near) resonances, but again backwards theory would assert there is a system close to that specified which exactly possesses the constructed~\im.
When the initial conditions possess multiple identifiable waves, then akin to quasi-periodic cases, embedding the system in the ensemble of all independent phase-shifts should systematically lead to a model expressed in interacting \text{nonlinear Schr\"odinger \pde{}s.}

\subsubsection{Fractional differential evolution in time}
\label{SSSfde}

Fractional derivative models are increasingly found in various engineering and science scenarios \cite[e.g.,][]{HongGuangSun2018}.
Here we address fractional time derivatives, \emph{not} fractional space derivatives.
In this case let the time differential operator~\alphaD\ for real fractional \(\alpha\in(0,2)\C \alpha\neq1\), be interpreted in the Caputo sense \cite[e.g.,][]{Gorenflo2008} and implicitly from initial time \(t=0\).
Let's reconsider the general solution~\cref{Emode} for the linearised cell-problem~\cref{EempdeX}.
From the general solution~\eqref{Egsmfde} derived in \ifJ Supplementary Material  \cite[Appendix~D,][]{Roberts2024a}\else\cref{Sgsmfde}\fi, the general solution to the \(m\)th-mode fractional differential equation (\fde)~\cref{Emode}, \(\alphaD a_m = \lambda_m a_m\)\,, is, in terms of parameter \(\mu_m:=(-\lambda_m)^{1/\alpha}\) and functions~\(\e_\alpha^{(k)}(t)\) defined by~\cref{Eealpha},
\begin{equation}
a_m(t)=\begin{cases}C_{m0}\e_\alpha^{(0)}(\mu_mt)
, &0<\alpha<1\,,
\\ C_{m0}\e_\alpha^{(0)}(\mu_mt),
+ C_{m1}\e_\alpha^{(-1)}(\mu_mt), & 1<\alpha<2\,.
\end{cases}
\label{Egsmfde}
\end{equation}
The free constants~\(C_{m0},C_{m1}\) happen to be proportional to the derivatives at time zero: \(C_{mk}\propto a_m^{(k)}(0^+)\).
Define the two functions \ifJ\cite[Appendix~D]{Roberts2024a}\else(\cref{Sgsmfde})\fi
\begin{align}&
\e_\alpha^{(0)}(t):=E_{\alpha,1}(-t^\alpha)
&&
\e_\alpha^{(-1)}(t):=tE_{\alpha,2}(-t^\alpha),
\label{Eealpha}
\end{align}  
in terms of the Mittag--Leffler function $E_{\alpha, \beta}(z) := \sum_{k=0}^\infty \frac{z^k}{\Gamma(\alpha k+\beta)}$ \cite[e.g.,][(A.1)]{Gorenflo2008}.

Arguably, the most important aspect of the general solution~\eqref{Egsmfde} is its behaviour for large time in the most common case of real negative eigenvalues~\(\lambda_m\).  
In this case~\(\mu_m=(-\lambda_m)^{1/\alpha}\) is real and positive and we need only consider the general solution for real positive arguments to~\(\e_\alpha^{(k)}\).
Since \(E_{\alpha,\beta}(z)\sim (-z)^{-1}/\Gamma(\beta-\alpha)\) as \(|z|\to\infty\) with \(|\arg(-z)|<\pi(1-\alpha/2)\),  \ifJ Supplementary Material \cite[Appendix~D,][]{Roberts2024a}\else\cref{Sgsmfde}\fi\ derives that \(\e_\alpha^{(k)}\), \cref{Eealpha}, generically decay algebraically to zero:
\begin{align}&
\e_\alpha^{(0)}(t)\sim \frac{t^{-\alpha}}{\Gamma(1-\alpha)}\,,
&&
\e_\alpha^{(-1)}(t)\sim \frac{t^{1-\alpha}}{\Gamma(2-\alpha)}\,,
&&\text{as }t\to+\infty\,.
\label{EfdeAsym2}
\end{align}
\begin{SCfigure}
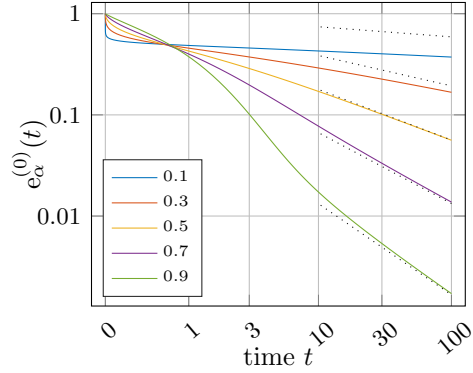
\centering
\caption{\label{FfdeSolns1}For the modal \fde~\eqref{Emode}, plot the general solution~\eqref{Egsmfde} component~\(\e_\alpha^{(0)}(t)\) for five \(\alpha\in(0,1)\) (as in legend).  
The time axis is quasi-log via an asinh scaling of the axis, and so clearly shows the long-time algebraic decay of this component in the general solution.  
The dotted lines are the large time asymptotes~\eqref{EfdeAsym2}. }
\inPlot{Figs/fdeSolns1}
\end{SCfigure}%
\cref{FfdeSolns1} plots the numerically computed solution component~\(\e_\alpha^{(0)}(t)\) for various \(0<\alpha<1\) with the large-time approximations~\eqref{EfdeAsym2}.
The quasi-log-log nature of the plot clearly exhibits the long-time algebraic decay of this component of the general solution~\cref{Egsmfde}, with both the decay and the approach to the asymptotes very slow for small~\(\alpha\).
For the larger case of \(\alpha=0.9\), \(\e_\alpha^{(0)}(t)\)~exhibits exponential-like decay for small time, say \(t<5\)---characteristic of the exponential decay for the \ode\ case of \(\alpha=1\)---before morphing to algebraic decay \text{for larger times, say \(t>10\)\,.}

\begin{figure}
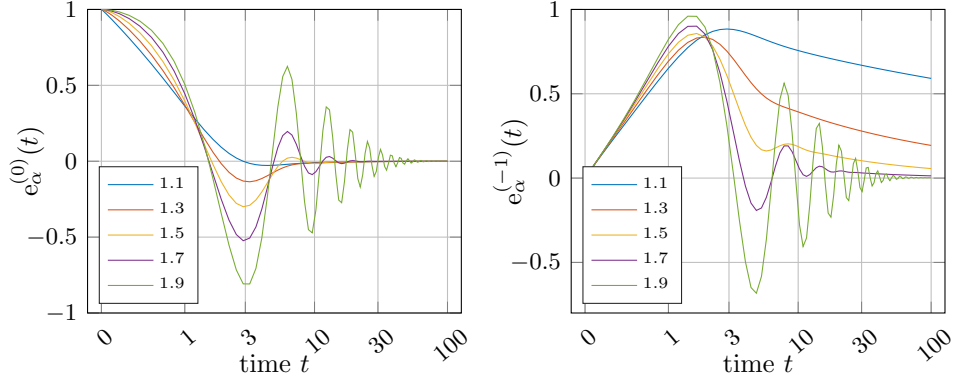
\centering
\caption{\label{FfdeSolns2}For the modal \fde~\eqref{Emode}, plots the two general solution~\eqref{Egsmfde} components for five \(\alpha\in(1,2)\) (see legend): (left)~\(\e_\alpha^{(0)}(t)\); (right)~\(\e_\alpha^{(-1)}(t)\).  The time axis is quasi-log via an asinh scaling of the axis.}
\begin{tabular}{@{}c@{}c@{}}
\inPlot{Figs/fdeSolns2}
&
\inPlot{Figs/fdeSolns3}
\end{tabular}
\end{figure}%
For the next range of exponents, \(1<\alpha<2\)\,, \cref{FfdeSolns2} plots the two solution components~\(\e_\alpha^{(0)}(t)\) and~\(\e_\alpha^{(-1)}(t)\)\,, \cref{Eealpha}.
These components oscillate some number of times---the number of oscillations increase as \(\alpha\to2\) (the wave case)---before at large time morphing into the eventual algebraic decay of the asymptotic~\eqref{EfdeAsym2}.

Now let's address what the above means for modelling and/or homogenisation of general fractional physical systems~\eqref{Egenpde}.
The clearest case is for linear autonomous systems~\eqref{Egenpde} and when the cell eigen-problem~\eqref{Eemeig} has a robust spectrum, such as obtained from a self-adjoint cell-operator~\cL.
Then the modal \fde{}s~\eqref{Emode} are robust under perturbations by macroscale spatial gradients.
Consequently, under such perturbations we expect each mode to still evolve roughly like \((-\lambda_m t)^{-\alpha-1+\alphahat}\) for large time, where \(\alphahat := \lceil\alpha\rceil ={}\)nearest integer\({}\geq\alpha\), and provided \(|\arg(-\lambda_m)|\leq\vartheta_{\max}<\pi(1-\alpha/2)\) for some angle bound~\(\vartheta_{\max}\).
Because of the eigenvalue factor~\(\lambda_m\), the modes with large~\(|\lambda_m|\) decay the quickest, albeit with the same exponent.  
The result is that we reasonably expect that the large time evolution is dominated by the modes with small~\(|\lambda_m|\).
Then, as in the exponential case of \cref{SSSsd}, we may reasonably construct and use a multi-modal, multi-continuum, homogenisation based upon choosing the \(M\)-modes with smallest~\(|\lambda_m|\).

The above argument reasonably justifies,  for general fractional~\(\alpha\), homogenisations like~\eqref{EE3Mevol,eqimphpde,Epde2,EEhcegManifold2}.

However, nonlinearities and non-autonomous forcing complicates the situation.
There is some extant theory of \im{}s for fractional differential systems \cite[e.g.,][]{Cong2016, Ma2015}, but none appears suitable to invoke for general systems~\eqref{Egenpde}.
I conjecture that adapting proposed backwards theory \cite[]{Roberts2018a, Hochs2019} would provide the most accessible route for supporting  \im\ homogenisation of nonlinear non-autonomous \fde\ systems.
However, it remains to be established whether or not a given \fde~system is generally `close' to a diffeomorphism of a constructible \fde\ system that also is in the \text{separated canonical form.}

\subsubsection{Improving spatial resolution}
\label{SSSisr}

\cref{SSSsd,SSSwls,SSSfde} focusses on the scenario when one chooses multi-modal multi-continuum homogenisations based upon selecting those modes with the longest lifetime or slowest evolution.
But there is significant interest in scenarios where the prime motivation is to better resolve spatial structures.
For two examples: \cite{Alavi2023} [p.2164] comment that ``Enriched continuum theories are required in such situations to capture the effect of spatially rapid fluctuations at the mesoscopic and macroscopic levels''; and \cite{Somnic2022} [p.8] wrote ``Classical continuum theory is not suitable when \ldots\ high strain gradients are observed in the domain''.
Hence this section discusses possibilities for choosing \im\ modes based upon the criterion of improving the spatial resolution irrespective of the time resolution.

Often the two criteria, spatial resolution and time resolution, are essentially equivalent: then the previous discussions apply.
This equivalence is generally the case for nonlinear systems because nonlinearity typically has many mode interactions which spread the dynamical energy among all modes, and so all modes are generally excited.

However, in linear systems (most of \cref{SoneDintro,Shceq,Selastic2d}) symmetries are more likely to be preserved.
Consider the \cref{EGsd} multi-modal multi-continua modelling of shear dispersion along a channel: \cite{Watt94b} found that the even modes across the channel interacted with the leading mean mode, but the odd modes did not. 
Consequently, choosing the leading two \emph{even} modes to form a bi-continua model (as listed in \cref{EGsd}) gives a model with much improved spatial resolution, when compared to that of the leading order model, but no improvement of the temporal resolution because the two-even-mode model neglects the gravest odd-mode (which has the slowest non-zero decay rate).
In general systems, quantitive estimates of the spatial resolution may be made for linear systems,  \cref{SSScsw} for example, by constructing to high-order in~\(\partial_x\), taking the Fourier transform which converts macroscale spatial derivatives to multiplication by the macroscale wavenumber~\(k\), that is, \(\partial_x \equiv \i k\)\,.
Then a Domb--Sykes plot \cite[e.g.]{Domb57, Hunter87} or Mercer--Roberts plot \cite[Appendix]{Mercer90} for series' in the wavenumber~\(k\), such as \cref{Fcas1dCs2}, predict the convergence limiting singularity in the complex \(k\)-plane, and hence the radius of convergence, say~\(k_*\).
The quantitative minimum resolved wavelength by the modelling is then deduced to be~\(2\pi/k_*\).
In contrast to traditional asymptotic homogenisations, which leave open the question of ``how small is small?" \cite[e.g.,][]{Somnic2022}, the \im~approach can empower us to quantitatively bound the \text{spatial resolution of a model.}

One way to view an \im\ for linear systems is as a perturbation from the cell-problem~\eqref{EempdeX} in small macroscale wavenumber~\(k\). 
A crucial characteristic of the modelling are the eigenvalues~\(\lambda(k)\) as a perturbative function of the wavenumber~\(k\).
Let's explore some fascinating aspects of such a multi-mode, multi-continuum, homogenisation via the mathematical structure of perturbative eigenvalue problems \cite[e.g.,][\S7.5]{Bender81}.

\begin{subequations}\label{EEGeig2}%
\begin{example}[simple spatial modelling]\label{EGeig2}
Consider modelling the macroscale evolution of the two interacting and differentially advecting components~\(u_j(t,x)\) governed by the coupled \pde{}s \cite[a variation of Exercise~7.4]{Roberts2014a}
\begin{align}&
u_{1t}=-u_{1x}+(-u_1+u_2),&&
u_{2t}=+u_{2x}+(+u_1-u_2).
\label{EGeig2pde}
\end{align}
Invoke the spatial Fourier transform by seeking solutions \(u_j=c_j\e^{\i kx+\lambda t}\).   
Then the system~\cref{EGeig2pde} reduces to the eigen-problem
\begin{equation}
\lambda\begin{bmatrix} c_1\\c_2 \end{bmatrix}
=\begin{bmatrix} -\i k-1 &1 \\1 & +\i k-1 \end{bmatrix}
\begin{bmatrix} c_1\\c_2 \end{bmatrix}.
\label{EGeig2p}
\end{equation}
The eigenvalues are thus
\begin{equation}
\lambda=-1\pm\sqrt{1-k^2}
\sim-1\pm(1-\tfrac12k^2-\tfrac18k^4-\tfrac1{16}k^6-\cdots).
\label{EGeig2e}
\end{equation}
Observe that the eigenvalues of this system are just two Riemann sheets of the one surface in the complex \(k\)-plane.
The slow \im~model corresponds to the small~\(\lambda\) sheet of \(\lambda =-1+\sqrt{1-k^2}
\sim -\tfrac12k^2 -\tfrac18k^4 -\tfrac1{16}k^6 -\cdots\).
Taking the inverse Fourier transform of this sheet gives the slow \im~model, in terms of a macroscale~\(U_0(t,x)\), to be some chosen truncation of
\begin{equation}
U_{0t} =\tfrac12U_{0xx} -\tfrac18U_{0xxxx} +\tfrac1{16}U_{0xxxxxx} +\cdots.
\label{EGeig2U}
\end{equation}
The question under discussion is the following: for derived models such as this~\cref{EGeig2U}, what is a quantitive limit on the spatial resolution?
That is, how small is `small wavenumber'?
As an example of the general case, the Fourier space~\cref{EGeig2e} shows that here the limit is due to square-root singularities\footnote{These are `physical' singularities as they happen to occur at real wavenumber.} at wavenumber \(k=\pm1\) (branch points), and so the series~\cref{EGeig2U} does not converge for \(|k|>1\)\,. 
Hence, here the sharp lower bound for spatial resolution of the model~\cref{EGeig2U} is that of spatial structures with wavelength\({}>2\pi\).
\qed\end{example}
\end{subequations}

The crucial property of eigen-problems for a given \(n\times n\) matrix which is a function of a parameter~\(k\) is that its eigenvalue function~\(\lambda(k)\) generally forms a single \(n\)-sheeted (Riemann) surface in complex~\(k\) \cite[e.g.,][p.350]{Bender81}.
Such an \(n\)-sheeted surface possesses various branch points at complex~\(k\) that are necessary to connect the sheets together into a unified whole (e.g., as in~\cref{EGeig2e}).
In some cases, often due to symmetries, the \(n\)-sheeted surface does partition into several disjoint sheeted surfaces.
These properties also generalise to many \pde\ operators of \text{interest in applications.}

The branch points in the Riemann sheets are singularities that limit the convergence of the perturbative expansion for~\(\lambda(k)\) in wavenumber~\(k\), and hence limit the spatial resolution of the macroscale homogenisation (e.g., \cref{EEGeig2}).
For example, the tri-continuum high-order homogenisation of \cref{SSScsw} shows the eigenvalue sheet of the homogenisation, for small heterogeneity~\(a\), has pole singularities at real wavenumbers \(k=\pm3/2\), pole singularities that limit the resolution of the homogenisation \text{in physical space.}

\paragraph{Regularisation is akin to Pad\'e approximation}
Recall the example regularisation of \cref{Shcegsmm}.
Let's recast the example gradient models~\cref{EEhcegManifold1} in Fourier space, effectively replace~\(U_0(t,x)\) by~\(\tU(t)\e^{\i(k_1x+k_2y)}\) for wavenumber vector \(k:=(k_1,k_2)\).
The derived and the two regularised models become, respectively (to two decimal digits)
\begin{equation*}
\alphaD{\tU_0}= \tU_0\cdot
\begin{cases}
-.53k_1^2  -.95k_2^2
+(.34k_1^4 -.14k_1^2k_2^2 +.015k_2^2)
-(-.90k_1^6 
\\\quad{}
+.28k_1^4k_2^2 -.0021k_1^2k_2^4 -.0038k_2^6)
+\Ord{|k|^7},
\quad\text{for \cref{EhcegManifold1U}};
\\[1ex]
\text{\large$\frac{ -.54k_1^2 -.95k_2^2 -.74k_1^2k_2^2 }
 {1+.63k_1^2+.02k_2^2}$} +\Ord{|k|^5},
\quad\text{for \cref{EhcegManifold1Uc}};
\\[2ex]
\text{\large$\frac{ -.54k_1^2 -.95k_2^2 -.74k_1^2k_2^2 -2.33k_1^4k_2^2}
 {1+.63k_1^2+.02k_2^2 +2.07k_1^4 +.007k_1^2k_2^2  +.004k_2^4}$} +\Ord{|k|^7},
\quad\text{for \cref{EhcegManifold1Ud}}.
\end{cases}
\end{equation*}
Evidently such regularised models in Fourier space are equivalent to rational function approximation of~\(\lambda(k)\), as one might obtain by (multivariate) Pad\'e approximation \cite[e.g.,][]{Cuyt86}.
Because of the well-established ability of Pad\'e approximation to improve the predictions of analytic functions \cite[e.g.,][]{Stahl97}, especially away from branch points, we generally expect such rational representations of~\(\lambda(k)\) to be useful over a wider range of wavenumbers than the series representation.
Hence, we generally expect corresponding regularised \pde\ models~\cref{EhcegManifold1Uc,EhcegManifold1Ud} to have improved spatial resolution compared to the constructed \pde~\cref{EhcegManifold1U}---as discussed for plasticity in the \text{survey by \cite{Bazant2002}.}

But such regularisation has two deficiencies. 
First, it cannot account for the dynamics associated with those Riemann sheets which, often by symmetry, happen to be disconnected from the sheet of the constructed slow manifold.
Shear dispersion in a channel or pipe are examples \cite[e.g.,][]{Watt94b}: the classic one-mode continuum model only `knows about' axisymmetric modes and their dynamics, and consequently there is no regularisation that can predict non-axisymmetric dynamics.
Second, regularisation cannot get past a branch-point located at or near physical values of real wavenumber~\(k\).
To remedy such deficiencies, we must instead add micromorphic modes \text{to the modelling.}

\paragraph{Which branch-points limit a chosen multi-continuum models?}
For linear systems, an \(M\)-mode, \(M\)-continuum homogenisation constructs an homogenisation whose eigenvalues~\(\lambda\) form an \(M\)-sheeted Riemann surface as an analytic function, say~\(\lambda^M(k)\), of wavenumber~\(k\) in the complex \(k\)-plane.
The original system also has a many-sheeted analytic eigenvalue function~\(\lambda(k)\).
The \(M\)-mode homogenisation approximates~\(\lambda(k)\) by analysing effects and interactions of the \(M\)~chosen modes---modes chosen from the eigenvalue spectrum of the cell-problem obtained at at zero wavenumber (\cref{SSSlbim}).
The systematic construction guarantees that the perturbation series for~\(\lambda^M(k)\) is identical with the series for the corresponding \(M\)~Riemann sheets of the original~\(\lambda(k)\).
Consequently, we expect an \(M\)-mode homogenisation to encompass smoothly all branch-point singularities in the original~\(\lambda(k)\) that arise through interactions among the chosen \(M\)-modes (whether physical singularities at real~\(k\), or unphysical at complex valued~\(k\)).
Branch-point singularities in~\(\lambda(k)\) that occur, physically or not, due to interactions between the \(M\)-chosen sheets and any/all of the \emph{other unchosen} sheets cannot be captured in the \(M\)-mode homogenisation, and so it is these singularities that generally limit the radius of convergence of the \(M\)-continuum homogenisation.
It is this restriction in wavenumber~\(k\) that then translates to a limit on the spatial macroscales potentially resolvable in a chosen \(M\)-mode, \(M\)-continuum homogenisation. 

Consequently, to obtain the best macroscale spatial resolution in a multi-mode multi-continua homogenisation, without considering temporal resolution, one needs to select the \(M\)-modes that interact and resolve the branch-points \emph{nearest} to wavenumber zero of the eigenvalue~\(\lambda(k)\) Riemann surface among the sheets of small eigenvalue~\(\lambda\).  
However, such selection is usually difficult because commonly the near branch-point singularities are at unphysical complex-valued wavenumbers for no known physical reason.
Unless there are known physical symmetries, the \emph{best practical general guide to improve macroscale spatial resolution appears to be to select those modes with the longest time scale}, as discussed \text{in \cref{SSSsd,SSSwls,SSSfde}.}

\subsection{Construct a chosen invariant manifold multi-continuum homogenised model}
\label{SSSccimmch}

The construction of multi-continuum homogenisation of, potentially nonlinear, systems in multiple large spatial dimensions relies on theory proven by \cite{Bunder2018a}, which in turn rests on general theory by \cite{Aulbach2000, Potzsche2006, Hochs2019}. 
One of the crucial theoretical results is the direct connection that if a derived approximation satisfies the embedding \pde~\cref{Eempde} to a residual of~\Ord{\gradx^{N+1}}, then the corresponding homogenisation is correct to an error~\Ord{\gradx^{N+1}} \cite[\S3]{Bunder2018a}. 
In multi-D space, a practical procedure to construct models to any chosen order of residual are, in this section, shown to be encompassed by variations to the iterative \cref{Pconstruct}, as applied for example to the homogenisation of 2-D heterogeneous elasticity (\cref{Selastic2d}).

A key step is to determine corrections for any given approximation.
For \(M\)~selected eigenvector modes (\cref{Simmmh}), choose a definition of the vector of local amplitudes~\(\Uv(t,\fx):=(U_0,\ldots,U_{M-1})\).
We seek an \im\ of the embedding \pde~\cref{Eempde} in the form \(\fu=v(\Uv,t,\fx,\thetav)\) such that \(\alphaD \Uv=G(\Uv,t,\fx)\) where the right-hand side dependence upon~\Uv\ implicitly involves its general gradients~\(\Uv_{\x_i}\C \Uv_{\x_i\x_j}\C\ldots\), and the explicit \((t,\fx)\) dependence are the slow, macroscale, variations arising from macroscale functional graduations in the problem~\cref{Eempde}.
For every given approximations~\(\tv,\tG\) to the \im~\(v,G\), define~\(\Res(\tv,\tG)\) to be the residual of the embedding \pde~\cref{Eempde}.  

\begin{lemma}\label{LgenHomo}
Compute corrections~\(v',G'\) to the approximations~\(\tv,\tG\) by solving a variant of the cell problem~\eqref{EempdeX}, linearised and forced by the residual, called the \emph{homological equation} and taking various forms:\footcite[e.g.,][\protect\cref{EEgenHomologic} is a more explicit form than, e.g., (4.3) of Potzsche \& Rasmussen]{Potzsche2006, Roberts2014a, Siettos2021, Martin2022}
\begin{subequations}\label{EEgenHomologic}%
\needspace{4\baselineskip}
\begin{itemize}
\item for \(\alpha=1\)\,, and generalising~\cref{E1dHomologic},
\begin{equation}
\cL v'
-\D t{v'}
-\sum_{m=0}^{M-1}\D{U_m}{v'}(\lambda_mU_m+\zeta_mU_{m+1})
-\sum_{m=0}^{M-1}v_mG'_m
=\Res(\tv,\tG);
\label{EgenHomologic1}
\end{equation}

\needspace{5\baselineskip}
\item for \(\alpha=2\),%
\footnote{In the case of \(\alpha=2\) and nonlinear systems, then both the homological equation and its solution becomes very complicated.   
In this case it is vastly simpler to recast as a system first-order in time, and then use~\cref{EgenHomologic1}. }
\begin{align}&
\cL v'
-\DD t{v'}
-\sum_{m=0}^{M-1}\D{U_m}{v'}(\lambda_mU_m+\zeta_mU_{m+1})
-2\sum_{m=0}^{M-1}\Dx{U_m}t{v'}\D t{U_m}
\nonumber\\&\qquad{}
+\sum_{m,n=0}^{M-1} \Dx{U_m}{U_n}{v'}\D t{U_m}\D t{U_n}
-\sumM v_mG'_m
=\Res(\tv,\tG).
\label{EgenHomologic2}
\end{align}

\needspace{6\baselineskip}
\item for the case of a slow \im~model with a significant spectral gap (\cref{Dgap}),
\quad and when the time dependence is correspondingly slow in the \pde~(in the right-hand side function~\(g\) of \cref{Egenpde,Eempde}), in practice one may use the following simplified form (at the cost of requiring more iterations):
\begin{equation}
\cL v'
-\sumM v_mG'_m
=\Res(\tv,\tG)
\quad(\text{for both }\alpha=1,2).
\label{EgenHomoSimp}
\end{equation} 

\end{itemize}%
\end{subequations}%
\end{lemma}
In \cref{LgenHomo}, as in \cref{SSScmmm}, the factors \((\D{U_m}{v'})U_n\) are Fre\'chet derivatives to be interpreted in the Calculus of Variations sense such that here it represents \(v'_{U_m}U_n +\sum_iv'_{U_{m\x_i}}U_{n\x_i} +\sum_{i,j}v'_{U_{m\x_i\x_j}}U_{n\x_i\x_j} +\cdots\) where these subscript-derivatives of~\(v'\) are done with respect to the subscript symbol \cite[]{Roberts88a}. 
Also, the derivative \(\Dn t\alpha{v'}\) denotes the partial time-derivative keeping all~\(U_m\) constant (as opposed to \(\alphaD\) which denotes the time derivative including the time evolution of all~\(U_m\)).
Recall that non-zero~\(\zeta_m\), appearing in~\cref{EgenHomologic1,EgenHomologic2},  arises in cases with generalised eigenvectors (\cref{Aeig}).

Differences in the homological equations~\eqref{EEgenHomologic} from the usually used cell-problems mostly lie in the terms with factors of~\(\D{U_m}{v'}\).
These terms arise in this systematic \im~framework through accounting for physical out-of-equilibrium effects, both explicitly in the left-hand side and also implicitly encoded within the right-hand side~\(\Res(\tv,\tGv)\).
Then \cref{Pconstruct}, with Step~8(b) using an homological equation from~\cref{EEgenHomologic}, applies to systematically construct \im~models in these \text{general scenarios.}

\begin{proof}[Proof of \cref{LgenHomo}] 
Let's consider the embedding \pde~\eqref{Eempde} in the form \(\alphaD \fu = \cL\fu +\Fv(\fu)\) for cell-operator~\cL\ defined in~\cref{Eemeig}, and \Fv~denotes all `perturbative' effects of~\eqref{Eempde}.
The complicated part of the embedding \pde~\eqref{Eempde} is the time derivative~\(\alphaD\fu\). 
Consider \(\alpha=1,2\) in turn, letting \(\D t{}\C\D{U_m}{}\) denote partial (Fre\'chet) derivatives keeping constant the other variables in the list~\(t,U_0,\ldots,U_{M-1}\)\,.

For \(\alpha=1\), and using \(\partial_t{U_m}=G_m\), \begin{equation*}
 \partial_tv
=\D tv +\sumM \D{U_m}v\D t{U_m}
=\D tv +\sumM \D{U_m}vG_m\,.
\end{equation*}
Given approximations~\(\tv,\tGv\) to an \im~model, seek an improved model \(v=\tv+v'\) and \(\Gv=\tGv+\Gv'\) for small (dashed) corrections.
Substitute into the governing equations to find the residual
\begin{align*}
\Res(v,\Gv)&
= \D t\tv+\D t{v'} +\sumM \D{U_m}{(\tv+v')}(\tG_m+G'_m) -\cL(\tv+v')-\Fv(\tv+v')
\\&\qquad\text{(upon neglecting products of small quantities)}
\\&
\approx  \D t\tv +\sumM \D{U_m}{\tv}\tG_m -\cL\tv-\Fv(\tv) 
\notJbreak{}
+\D t{v'} +\sumM \left[\D{U_m}{v'}\tG_m +\D{U_m}{\tv}G'_m \right] -\cL v'
\\&\qquad\text{(using leading order \(\tv,\tGv\) when multiplied by small corrections)}
\\&
\approx  \Res(\tv,\tGv) 
+\D t{v'} +\sumM \left[\D{U_m}{v'}(\lambda_mU_m+\zeta_mU_{m+1}) +v_mG'_m \right] -\cL v'.
\end{align*}
To seek corrections that reduce the residual, setting this right-hand expression to zero leads to the homological equation~\cref{EgenHomologic1}.

For \(\alpha=2\), and using \(\partial_t^2{U_m}=G_m\)\,, \begin{align*}&
 \partial_t^2v
=\partial_t\left[\D tv +\sumM \D{U_m}v\D t{U_m}\right]
\\&
= \DD tv   
+2\sumM \Dx t{U_m}v\D t{U_m}
+\sum_{m,n=0}^{M-1} \Dx{U_m}{U_n}v\D t{U_m}\D t{U_n}
+\sumM \D{U_m}vG_m\,.
\end{align*}
Akin to the \(\alpha=1\) substitute \(v=\tv+v'\) and \(\Gv=\tGv+\Gv'\) into the governing equations to find the residual
\begin{align*}
\Res(v,\Gv)&
= \D t\tv+\D t{v'} +2\sumM \Dx t{U_m}{(\tv+v')}\D t{U_m}
+\sum_{m,n=0}^{M-1} \Dx{U_m}{U_n}{(\tv+v')}\D t{U_m}\D t{U_n}
\\&\quad{}
+\sumM \D{U_m}{(\tv+v')}(\tG_m+G'_m) 
-\cL(\tv+v')-\Fv(\tv+v')
\\&\qquad\text{(then neglecting products of small quantities)}
\\&
\approx  \D t\tv +2\sumM \Dx t{U_m}{\tv}\D t{U_m} 
+\sum_{m,n=0}^{M-1} \Dx{U_m}{U_n}{\tv}\D t{U_m}\D t{U_n}
\notJbreak{}
+\sumM \D{U_m}{\tv}\tG_m -\cL\tv-\Fv(\tv) 
\\&\quad{}
+\D t{v'} +2\sumM \Dx t{U_m}{v'}\partial_t{U_m} 
+\sum_{m,n=0}^{M-1} \Dx{U_m}{U_n}{v'}\D t{U_m}\D t{U_n}
\\&\quad{}
+\sumM \left[\D{U_m}{v'}\tG_m +\D{U_m}{\tv}G'_m \right] -\cL v'
\\&\qquad\text{(using leading order \(\tv,\tGv\) when multiplied by small corrections)}
\\&
\approx  \Res(\tv,\tGv) 
+\D t{v'} +2\sumM \Dx t{U_m}{v'}\D t{U_m} 
+\sum_{m,n=0}^{M-1} \Dx{U_m}{U_n}{v'}\D t{U_m}\D t{U_n}
\\&\quad{}
+\sumM \left[\D{U_m}{v'}(\lambda_mU_m+\zeta_mU_{m+1}) +v_mG'_m \right] 
-\cL v'.
\end{align*}
To seek corrections that reduce the residual, setting this right-hand expression to zero leads to the homological equation~\cref{EgenHomologic2} (in linear problems the sum \(\sum_{m,n=0}^{M-1}\) is absent as \(\Dx{U_m}{U_n}{v}=0\)).

Three practical simplifications are often feasible for solving the homological equations~\eqref{EgenHomologic1,EgenHomologic2}, and these simplifications lead to~\cref{EgenHomoSimp} which is much easier to solve, and vastly quicker for humans to code solution methods.
\begin{enumerate}
\item Neglect the terms in~\(\zeta_m\) arising from generalised eigenvectors.
This neglect generally increases the number of required iterations by no more than a factor equal to the largest multiplicity of the chosen~\(M\) eigenvalues. 

\item When constructing slow \im~models to errors of order~\(p\) in nonlinearity, provided the spectral gap~\((\beta_L,\beta_\ell )\) (\cref{Fnfseval,Fnfreval0,Dgap}) is big enough so that the ratio \(\beta_\ell /\beta_L>p\)\,, then one may neglect the term in~\(\lambda_m\) as well.
Then each iteration reduces the numerical error in \im~coefficients by a factor of at most~\(p\beta_L/\beta_\ell<1\).

\item When the explicit time variation in the right-hand side~\(g\) is smooth over macroscale times, then the time-operator terms in~\(\Dn t\alpha{v'}\) may be neglected: the time variations are then accounted in the updates via operators~\({\alphaD}^n\) for various~\(n\) up to some chosen order \cite[e.g.,][]{Mercer90}.
The next \cref{SSiftf} discusses cases involving rapid microscale time variations in~\(g\).
\end{enumerate}
\hfill\end{proof}

The iterative construction terminates (Step~8 of \cref{Pconstruct}) when the zero-test of the residual is that all coefficients are smaller than some set small numerical threshold such as~\(10^{-7}\).
The elasticity homogenisation \ifJ of Supplementary Code \cite[Appendix~C,][]{Roberts2024a}\else code of \cref{cas2d}\fi\ generally invokes two of these simplifications, and so \text{solves~\eqref{EgenHomoSimp} for updates.}

\subsection{Intricacies of fast-time fluctuations}
\label{SSiftf}

Consider the scenario when the explicit time variations in the right-hand side of \pde~\cref{Egenpde} are on the same or faster timescale than that of modes chosen to be in the \im~model: the time variations could be deterministic or stochastic.
The underpinning dynamical systems theory is non-autonomous and so still applies: the \im{}s~\cM\ still exist and may be emergent \cite[e.g., Prop.~3.6 of ][Thm.~2.18]{Aulbach2000, Hochs2019}; and~\cM\ may be constructed to arbitrary order \cite[e.g.,][]{Potzsche2006}.

However, with explicit fast-time fluctuations the necessary algebraic details of approximating~\cM\ typically undergo a `combinatorial explosion' due to a rapidly increasing number of mode interactions.  
For example, see the finite-dimensional stochastic case of \cite{Roberts05c}.
Due to this `explosion' I only address systems written as first-order in time (\(\alpha=1\)): rewrite second-order in time systems as first-order.
In the class of systems where the residual of the governing \pde{}s contain explicit fast-time dependence, then in the homological equations~\eqref{EgenHomologic1} one must \emph{not} neglect~\(-\D t{v'}\) on the left-hand side. 
Crucial to the interpretation of the homological equations~\eqref{EgenHomologic1} is that this time partial derivative operator is done keeping constant both space~\xv,\zv, and also constant are all~\(U_m\), \text{and their gradients.}

With rapid time variation in the right-hand side, we generally need to solve~\eqref{EgenHomologic1} for every mode~\(n\) \cite[e.g.,][\S2 and \S3.2 resp.]{Roberts05c, Roberts06k}. 
Define mode correction \(v'_n:=\left<w_n,v'\right>\), and define residual component \(\Res_n(t):=\big<w_n,\Res(\tv,\tG)\big>\), and then take \(\left<w_n,\cdot\right>\) of~\eqref{EgenHomologic1} (omitting `off-diagonal'~\(\zeta_m\) as in~\eqref{EgenHomoSimp}) to require that, for every mode \(n=0,1,2,\ldots\)\,,
\begin{equation}
-\D t{v'_n}
+\lambda_nv'_n
-\sum_{m=0}^{M-1}\D{U_m}{v'_n}\lambda_mU_m
-G'_n =\Res_n(t)
\label{EmodeHomologic}
\end{equation}
(where \(G'_n=0\) for every \(n\geq M\)).
Almost always we express the residual as a polynomial in the amplitudes~\(U_m\) so that the mode residual is of the form of a sum \(\sum_\qv \Res_{n,\qv}(t) \Uv^\qv\) where \(\Uv^\qv := U_0^{q_0} \cdots U_{M-1}^{q_{M-1}}\).
Then we seek solutions of~\cref{EmodeHomologic} as the sum \(v'_n=\sum_\qv v'_{n,\qv}(t) \Uv^\qv\) and \(G'_n=\sum_\qv G'_{n,\qv}(t) \Uv^\qv\)\,, which requires
\begin{equation}
-\D t{v'_{n,\qv}}
+\mu_{n,\qv}v'_{n,\qv}
-G'_{n,\qv} =\Res_{n,\qv}(t),
\quad\text{where } 
\mu_{n,\qv}:=\lambda_n-\sum_{m=0}^{M-1} q_m\lambda_m\,.
\label{EmodeTerm}
\end{equation}
Such constant coefficient differential equations in time for mode components are sometimes straightforward to solve.
For example, in the simplest scenario of sinusoidal forcing, the solutions of the modal updates~\eqref{EmodeHomologic} are correspondingly sinusoidal \cite[e.g.,][]{Touze2006, Touze2021, Roberts2012a},  with resonant terms (otherwise secular-generating) being assigned to~\(G'_{n,\qv}\) (provided \(n<M\)).

In scenarios where the fast-time fluctuations in the governing \pde{}s~\eqref{Egenpde}, and hence in the mode residuals~\(\Res_n(t)\), are stochastic or as yet unknown control or other forcing, then dissipative and wave-like scenarios are quite different.  
In the dissipative case, \cref{SSSsd}, the solutions of modal updates~\eqref{EmodeHomologic} involve exponentially-decaying history integrals of the fast-time fluctuations \cite[e.g.,][\S3.3]{Arnold98, Roberts05c, Roberts06k}. 
Such history integrals are: 
akin to those arising in the Mori--Zwanzig formalism \cite[e.g.,][]{Chorin2006, Stinis05};
and are also another source of history dependence \cite[e.g.,][]{Eggersmann2019} frequently unrecognised in modelling.
The rate of exponential decay depends upon the mode number~\(n\) and the powers of all~\(U_m\) in each term, so a combinatorial explosion is usual.
In wave-like systems (\cref{SSSwls})  the solutions of modal updates~\eqref{EmodeHomologic} are much more delicate, and so maintaining accuracy requires considerable care \cite[e.g.,][\S5]{Roberts06k}.

Nonetheless, many non-autonomous systems with fast-time fluctuations are able, via \cref{Pconstruct} with~\cref{EmodeHomologic}, to be homogenised within this framework.


\section{Example in 2D elasticity homogenisation}
\label{Selastic2d}

This section applies the framework of \cref{Sgentheory} to rigorously construct and support a tri-continuum, three-mode, homogenisation of \text{2-D} elasticity for an example heterogeneous Young's modulus.
To potentially analyse quite general elastic heterogeneities, we describe here a combined  algebraic-numerical approach.
The sub-cell problems, solving homological equations~\cref{EEgenHomologic}, are computed numerically.
The macroscale variations themselves, and their influence of sub-cell structures, are described algebraically in terms of macroscale gradients~\(U_{mx}\C U_{my}\C U_{mxx}\C \ldots\), and other parameters.
Computer algebra holds the algebraic macroscale expressions united with coefficients that are microscale numerical arrays of 2-D sub-cell structures.

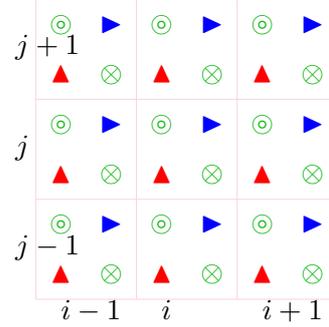
\begin{SCfigure}[1.2]
\centering
\caption{\label{Fmicrogrid}a small part of the microscale grid, of spacing~\(\delta x\) and~\(\delta y\), used to code 2-D elasticity.  The grid is staggered on the microscale: \uSym, horizontal displacements and velocities; \vSym, vertical displacements and velocities; \oSym, \xSym, components of strain and stress tensor~\eqref{EEstress}.
}
\setlength{\unitlength}{4ex}
\def\N{3}\def\Np{4}\def\NN{6}
\def\M{3}\def\Mp{4}\def\MM{6}
\begin{picture}(\NN.5,\MM.5)
\def\symBox#1{\makebox(0,0){#1}}
\put(1,0){
{\color{magenta!20}
    \multiput(-0.5,0.5)(0,2){\Mp}{\line(1,0){\NN}}
    \multiput(-0.5,0.5)(2,0){\Np}{\line(0,1){\MM}} 
    }
\multiput(1,2)(2,0){\N}{\multiput(0,0)(0,2){\M}{\uSym}}
\multiput(0,1)(2,0){\N}{\multiput(0,0)(0,2){\M}{\vSym}}
\multiput(0,2)(2,0){\N}{\multiput(0,0)(0,2){\M}{\oSym}}
\multiput(1,1)(2,0){\N}{\multiput(0,0)(0,2){\M}{\xSym}}
\setcounter{i}0
\multiput(0,0.1)(2,0){\N}{%
    $i\ifcase\value{i}-1\or\or+1\fi$\stepcounter{i}}
\setcounter{i}0
\multiput(-0.9,1.4)(0,2){\M}{%
    $j\ifcase\value{i}-1\or\or+1\fi$\stepcounter{i}}
}
\end{picture}
\end{SCfigure}

We adopt a simple robust microscale discretisation of the equations for heterogeneous isotropic {2-D}~elasticity.  
On a staggered microscale \(xy\)-grid of spacing~\(\delta x\) and~\(\delta y\), \cref{Fmicrogrid} shows a fragment, define the material displacements: 
\uSym,~horizontal~\(u_{ij}(t)\); 
\vSym,~vertical~\(v_{ij}(t)\).
The adopted microgrid elasticity uses centred finite differences, \(\delta_i\)~and~\(\delta_j\), to compute stresses at the shown staggered microscale \text{grid-points (\cref{Fmicrogrid}):}
\begin{subequations}\label{EEstress}%
\begin{align}&
\text{\xSym}\quad \sigma_{ij}^{xy}
:=\mu_{ij}\big[\delta_iv_{ij}/\delta x 
+\delta_ju_{ij}/\delta y\big];
\\&
\text{\oSym}\quad \sigma_{ij}^{xx}
:=(\lambda_{ij}+2\mu_{ij}){\delta_iu_{ij}}/{\delta x} 
+\lambda_{ij}{\delta_jv_{ij}}/{\delta y}\,;
\\&
\text{\oSym}\quad \sigma_{ij}^{yy}
:=\lambda_{ij}{\delta_iu_{ij}}/{\delta x} 
+(\lambda_{ij}+2\mu_{ij}){\delta_jv_{ij}}/{\delta y}\,,
\end{align}
\end{subequations}
where \(\lambda_{ij},\mu_{ij}\) denote heterogeneous Lam\'e parameters.
Then centred finite differences compute the following (non-dimensional) acceleration \ode{}s
\begin{subequations}\label{EEaccel}%
\begin{align}
&\text{\uSym}\quad 
 \partial_t^2 u_{ij}={\delta_i\sigma_{ij}^{xx}}/{\delta x} 
+{\delta_j\sigma_{ij}^{xy}}/{\delta y}\,,
\\&
\text{\vSym}\quad  
 \partial_t^2 v_{ij}={\delta_i\sigma_{ij}^{xy}}/{\delta x} 
+{\delta_j\sigma_{ij}^{yy}}/{\delta y}\,.
\end{align}
\end{subequations}
The Lam\'e parameters appearing in the stresses~\eqref{EEstress} are defined as
\begin{align}&
\lambda_{ij}:=\frac{\nu_{ij} E_{ij}}{(1+\nu_{ij})(1-2\nu_{ij})}\,,
&&
\mu_{ij}:=\frac{E_{ij}}{2(1+\nu_{ij})}\,,
\label{ELame}
\end{align}
in terms of heterogeneous Young's modulus~\(E_{ij}\) and Poisson ratio~\(\nu_{ij}\).

\subsection{Example cell problem}
\label{SSecp}

Suppose the microscale heterogeneity is reflected in the Lam\'e parameters~\(\lambda_{ij}\) and~\(\mu_{ij}\) being \(n_x\)-periodic in~\(i\) and \(n_y\)-periodic in~\(j\).  
That is, the material is \(\ell_x:=n_x\delta x\) periodic in the \(x\)-direction and \(\ell_y:=n_y\delta y\) periodic in the \(y\)-direction. 

We choose to non-dimensionlise on the microscale cell size, choosing the length scale so that each cell is \(1\)-periodic in~\(x,y\).
\begin{SCfigure}
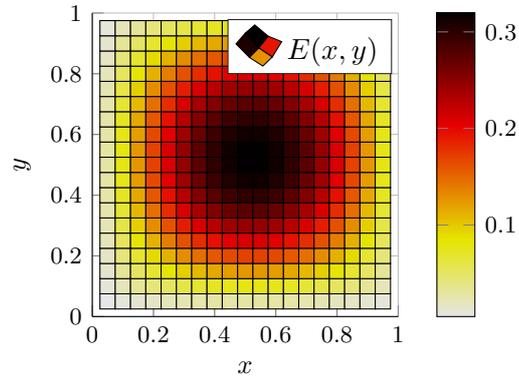
\centering
\caption{\label{Fspectra6n10E}example of one period, one cell, of the \text{2-D} microscale heterogeneity in the Young's modulus \(E = \big[0.01+|\sin(\pi x)\sin(\pi y)|\big]/\pi\).
The non-dimensionalised microscale periods are \(\ell_x=\ell_y=1\).
Data is plotted on the \(2n_x\times2n_y\) microscale staggered grid points. 
}
\inPlot{Figs/spectra6n10E}
\end{SCfigure}%
\cref{Fspectra6n10E} plots one such non-dimensional microscale cell of the example Young's modulus used for the case reported here with \(n_x=n_y=10\).
The material is made up of these cells repeating indefinitely in the \(xy\)-plane.
This microscale heterogeneity is similar to the \S3--5 example of \cite{Sarhil2024} (diagrams phase-shifted in each cell).
This heterogeneity is also similar to the `spinning top' metamaterials explored in \S7 of \cite{Milton2007}: but, in contrast, whereas Milton \& Willis [p.874] ``\emph{believe} the above formulation should be a good approximation'', our approach provides a rigorously proven homogenised model (\cref{Sgentheory}).
The homogenisation describes the evolution of elastic waves which have wavelengths significantly larger than one; that is, in this non-dimensionalisation the macroscale wavenumbers are those significantly \text{less than~\(2\pi\). }

In principle we should now embed the system~\cref{EEstress,EEaccel} and then linearise to obtain and analyse its cell-problem~\cref{EempdeX} with its periodic boundary conditions as justified by \cref{Spsem2dh}.
However, in practice let's obtain a preliminary physical understanding of general cell solutions by exploring the eigenvalues and eigenvectors of the cell eigen-problem~\eqref{Eemeig} corresponding directly to the physical~\cref{EEstress,EEaccel}.

\begin{figure}
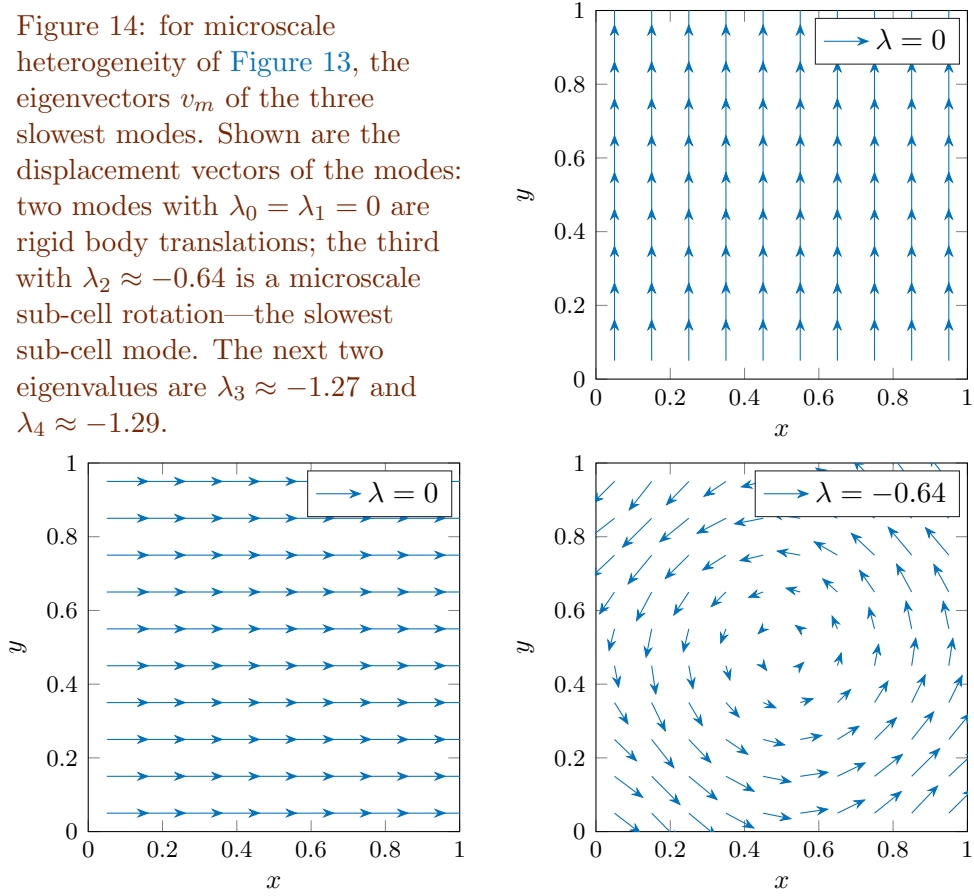
\centering
\def\extraAxisOptions{small,axis equal image,width=75mm}
\begin{tabular}{@{}cc@{}}
\parbox{0.46\linewidth}{%
\caption{\label{Fspectra6n10ks}\raggedright%
for microscale heterogeneity of \cref{Fspectra6n10E}, the eigenvectors~\(v_m\) of the three slowest modes.
Shown are the displacement vectors of the modes: two modes with \(\lambda_0=\lambda_1=0\) are rigid body translations; the third with \(\lambda_2 \approx -0.64\) is a microscale sub-cell rotation---the slowest sub-cell mode.  
The next two eigenvalues are \(\lambda_3 \approx -1.27\) and \(\lambda_4 \approx -1.29\).}
}&
\raisebox{-0.5\height}{\inPlot{Figs/spectra6n10k1}}\\
\inPlot{Figs/spectra6n10k2}&
\inPlot{Figs/spectra6n10k3}
\end{tabular}
\end{figure}

For the particular example of \cref{Fspectra6n10E}, with microgrid \(n_x=n_y=10\) in each cell, \ifJ Supplementary Code \cite[Appendix~C,][]{Roberts2024a}\else\cref{cas2d}\fi\ computes the \(200\times200\) Jacobian matrix of the right-hand side of the elasticity equations~\cref{EEaccel}, and then computes its eigenvalues and eigenvectors.
The smallest five eigenvalues are \(\lambda_0=\lambda_1=0\C
\lambda_2=-0.6431\C
\lambda_3=-1.2653\C
\lambda_4=-1.2886\C\ldots\)\,.
\cref{Fspectra6n10ks} plots the displacement fields of eigenvectors corresponding to the three \text{smallest eigenvalues.}

To form a macroscale model to be valid for all timescales longer than some threshold, in this problem with second-order derivatives in time, we must resolve all waves with corresponding frequency smaller than some threshold.
In this non-dimensional example, the frequencies of the waves, \(\omega_m=\sqrt{-\lambda_m}\)\,, corresponding to the computed spectral modes (\cref{Fspectra6n10ks}) are  \(\omega_0=\omega_1=0\C
\omega_2=0.8019\C
\omega_3=1.1248\C
\omega_4=1.1352\C\ldots\)\,.
\begin{itemize}
\item The usual basic homogenisation model is based upon averaging over a cell, and these averages are usually done with a \emph{constant} weight function.  
Such constant weights \emph{implicitly} correspond to the two sub-cell uniform displacement eigenvectors of \cref{Fspectra6n10ks}. 

In the multi-continuum framework proposed and developed herein, this paragraph arrives at the same basis as the usual homogenisation, but by a different rationale (\cref{SSSwls}). 
The rationale is to choose to resolve all dynamics with a timescale longer than some threshold, say we choose a (non-dimensional) threshold approximately~\(5\); that is,
we choose frequencies smaller than \(1/5=0.2\) (roughly).
Here the two `rigid body' cell-modes of \cref{Fspectra6n10ks} corresponding to \(\omega_0=\omega_1=0\) are the only modes with frequency \(\omega_m<0.2\).
Hence \cref{Se2dHomo,SSShoh} choose a two-mode model that constructs, respectively, traditional and higher-order \text{gradient homogenisations.}

\item However, here the third frequency~\(\omega_2\approx0.8\) is separated by a gap from the higher frequencies~\(\omega_3,\omega_4,\ldots\geq1.1\) (albeit not a big gap).
Hence, when we aim to resolve dynamics on a timescale longer than the shorter threshold of~\(1.0\) say, then we need to choose a homogenisation based upon the \emph{three} sub-cell modes of smallest frequency.
The third mode, plotted in \cref{Fspectra6n10ks}, corresponds to rotations of the `hard' core in the centre of each cell (\cref{Fspectra6n10E}).
In the macroscale modelling, this third mode represents the relatively slow vibration of such sub-cell rotations and how they interact and evolve over a large spatial domain \cite[cf.~\S7 of][]{Milton2007}.  

\cref{SSStchm} constructs this three-mode tri-continuum homogenisation, and so gives a rigorous guiding-centre model for macroscale elastic wave dynamics in~2-D.

\item We could choose more `slow' modes, by increasing parameter~\(M\) in the code, to form a multi-continuum homogenisation with more modes resolving the dynamics of more microscale physics.
\end{itemize}

\subsection{Phase-shift embedding of the 2-D heterogeneity}
\label{Spsem2dh}

As in previous sections, a multi-continuum homogenisation is rigorously achieved via embedding the physical system in the ensemble of all phase-shifts of the heterogeneity.
The embedding of~\eqref{EEstress,EEaccel} is done as in \cref{Spse} by considering horizontal displacement field~\(\fu_{ijkl}(t)\) on a 4-D `spatial' lattice in \(xy\theta\phi\)-space, and similarly for the vertical displacement~\vf, and other fields~\(\sigma^{xx},\sigma^{xy},\sigma^{yy}\).
For simplicity, let the sub-cell lattice spacings \(\delta\theta=\delta x\) and \(\delta\phi=\delta y\)\,.
Then the spatial domain for the 4-D lattice is \(\RR^2\times[0,\ell_x)\times[0,\ell_y)\), with all fields \(\ell_x,\ell_y\)-periodic in~\(\theta,\phi\) respectively.
In terms of the shift operator~\(E_i\) defined so that \(E_i^r\fu_{ijkl}:=\fu_{i+r,jkl}\) and similarly for~\(j,k,l\), let's define the 2-D centred difference operator along the diagonals in the \(ijkl\)-lattice:
\begin{equation}
\deltab_{pq}:=\E+_p\E+_q-\E-_p\E-_q\,,
\quad\text{for }p,q\in\{i,j,k,l\}\,.
\label{EdiffOps}
\end{equation}
Then consider the following system that embeds \cref{EEstress,EEaccel}:
\begin{subequations}\label{EEembedE}%
\begin{align}&
\text{\xSym}\quad \sigma_{ijkl}^{xy}
:=\mu_{kl}\big[\deltab_{jl}\fu_{ijkl}/\delta y 
+\deltab_{ik}\vf_{ijkl}/\delta x\big];
\\&
\text{\oSym}\quad \sigma_{ijkl}^{xx}
:=(\lambda_{kl}+2\mu_{kl}){\deltab_{ik}\fu_{ijkl}}/{\delta x} 
+\lambda_{kl}{\deltab_{jl}\vf_{ijkl}}/{\delta y}\,;
\\&
\text{\oSym}\quad \sigma_{ijkl}^{yy}
:=\lambda_{kl}{\deltab_{ik}\fu_{ijkl}}/{\delta x} 
+(\lambda_{kl}+2\mu_{kl}){\deltab_{jl}\vf_{ijkl}}/{\delta y}\,;
\\
&\text{\uSym}\quad 
 \partial_t^2 \fu_{ijkl}={\deltab_{ik}\sigma_{ijkl}^{xx}}/{\delta x} 
+{\deltab_{jl}\sigma_{ijkl}^{xy}}/{\delta y}\,;
\\&
\text{\vSym}\quad  
 \partial_t^2 \vf_{ijkl}={\deltab_{ik}\sigma_{ijkl}^{xy}}/{\delta x} 
+{\deltab_{jl}\sigma_{ijkl}^{yy}}/{\delta y}\,.
\end{align}
\end{subequations}

By the form of the embedding~\cref{EEembedE} every solution~\(\fu_{ijkl}(t), \vf_{ijkl}(t)\) to~\cref{EEembedE} gives rise to solutions of the original physical system \cref{EEstress,EEaccel}, for every phase-shift of the elasticity parameters (\cref{lemeqv}).
To see this, for every~\(k',l'\) define \(u'_{ij} := \fu_{ij,i+k',j+l'}\) and similarly for the other fields.  
Then~\cref{EEembedE}, satisfied by~\(\fu_{ijkl}(t), \vf_{ijkl}(t)\), reduces to the original \cref{EEstress,EEaccel} for~\(u'_{ij}\) with Lam\'e parameters~\(\lambda_{ij},\mu_{i,j}\) replaced by their phase shifts~\(\lambda_{i+k',j+l'} \C \mu_{i+k',j+l'}\).
In particular, for \(k'=l'=0\) this defined~\(u'_{ij}\), with corresponding other fields, satisfies the \text{original \cref{EEstress,EEaccel}.}

Similarly, an ensemble of all solutions to the ensemble of phase-shifted problems \cref{EEstress,EEaccel} forms a solution to the embedding~\cref{EEembedE} (\cref{lemcon}).
That is, the embedding system~\cref{EEembedE}, homogeneous in~\(i,j\) (effectively in~\(x,y\)), is equivalent to the ensemble of phase-shifts of the given heterogeneous \text{system~\cref{EEstress,EEaccel}. }

\subsection{Basis of invariant manifolds}
\label{SSbim}

The \im, multi-continuum, micromorphic, framework wraps around whatever microscale code a user supplies---here it is the embedded microscale system~\eqref{EEembedE}.
The 4-D lattice system~\eqref{EEembedE} is not a \pde, nonetheless the same framework~\cref{EempdeX} of \cref{SSSlbim} applies.
The reason is that the operators and functions in the general form~\cref{EempdeX} may be encompassed by theory \cite[e.g.,][]{Roberts2020a} even when microscale-nonlocal in space.
Indeed, extant forward theory of \im{}s\footcite[e.g.,][]{Carr81, Bates98, Aulbach2000, Chekroun2015a} is easier to apply to such spatially discrete systems because the difference operators are \emph{bounded}, whereas  the usual alternative of derivatives are troublesome \text{unbounded operators.}

In the practical construction of \im\ models \ifJ\cite[see Supplementary Code,][Appendix~C]{Roberts2024a}\else(see \cref{cas2d})\fi\, the spatial differences may be written in terms of shift operators when convenient (as in~\eqref{EEembedE}), or equivalent differential operators when that is convenient.

For macroscale homogenisation, the valid scenarios are that variations in~\(x,y\), and so also in~\(i,j\), are gradual in some useful sense.
The \emph{basis} of the homogenisation is the case where there are effectively no variations in~\(x,y\), hence the microscale shifts \(E_i\fu_{ijkl}\approx \fu_{ijkl}\) and similarly for all fields, that is, \(E_i,E_j\mapsto1\)\,.
In this base case, the centred difference operators \(\deltab_{ik}\mapsto\delta_k\) and \(\deltab_{jl}\mapsto\delta_l\).
Consequently, in this base case, at each and every cross section, parametrised by~\(x,y\) or equivalently indexed by~\(i,j\) which I omit for brevity, the embedding system~\eqref{EEembedE} reduces to the cell-problem%
\begin{subequations}\label{EEcellE}%
\begin{align}&
\text{\xSym}\quad \sigma_{kl}^{xy}
:=\mu_{kl}\big[ \delta_{k}\vf_{kl}/\delta x
+\delta_{l}\fu_{kl}/\delta y \big];
\\&
\text{\oSym}\quad \sigma_{kl}^{xx}
:=(\lambda_{kl}+2\mu_{kl}){\delta_{k}\fu_{kl}}/{\delta x} 
+\lambda_{kl}{\delta_{l}\vf_{kl}}/{\delta y}\,;
\\&
\text{\oSym}\quad \sigma_{kl}^{yy}
:=\lambda_{kl}{\delta_{k}\fu_{kl}}/{\delta x} 
+(\lambda_{kl}+2\mu_{kl}){\delta_{l}\vf_{kl}}/{\delta y}\,;
\\
&\text{\uSym}\quad 
 \partial_t^2 \fu_{kl}={\delta_{k}\sigma_{kl}^{xx}}/{\delta x} 
+{\delta_{l}\sigma_{kl}^{xy}}/{\delta y}\,;
\\&
\text{\vSym}\quad  
 \partial_t^2 \vf_{kl}={\delta_{k}\sigma_{kl}^{xy}}/{\delta x} 
+{\delta_{l}\sigma_{kl}^{yy}}/{\delta y}\,.
\end{align}
\end{subequations}
All fields are to be \(n_x,n_y\)-periodic in~\(k,l\), respectively.
\cref{SSecp} discusses this cell problem: recall that a general solution of~\cref{EEcellE} is a linear combination of modes with frequencies~\(0=\omega_0 \leq \omega_1 \leq \omega_2 \leq\cdots\)\,.
Here we choose to homogenise based upon the cases of \emph{two} and \emph{three} modes of lowest frequency.
The corresponding bi- and tri-continuum homogenisations are constructed as \im{}s of~\cref{EEembedE} as a regular perturbation that accounts for large length-scale modulation across the cells in the \(x,y\)~variables (equivalently in indices~\(i,j\)).

\subsection{Construct multi-continuum homogenised models}
\label{Se2cmmh}

The computer algebra of \ifJ  Supplementary Code \cite[Appendix~C,][]{Roberts2024a}\else\cref{cas2d}\fi\ constructs multi-modal, multi-continuum homogenised models for this example of 2-D elasticity.
One chooses and sets the desired number of modes~\(M\), and the desired order of error in macroscale gradients~\(N\).
We discuss three cases all constructed by the same code that implements \cref{Pconstruct}, and all justified by \cref{Pgensm} as guiding-centre (\cref{Dim}) slow manifold models.

\subsubsection{Classic 2-D homogenisation}
\label{Se2dHomo}

Here base the homogenisation upon perturbing the two zero frequencies \(\omega_1=\omega_2=0\) of the two eigenvalues \(\lambda_0=\lambda_1=0\).
One correspondingly sets the choice \(M=2\) in the constructive \ifJ Supplementary Code \cite[Appendix~C,][]{Roberts2024a}\else code of \cref{cas2d}\fi.
The two amplitudes, order-parameters, macroscale variables~\(U_0,U_1\) are here defined to be the average over a cell of the \(x,y\)-direction displacements, respectively.
These variables are \emph{not} defined like this via assumed cell-averaging, because here there is no such assumption.
In contrast, we choose to \emph{define}~\(U_0,U_1\) to measure these cell modes because the physics of macroscale conservation laws are often best expressed via \emph{unweighted} spatial integrals, and so the effects of such conservation laws are best seen on the macroscale by using macroscale amplitudes which are themselves unweighted integrals/sums over the material cells.

With \(M=2\) modes, the \ifJ Supplementary Code \cite[Appendix~C,][]{Roberts2024a}\else code of \cref{cas2d}\fi\ constructs a slow \im\ here to be  \((\fu,\vf)=v_0(\theta,\phi)U_0(t,x,y)+v_1(\theta,\phi)U_1(t,x,y)+\cdots\) in terms of the two leading eigenvectors~\(v_0,v_1\) plotted in \cref{Fspectra6n10ks}, and where the ellipsis represents some computed corrections in gradients of~\(U_0,U_1\) which for simplicity are not recorded here.  
As in \cref{Spse}, such an \im\ of the phase-shifted embedding problems leads to the spatial displacement fields of the original problem to be \((u,v)=v_0(x,y)U_0(t,x,y)+v_1(x,y)U_1(t,x,y)+\cdots\)\,.
The evolution of the amplitudes~\(U_0,U_1\) then give the homogenised model:  to three decimal places it is\begin{subequations}\label{EE2dBase}%
\begin{align}
\DD t{U_0}&
=0.144\DD x{U_0}+0.020\DD y{U_0}+0.092\Dx xy{U_1}
+\Ord{\grad^3},
\label{E2dBasex}
\\
\DD t{U_1}&
=0.020\DD x{U_1}+0.144\DD y{U_1}+0.092\Dx xy{U_0}
+\Ord{\grad^3}.
\label{E2dBasey}
\end{align}
\end{subequations}
The coefficients of these \pde{}s are the effective elastic moduli for long-waves in the material.
The symmetry of the heterogeneous cells in the \(x,y\)-directions (e.g., \cref{Fspectra6n10E}) results in the symmetry in coefficients, and hence in the effective elastic tensor, apparent in the homogenisation~\cref{EE2dBase}.
However, the coupled macroscale wave-\pde{}s~\cref{EE2dBase} are anisotropic due to the square cells distinguishing \text{among various directions.}

The effective wave equations~\cref{EE2dBase} for long-waves will agree with \pde{}s as obtained by every good classic homogenisation method:  the second-oder gradient terms on the right-hand side arise from solving the same cell problem with the same lower order forcing and same solvability condition, so the model~\pde\ at this leading order will be the same.
Such agreement is because the out-of-equilibrium chain rule, \(\partial_t\fu=\D{U_m}\fu\cdot \D t{U_m}\) et al., has no effect at this leading order (\cref{LgenHomo}). 

Our dynamical systems approach empowers us to improve such a basic homogenisation by straightforwardly accounting for more sub-cell physics. 
The next \cref{SSShoh,SSStchm} describe two different analytic ways to do so.

\paragraph{Quantitatively estimate errors}
The asymptotic errors~\Ord{\grad^3} in \pde{}s~\cref{EE2dBase}, as also for the errors in~\cref{EE2dO4,EE2dM3}, could be quantitatively estimated. 
\cite{Roberts2016a} developed a general mathematical theory for reduced-order multiscale modelling and homogenisation in multiple spatial dimensions, and their expression~(52) is an explicit novel general formula for the remainder error that applies to~\cref{EE2dBase}.
However, the expression is sufficiently complicated that we leave this aspect for further investigation.
Instead I suggest that the fourth-order gradient terms in the next model~\cref{EE2dO4} give a practical estimate of the local error in any predictions made by the second-order model~\cref{EE2dBase}.

\subsubsection{Higher-order gradient homogenisation}
\label{SSShoh}

It is straightforward in this framework to proceed to higher-order in macroscale gradients---a more rigorous route than the second-gradient homogenisation heuristics of \cite{Forest2011}~[\S2].  
Such higher-orders account for more physical interactions in the sub-cell dynamics and how these affect the evolution over the macroscale.

For example, \ifJ Supplementary Code \cite[Appendix~C,][]{Roberts2024a}\else\cref{cas2d}\fi\ readily constructs the fourth-order model, here reported to two significant digits, and neglecting terms with numerically small \text{coefficients\({}<0.001\):}
\begin{subequations}\label{EE2dO4}%
\begin{align}
\DD t{U_0}&
=0.144\DD x{U_0} +0.020\DD y{U_0} +0.092\Dx xy{U_1}
\nonumber\\&\quad{}
+0.0037\Dn x4{U_0} 
+0.0042\frac{\partial^4U_0}{\partial x^2\partial y^2} 
+0.0042\frac{\partial^4U_1}{\partial x\partial y^3}
+\Ord{\partial_x^5+\partial_y^5},
\label{E2dO4x}
\\
\DD t{U_1}&
=0.020\DD x{U_1}+0.144\DD y{U_1}+0.092\Dx xy{U_0}
\nonumber\\&\quad{}
+0.0037\Dn y4{U_1} 
+0.0042\frac{\partial^4U_1}{\partial x^2\partial y^2} 
+0.0042\frac{\partial^4U_0}{\partial x^3\partial y}
+\Ord{\partial_x^5+\partial_y^5}.
\label{E2dO4y}
\end{align}
\end{subequations}
Physically, the fourth-order derivative terms on the second lines of~\cref{E2dO4x,E2dO4y} characterise the physics of macroscale anisotropic wave-dispersion in the homogenised model of this heterogeneous material.

\subsubsection{Tri-continuum homogenised model}
\label{SSStchm}

The alternative way to account for more physical interactions in the sub-cell dynamics is to retain more sub-cell modes to form a multi-continuum model.
For this example, let's choose to retain the gravest three modes.

Here the three modes correspond to eigenvalues \(\lambda_0=\lambda_1=0\) and \(\lambda_2=-0.6431\), and are separated from all the other eigenvalues headed by \(\lambda_3=-1.2653\).  
\cref{Fspectra6n10ks} plots the physical structure of the three eigenmodes corresponding to these three gravest eigenvalues.
The three eigenmodes are two of displacement in \(x,y\)-directions, and a sub-cell mode representing rotation of the `hard'-centre of each cell.
Because of this third rotational mode, this model is a \emph{micropolar continua} in the sense of Eringen \cite[p.4]{Maugin2010a}.
\cite{Combescure2022} would call the model derived here a ``Cosserat model [because it is] adding only local rotational degrees of freedom to the kinematic field''.
Importantly, these modes are \emph{not} assumptions we impose onto the physics, instead these are modes that the sub-cell physics informs us \emph{are the appropriate} sub-cell structures. 

With \(M=3\), and using the approximate homological equation~\cref{EgenHomoSimp}, \ifJ Supplementary Code \cite[Appendix~C,][]{Roberts2024a}\else\cref{cas2d}\fi\ takes 32~iterations to construct, with residual of~\(10^{-7}\), a slow \im\ here to be  \((\fu,\vf)=v_0(\theta,\phi)U_0(t,x,y)+v_1(\theta,\phi)U_1(t,x,y)+v_2(\theta,\phi)U_2(t,x,y)+\cdots\) in terms of the three leading eigenvectors~\(v_0,v_1,v_2\) plotted in \cref{Fspectra6n10ks}.
The ellipsis represents some computed corrections in spatial gradients of~\(U_0,U_1,U_2\) which for simplicity are not recorded here.  
To two significant digits, and neglecting terms with numerically small coefficients\({}<0.001\), the correspondingly constructed three-mode, tri-continuum, second-order homogenised evolution is governed by the following \pde{}s:
\begin{subequations}\label{EE2dM3}%
\begin{align}
\DD t{U_0}&
=0.144\DD x{U_0}+0.020\DD y{U_0}+0.092\Dx xy{U_1}
\notJbreak{}
+0.0021\DD x{U_2}
-0.0024\Dx xy{U_2}
+\Ord{\grad^3},
\label{E2dM3x}
\\
\DD t{U_1}&
=0.020\DD x{U_1}+0.144\DD y{U_1}+0.092\Dx xy{U_0}
\notJbreak{}
+0.0021\DD y{U_2}
-0.0024\Dx xy{U_2}
+\Ord{\grad^3},
\label{E2dM3y}
\\
\DD t{U_2}&
=-0.643\,U_2
+0.0032\left(\DD x{U_0}+\DD y{U_1}\right)
-0.0023\Dx xy{(U_0+U_1)}
\nonumber\\&\quad{}
-0.031\left(\DD x{U_2}+\DD y{U_2}\right)
+\Ord{\grad^3}.
\label{E2dM3r}
\end{align}
\end{subequations}
The above is a mathematically and physically rigorous tri-continuum model for the homogenised dynamics of the 2-D heterogeneous material.

As well as resolving shorter time scales than~\cref{EE2dBase,EE2dO4}, this three-mode tri-continuum homogenisation contributes to the wave dispersion resolved by the fourth-order model~\cref{EE2dO4}.
To see this use a quasi-adiabatic, quasi-static, approximation to the rotational mode~\cref{E2dM3r}, that \(0\approx-0.64\,U_2 +0.0032(U_{0xx}+U_{1yy}) -0.0023(U_{0xy} +U_{1xy})\), to thence deduce the sub-cell rotational  amplitude
\(U_2\approx +0.0049(U_{0xx}+U_{1yy})
-0.0035(U_{0xy}+U_{1xy})\),
in terms of the local curvatures of the macroscale mean displacements~\(U_0,U_1\).
Substituting this into~\cref{E2dM3x,E2dM3y} leads to fourth-order terms in the form of~\cref{EE2dO4}.  
This contribution does not complete the fourth-order terms in~\cref{EE2dO4} because, to be consistent, the \pde{}s~\cref{E2dM3x,E2dM3y} should also have their fourth-order terms constructed as in \cref{SSShoh}.

Physically, that this quasi-adiabatic approximation does not generate any second-order effects in the displacements thus indicates that here there is a relatively weak interaction from sub-cell rotation back into the two principal displacement modes.
The significance of this tri-continuum homogenisation is that it is valid on shorter timescales than the bi-continuum~\cref{EE2dBase,EE2dO4}, and potentially valid on shorter length-scales.
Lastly, as \cite{Maugin2010a} comments, resolving such ``a rotating microstructure allows for the introduction of wave modes of rotation of the optical type'' into the macroscale model~\cref{EE2dM3}.

\section{Conclusion}

This article synthesises a novel, systematic, rigorous and practical approach to creating families of multi-continuum, micromorphic, macroscale homogenised models of microscale heterogeneity in mechanics.
\cref{SoneDintro,Shceq} develops the basics of the approach in 1-D space, whereas \cref{Sgentheory,Selastic2d} addresses the complexities of general nonlinear systems in multi-D space, potentially with quasi-periodic microscale, underpinned by nonlinear dynamical systems theory.
We have developed the techniques called for in the recent review by \cite{Fronk2023} [\S2.4~Future Work] who comment  ``In terms of analysis, generalizing the analysis of Lamb modes to \emph{periodic} nonlinear plates is non-trivial, and will likely require a hybrid analytical-numerical approach'':  \cref{Sgentheory} provides such a nonlinear homogenisation framework, including quasi-periodicity, and \cref{Selastic2d} uses an example hybrid analytical-numerical approach.

An outstanding issue for multi-continua homogenisation is to decide on a suitable set of microscale structures.
Here the proposed rationale (\cref{SSSmmmcme,SSSsd,SSSwls}) is to choose to retain modes slower than a chosen threshold---a threshold selected according to the time-scales required for the intended application.   
An alternative is to choose modes according to the required space-scales of the application (\cref{SSSisr}).  
In many cases these two rationale are more-or-less equivalent.  
But sometimes, due to physical or geometric symmetries, the two rationale are different enough to justify a different set of modes for each scenario.
The systematic approach developed here simplifies considerably much of the previous difficulty in choosing an appropriate model for a given microstructure.

The approach encompasses the homogenisation of nonlinear systems because it uses nonlinear dynamical systems theory (via \cref{Pgenft,Pgenabt,Pgensm}).  
\cref{SShhn} discusses an example, and illustrates that in modelling or homogenising \emph{nonlinear} systems the nature of the time evolution operator~\alphaD\ significantly affects the homogenisation through significant changes in the algebra.

The theory and results are \emph{not} limited to the traditional scale separation that the micro:macro scale ratios \(\ell/L\to0\) (\cref{SoneDintro,Sgentheory}).
The theory and results apply at finite scale separation of real physics and engineering applications (e.g., \cref{Selastic2d}).
Indeed our systematic dynamical systems framework, coupled with high-order computer algebra, may predict a sharp quantitative limit to the spatial resolution of an homogenisation (e.g., \cref{SSScsw}).

\cref{Shcegbmm} verifies that our sound approach to modelling is transitive: the slow manifold homogenisation of a two-mode bi-continuum homogenisation of some physical system is the same as the slow manifold homogenisation of the physical system.

A further advantage of this nonlinear dynamical systems approach to homogenisation is that the associated theory and methods provides a sound route for correctly modelling initial conditions and domain boundary conditions for the homogenised model.
The issue of initial conditions and boundary conditions is often trivialised by the mathematical scale separation limit ``\(\ell/L\to0\)'', but is nontrivial at the finite scale separation of real applications \cite[e.g.,][pp.1034--5 and p.221 resp.]{Boutin1996, AbdoulAnziz2018}. 
Methodology to construct the correct initial conditions for reduced order models may be based upon solving a dual problem \cite[]{Roberts97b, Roberts89b}.
These methods then generalise to potentially provide correct domain boundary conditions for the homogenisation variables \cite[]{Roberts92c}, as explored for some elasticity examples \cite[]{Roberts93, Chen2016}.
Future research is needed to develop these techniques to encompass the scenarios addressed in this article, and thus to complete the homogenisation methodology.

Such correct projection of initial conditions will also inform a user of how to project forcing and uncertainty into an homogenised model \cite[e.g.,][\S12.4]{Roberts2014a}.

\paragraph{Acknowledgements}
I thank reviewers for their constructive comments, and
I thank colleagues Judy Bunder, Pavel Bedrikovetsky, Yannis Kevrekidis, and Thien Tran-Duc  for their encouragement, and Arthur Norman and colleagues who maintain and develop the computer algebra software Reduce.  
The Australian Research Council Discovery Project grant DP220103156 helped support \text{this research.}

\ifJ
\else\appendix

\section{Computer algebra construction}
\label{cas}

 
This Reduce-algebra code constructs any chosen
multi-continuum, micromorphic, invariant manifold
homogenisation of heterogeneous diffusion~\cref{Ehdifpde}
discussed by \cref{SoneDintro}. Define main parameters:
\begin{reduce}
mm:=3;     
theCase:=1;
cases:={ {3,3},{5,5},{3,21},{21,4},{3,31},{31,4} }$
alpha:=2;  
nonlinearCase:=0;  
\end{reduce}
Non-dimensionalise the problem so the heterogeneity is
\(2\pi\)-periodic---a variety of nonlinear trigonometric
dependence is possible, but here code~\cref{EGkappa}. The
`strength' of the heterogeneity must be parametrised by
variable~\(a\). Throughout, let~\verb|q| denote~\(\theta\).
\begin{reduce}
kappa:=1/(1+a*cos(q));
\end{reduce}

Optionally allow the strength to vary with macroscale~\(x\).
By default it varies gradually in~\(x\) as all derivatives
in~\(x\) are counted by~\(d\).   
\begin{reduce}
if 0 then depend a,x; 
\end{reduce}

Optionally, introduce heterogeneous nonlinearity
\cref{EhHetNon,egkappaeta} and set order of nonlinear error:
here~\Ord{\gamma^2} would resolve quadratic terms.  It
appears that multi-mode, nonlinear, wave-like systems need
more careful solution of homological equation.
\begin{reduce}
for p:=nonlinearCase:nonlinearCase do let gamma^(p+1)=>0;
if nonlinearCase>0 then begin factor gamma,vv,c1,c2; 
    if alpha>1 and mm>1
    then rederr("Need alpha=1 for nonlinear multi-mode ****"); 
end; 
eta:=c1*cos(q)+c2*sin(2*q);
\end{reduce}

Example high-order computation times are (\(\gamma=0\)): 
case~3, \verb|{3,21}|, uses 23~iterations in 20~secs; 
case~5, \verb|{3,31}|, uses 33~iterations in 155~secs;
case~4, \verb|{21,4}|, uses 24~iterations in 2~secs;
case~6, \verb|{31,4}|, uses 34~iterations in 3~secs.
That is, it terminates in \(N+P-1\) iterations.

Extract the orders of error from the case.
\begin{reduce}
ordd:=part(cases,theCase,1); 
orda:=part(cases,theCase,2); 
if ordd<2 then let d^2=>0;
if orda<2 then let a^2=>0;
\end{reduce}
Improve formatting of written output:
\begin{reduce}
on div; off allfac; on revpri;
factor d,df,uu;
\end{reduce}

Write approximations to the slow manifold homogenisation of
the embedding \pde~\cref{Eemdifpde}, such
as~\cref{EE3Mevol,eqimphpde,EEhetnonIM,EEhetnonIM2}, in
terms of `modal' fields~\(U_i(t,x)\), denoted
by~\verb|uu(i)|, that evolves according to \(\alphaD
t{U_i}=\verb|dudt(i)|\) for whatever approximation
\verb|dudt(i)| contains.  In the case of \(\alpha=2\),
define \(V_i:=\D t{U_i}\).
\begin{reduce}
array dudt(mm-1);
operator uu;  depend uu,t,x;
if alpha=1 
then let { df(uu(~i),t)=>   dudt(i)
    , df(uu(~i),t,x)   =>df(dudt(i),x)
    , df(uu(~i),t,x,2) =>df(dudt(i),x,2)
    } 
else begin operator vv; depend vv,t,x;
    let { df(uu(~i),t)=>   vv(i)
    ,     df(vv(~i),t)=>   dudt(i)
    , df(vv(~i),t,x)  =>df(dudt(i),x)
    , df(vv(~i),t,x,2)=>df(dudt(i),x,2)
    }
end;
\end{reduce}

\subsection{Iteration systematically constructs multi-modal model}
\label{seccaisc1d}

We iteratively construct the improved homogenizations
\cref{EE3Mevol,eqimphpde,EEhetnonIM,EEhetnonIM2}. Recall
that parameters~\verb|ordd| and~\verb|orda| specify the
orders of error.  
\begin{reduce}
write "
Second, Iteratively Construct
-----------------------------";
maxit:=99;
\end{reduce}

Expand diffusivity~\(\kappa\) as Taylor series in
heterogeneity parameter~\(a\).
\begin{reduce}
if sub(a=0,kappa) neq 1 
then rederr("kappa for a=0 must be scaled to one");
kappa:=taylor(kappa,a,0,orda);
kappa:=trigsimp( taylortostandard(kappa) ,combine)$
\end{reduce}

Start the iteration from the base multi-modal approximation
for the field and its evolution, given the
eigenmodes~\(v_m(\theta)\) are~\(1\C \sin\theta\C
\cos\theta\C \sin2\theta\C \ldots\)\,.  Store the
corresponding eigenvalues~\(\lambda_m\) of the modes in
array~\verb|lams|.
\begin{reduce}
u:=(for k:=0:mm-1 sum uu(k)
    *(if evenp(k) then cos(k/2*q) else sin((k+1)/2*q)));
array lams(mm);
for k:=0:mm-1 do lams(k) := 
    -(if evenp(k) then k/2 else (k+1)/2 )^2;
for k:=0:mm-1 do write dudt(k) := lams(k)*uu(k);
\end{reduce}
Iteratively seek solution to the specified orders of errors.
\begin{reduce}
for it:=1:maxit do begin write "
**** ITERATION ",it;
\end{reduce}

Progressively truncate the order of the order parameter so
that we control the residuals better: the bounds in these
if-statement are the aimed for ultimate order of errors.  
\begin{reduce}
  if it<ordd then let d^(it+1)=>0;
  if it<orda then let a^(it+1)=>0;
\end{reduce}

Compute the \pde\ residual via the flux, optionally
including heterogeneous nonlinear advection. 
\begin{reduce}
  flux:=-kappa*(d*df(u,x)+df(u,q));
  if gamma neq 0 then flux:=flux+gamma*eta*u^2/2;
  pde:=df(u,t,alpha)+d*df(flux,x)+df(flux,q);
  pde:=trigsimp(pde,combine);
\end{reduce}
Trace print either the residual expression or its length.
\begin{reduce}
  if length(pde)<20 then write pde:=pde
  else write lengthpde:=length(pde);
\end{reduce}

Update the evolution via the solvability conditions
over~\(\theta\), simultaneously removing the `resonant'
terms from the right-hand side of the homological
equation~\cref{E1dHomologic}.
\begin{reduce}
  rhs:=pde+( gd:=-(pde where {sin(~a)=>0,cos(~a)=>0}) );
  dudt(0):=dudt(0)+gd;
  for k:=1:mm-1 do begin
    vk:=(if evenp(k) then cos(k/2*q) else sin((k+1)/2*q));
    dudt(k):=dudt(k)+( gd:=-coeffn(rhs,vk,1) );
    rhs:=rhs+gd*vk;
  end;
\end{reduce}
For solving for corrections to the invariant manifold field
from the residual we first need an operator~\verb|homolog|
to account for all factors of amplitude~\(U_k\) in the
homological equation~\cref{E1dHomologic}. Define these once
only in the first iteration.
\begin{reduce}
if it=1 then begin
    operator homolog; linear homolog;
    depend q,qu;  depend uu,qu;
    let { homolog(~~a*uu(~m)^~~p,qu,~l) 
            => uu(m)^p*homolog(a,qu,l+p*lams(m))
        , homolog(~~a*df(uu(~m),x)^~~p,qu,~l) 
            => df(uu(m),x)^p*homolog(a,qu,l+p*lams(m))
        , homolog(~~a*df(uu(~m),x,~n)^~~p,qu,~l) 
            => df(uu(m),x,n)^p*homolog(a,qu,l+p*lams(m))
        };
\end{reduce}
Also setup these~\verb|intq| transformations to solve the
components of the homological equation~\cref{E1dHomologic}.
\begin{reduce}
    intq:={ homolog(sin(~~n*q),qu,~lam)=>sin(n*q)/(lam+n^2)
          , homolog(cos(~~n*q),qu,~lam)=>cos(n*q)/(lam+n^2)
          };
end;
\end{reduce}
With the above solvable right-hand side, now update the
multi-continuum, multi-modal, invariant manifold field.
\begin{reduce}
  rhs:=homolog(rhs,qu,0);
  u:=u+( ud:=-(rhs where intq) );
\end{reduce}

Finish the loop when the residual of the \pde\ is zero to
the specified order of error,
\begin{reduce}
  if pde=0 then write "Success: ",it:=it+100000;
end;
showtime;
if pde neq 0 then rederr("Iteration Failed to Converge");
\end{reduce}

\subsection{Post-process}

Post-process writes out the multi-continuum, multi-modal,
evolution, and the corresponding invariant manifold, to one
order lower, if not too long. 
\begin{reduce}
uLow:=(u*a*d)/a/d$
if length(uLow)<20 then write uLow:=uLow;
for i:=0:mm-1 do write "Dt",alpha," U",i," = ",dudt(i);
\end{reduce}

\paragraph{Optional high-order in heterogeneity}
Optionally output selected coefficients for power series
analysis.
\begin{reduce}
if orda>10 then begin
    c0:=coeffn(dudt(0),df(uu(0),x,2),1)/d^2$
    c1:=coeffn(dudt(1),uu(1),1)$
    c2:=coeffn(dudt(2),uu(2),1)$
    on rounded;
    out "cas1dCs.m";
    write "
    cs0 = ",coeff(c0,a),"
    cs1 = ",coeff(c1,a),"
    cs2 = ",coeff(c2,a),"
    for i=1:3, figure(i)
    figname=['Figs/' mfilename num2str(i)]
    clf, set(gca,'position',[.2 .2 .54 .54])
    switch i
    case 1, radiusConverge(cs0)
    case 2, radiusConverge(cs1)
    case 3, radiusConverge(cs2)
    end
    exportgraphics(gcf,[figname '.pdf'] ,'ContentType','vector')
    matlab2tikz([figname '.tex']...
    ,'showInfo',false,'showWarnings',false,'parseStrings',false ...
    ,'extraCode',['\tikzsetnextfilename{' figname '}'] ...
    ,'extraAxisOptions','max space between ticks={50}' )
    end
    "$
    shut "cas1dCs.m";
    off rounded;
end;
\end{reduce}

\paragraph{Optional high-order spatial structure}
Optionally apply discovered operator that simplifies
evolution of most~\(a^2,a^3\) terms.  Define
\(\verb|dinv|:=(1+4/9\partial_x^2)^{-1}\).
\begin{reduce}
if ordd>10 and orda=4 then begin
    write "Multiply a^2,a^3 by operator  D = 1+4/9*df(,x,x)";
    procedure dop(f); f+4/9*d^2*df(f,x,2);
    operator dinv; linear dinv;
    for i:=0:mm-1 do write "D*Dt",alpha," U",i," = "
        , coeffn(dudt(i),a,0)
        +a*coeffn(dudt(i),a,1)
        +a^2*dinv(dop(coeffn(dudt(i),a,2)),x)
        +a^3*dinv(dop(coeffn(dudt(i),a,3)),x);
end;
\end{reduce}

Computer algebra fin.
\begin{reduce}
end;
\end{reduce}



\section{Computer algebra construction of high contrast laminate}
\label{cashc}


This Reduce-algebra code constructs any chosen
multi-continuum, micromorphic, invariant manifold
homogenisation of the high-contrast example
of~\cref{Ehdifpde2} discussed by \cref{Shceq}. Write
approximations, such as
\cref{EEhcegManifold1,EEhcegManifold2}, to the slow manifold
homogenisation of the embedding \pde~\cref{Eemdifpde2} in
terms of `modal' macro-scale fields~\(U_i(t,x,y)\), denoted
by~\verb|uu(i,0,0)|, that evolves according to \(\Dn
t\alpha{U_i} =\verb|dudt(i)|\) for whatever \verb|dudt(i)|
happens to be constructed, and for either \(\alpha=1\) for
diffusive case or~\(\alpha=2\) for wave case.  

In principle, for linear systems and whenever the time
evolution operator~\alphaD\ commutes with spatial
derivatives, then the construction is also valid and gives
the dynamic homogenisation \(\alphaD{U_i} =\verb|dudt(i)|\).
 
Set the desired derivative parameter~\(\alpha\), the
required dimensionality~\(M\) of the invariant manifold
multi-modal multi-continuum model, and the order of error in
spatial gradients, \(\verb|ordd|=N+1\)\,.  Optionally
include nonlinear advection \(-\gamma uu_x\), by setting
appropriate truncation in~\verb|gamma|.
\begin{reduce}
alpha:=1;
mm:=2;
ordd:=5;
let gamma=>0;
\end{reduce}
The series approximations are only computed approximately,
so set a desired relative error for the numerical
coefficients. And use \verb|zeroSmall| to eliminate
negligible terms.
\begin{reduce}
relTolerance:=1e-8;
\end{reduce}

The macroscale amplitudes, order parameters, depend upon
time (depending upon~\(\alpha\)) and space~\((x,y)\).
Let \verb|uu(i,m,n)| denote \(\partial^{m+n} U_i 
/\partial x^m \partial y^n\).
\begin{reduce}
array dudt(mm-1);
operator uu;  depend uu,x,y,t;
let { df(uu(~i,~m,~n),x)  => uu(i,m+1,n)
    , df(uu(~i,~m,~n),y)  => uu(i,m,n+1)
    , df(uu(~i,~m,~n),y,2)=> uu(i,m,n+2)
    };
if alpha=1 
then let df(uu(~i,~m,~n),t)  => df( dudt(i), x,m,y,n)
else let df(uu(~i,~m,~n),t,2)=> df( dudt(i), x,m,y,n);
\end{reduce}

Improve printed output, and use numerical, floating-point,
coefficients.
\begin{reduce}
on div; off allfac; on revpri;
on rounded; print_precision 4$
factor d;
in "zeroSmall.txt"$
tolerance:=relTolerance$
\end{reduce}

\subsection{Parametrise a non-dimensional case}
For simplicity, non-dimensionalise space and time so that
\(\kappa_1=1\) and microscale periodicity is \(\ell=2\pi\).
Choose a layer that is say~6\% of the microscale so there
are a reasonable nine sub-microscale points across the
layer.
\begin{reduce}
ell:=2*pi;
eta:=0.06*ell; 
chi:=1;        
kappa0:=eta/chi/ell; 
\end{reduce}
Choose to discretise the \(\theta\)-structure across a cell
in terms of values at \(n\)~equi-spaced points, use \verb|q|
to denote~\(\theta\). Set the microscale grid on the
interval \([-\ell/2,+\ell/2] =[-\pi,\pi]\), and periodic: 
\begin{reduce}
n:=128;    
procedure q(j); ((j-1/2)/n-1/2)*ell;
dq:=q(2)-q(1);
\end{reduce}
Define a procedure to write out the parameters as a comment:
used to label various output files.
\begin{reduce}
procedure writeTheCase; << write "
    ," alpha=",alpha ,", Nmodes=",mm ,", n=",n 
    ,", ell=",ell ,", eta/ell=",eta/ell ,", chi=",chi;
    write "\tikzsetfigurename{Figs/hceg2Manifold",mm,"}"; 
    >>;
\end{reduce}
Set some elementary matrices, and then define the discrete
operator of heterogeneous interaction across a cell.
\begin{reduce}
matrix zeros(n,1),ones(n,1),Id(n,n),ee(n,n);
for j:=1:n do ones(j,1):=1; 
for j:=1:n do Id(j,j):=1;   
\end{reduce}
Matrix to shift vector up, circularly: \(Eu_j:=u_{j+1}\);
and \(E^T=E^{-1}\) is shift vector down,
\(E^Tu_j:=u_{j-1}\).
\begin{reduce}
for j:=1:n do ee(j,if j<n then j+1 else 1):=1; 
\end{reduce}
Set the cell coefficient matrix, and compute its norm (max
column sum of abs). The matrix 
\begin{equation*}
\cK:=\begin{bmatrix}
-\kappa_{1/2}-\kappa_{3/2}&\kappa_{3/2}&\cdots&\kappa_{1/2}
\\\kappa_{3/2}&-\kappa_{3/2}-\kappa_{5/2}&\cdots&0
\\0&\ddots&\ddots&\vdots
\\\kappa_{1/2}&\cdots&\kappa_{n-1/2}&-\kappa_{n-1/2}-\kappa_{1/2}
\end{bmatrix}
\end{equation*}
where \(\kappa_{j-1/2}\) is the diffusivity\slash elasticity
at \(\theta_{j-1/2}=\theta_j-\d \theta/2\). So diagonal of
matrix \verb|emkapep| is the coefficients `shifted' down a
half, that is~\(E^{-1/2}\kappa E^{1/2}\).  Set the
heterogeneous matrix accordingly, \(\cK := \delta \kappa
\delta /\d \theta^2\).
\begin{reduce}
matrix kk(n,n),emkapep(n,n);
for j:=1:n do emkapep(j,j) :=
    (if abs(cos(q(j-1/2)/2)*2)<eta/2 then kappa0 else 1);
kk := (ee-Id)*emkapep*(Id-tp ee)/dq^2$
normkk := max(for i:=1:n collect (for j:=1:n sum abs(kk(i,j))));
\end{reduce}

\subsection{Spectrum of cell problem}
Find the gravest eigenvalues and eigenvectors.  The first is
the constant eigenvector corresponding to zero eigenvalue.
All eigenvectors are found as unit vectors, and then
normalised to have unit mean-square in the integral norm.
\begin{reduce}
write "Finding leading eigenvalues and eigenvectors";
array evecs(mm),evals(mm);
evecs(0):=ones$
evals(0):=0;
\end{reduce}
The \(k\)th~eigenvector and eigenvalue are found by
numerical iteration, starting from approximations that
assume a thin low-diffusive layer is around
\(\theta=\pm\ell/2=\pm\pi\). 
\begin{reduce}
matrix evec(n,1);
for k:=1:mm do begin
\end{reduce}
Initialise rough approximate eigenvectors. 
\begin{reduce}
    for j:=1:n do evec(j,1):=
        (if evenp(k) then cos(k/2*q(j)) else sin(k/2*q(j)));
\end{reduce}
This iteration should quadratically converge to an
eigenvalue-eigenvector pair, and should not miss any
grave modes. 
\begin{reduce}
    for it:=1:9 do begin 
        if it>1 then evec:=1/(kk-eval*Id)*evec; 
        evec:=evec/sqrt(part(tp evec*evec,1,0));
        eval:=(tp evec)*kk*evec; 
        res:=(kk-eval*Id)*evec;  
        write normres:=for j:=1:n sum abs(res(j,1));
        if normres<1e-8*normkk 
        then write "success: ",it:=it+1000;
    end;
    if normres>1e-8*normkk then rederr("failed eigenvalue iteration");
    evecs(k):=evec*sqrt(2*pi/dq); 
    write evals(k):=eval(1,1);
    for l:=k step -1 until 1 do 
        if evals(l)>evals(l-1) then begin
            tmp:=evals(l); evals(l):=evals(l-1); evals(l-1):=tmp;
            tmp:=evecs(l); evecs(l):=evecs(l-1); evecs(l-1):=tmp;
        end;
end;
\end{reduce}

Optionally plot the eigenvectors. Also check that the
pattern of signs in the eigenvectors is as expected.
\begin{reduce}
array xv(mm);
for k:=0:mm do xv(k):=
    for j:=1:n collect {q(j),part(evecs(k),j,0)};
if 0 then begin
    if mm=1 then plot(xv(0),xv(1));
    if mm=2 then plot(xv(0),xv(1),xv(2));
    if mm=3 then plot(xv(0),xv(1),xv(2),xv(3));
end;
for k:=1:mm do begin
    v:=evecs(k);
    write signs := for j:=3:n-1 sum 
        abs(sign(v(j,1))-sign(v(j-1,1)))/2;
    if signs neq k then rederr("wrong pattern of signs");
end;
\end{reduce}
Output pgfplots commands to subsequently draw the
eigenvectors in \LaTeX\ via pgfplots.
\begin{reduce}
out "hcegEvecs.tex";
writeTheCase();
write"\begin{tikzpicture}
\makeatletter\let\gob\@gobble\makeatother
\begin{axis}[no marks,
    xlabel={microscale $\theta$},ylabel={eigenvector},
    legend pos=south east,legend style={font=\tiny}]";
for k:=0:mm do begin
    write "\addplot+[",if k=mm then "dashed" else ""
        ,"] plot coordinates {";
    foreach qu in xv(k) do write 
        "    (",part(qu,1),",",part(qu,2),")\gob"$
    write "};
    \addlegendentry{$",evals(k),"$}"$
end;
write "\end{axis}
\end{tikzpicture}";
shut "hcegEvecs.tex";
\end{reduce}

\subsection{Form inverse for homological updates}
\label{SSfifhu}
Corrections to the approximate invariant manifold shape and
evolution are done by the simplified the homological
equation~\cref{EgenHomoSimp}, a physics-informed linear
system \cite[e.g.,][\S5.3]{Roberts2014a}. The following
inverse would be `wrong' for multi-mode updates involving
non-zero eigenvalues in the manifold, but as part of an
adiabatic iteration it is good enough upon more iterations.
\begin{reduce}
matrix llinv(n+mm,n+mm),rhs0(n+mm,1);
for i:=1:n do for j:=1:n do llinv(i,j):=kk(i,j);
for i:=1:mm do begin
    v:=evecs(i-1); 
    for j:=1:n do llinv(n+i,j):=-v(j,1);
    for j:=1:n do llinv(j,n+i):=-v(j,1);
end;
llinv:=1/llinv$
write "Formed Linv OK";
off roundbf; write "--- I turn off RoundBF";
\end{reduce}

\subsection{Iteration systematically constructs multi-modal model}
\label{seccaischc}

Second, we iteratively construct the improved
homogenizations~\cref{EEhcegManifold1,EEhcegManifold2}.     
\begin{reduce}
write "
Iteratively Construct
---------------------";
maxit:=99;    
\end{reduce}

Start the iteration from the trivial approximation for the
field and its evolution.
\begin{reduce}
u:=for i:=0:mm-1 sum uu(i,0,0)*evecs(i)$
for i:=0:mm-1 do write dudt(i):=evals(i)*uu(i,0,0);
\end{reduce}
Seek solution to the specified orders of errors (via the
instant evaluation property of the for loop index).  Use
variable~\(d\) to count the number of \(x,y\)-derivatives.
\begin{reduce}
for it:=ordd:ordd do let d^it=>0;
for it:=1:maxit do begin 
\end{reduce}

Compute the \pde\ residual via the flux, trace printing the
length of the residual expression.  The embedding
\pde~\cref{Eemdifpde} is discretised across the cell in
terms of mean~\(\mu\) and difference~\(\delta\) operators,
but retaining algebraic \(x,y,t\)-dependence, as
\begin{align*}
x\text{-flux }F&:=
-\kappa\left(\mu\D x\fu +\frac{\delta\fu}{\d \theta}\right),
& \text{residual}&:=
\D t\fu+\mu\D xF+\frac{\delta F}{\d \theta}+\kappa\DD y\fu\,,
\end{align*}
for micro-grid spacing \(\d\theta:=2\pi/n\)\,. Here compute
\(E^{-1/2}F\), and write out the overall size of the
coefficients on the residual.
\begin{reduce}
    emflux:=-emkapep*((Id+tp ee)/2*d*df(u,x)+(Id-tp ee)*u/dq);
    pde:=df(u,t,alpha)+(Id+ee)/2*d*df(emflux,x)+(ee-Id)*emflux/dq
        -emkapep*df(u,y,2)*d^2;
\end{reduce}
Optionally add in example nonlinearity (when not 
\verb|gamma=>0|).
\begin{reduce}
    for j:=1:n do pde(j,1):=pde(j,1)+gamma*u(j,1)*(d*df(u(j,1),x)
        +(u(mod(j,n)+1,1)-u(mod(j-2,n)+1,1))/(2*dq) );
\end{reduce}
Compute a norm of the residual, assuming the magnitude of
coefficients in derivative~\(\grad^m\) typically grows
like~\(R^m\) for factor~\(R\approx 2\) or~\(3\), so the norm
is weighted accordingly (via~\verb|d|).   Make the
convergence test relative to the maximum norm found to date
during the iterations.
\begin{reduce}
    rr := part({2,3,2,2,2},mm);
    maxnorm := (if it=1 then 1 else max(maxnorm,normpde));
    tmp:=(pde where{ d=>1/rr, uu(~i,~m,~n)=>uu^i*xx^m*yy^m, gamma=>1 });
    write "Iteration ",it,",  ",
    normpde:=max(abs( for j:=1:n join (
        foreach lu in coeff(tmp(j,1),uu) join (
        foreach lx in coeff(lu      ,xx) join (
                      coeff(lx      ,yy) ))) ));
\end{reduce}

Update the manifold and evolution quasi-adiabatically,
zeroing coefficients smaller than \verb|tolerance|.
\begin{reduce}
    for j:=1:n do rhs0(j,1):=pde(j,1);
    for j:=1:mm do rhs0(n+j,1):=0;
    ugd:=llinv*rhs0;
    for j:=1:n do u(j,1):=u(j,1)+ugd(j,1);
    u:=zeroSmall(u);
    for j:=1:mm do 
        dudt(j-1):=zeroSmall(dudt(j-1)+ugd(n+j,1));
\end{reduce}

Finish the loop when the residuals of the \pde\ are zero to
the specified error tolerance,
\begin{reduce}
    if normpde<relTolerance*maxnorm 
    then write "Success: ",it:=it+100000;
end;
showtime;
if normpde>relTolerance*maxnorm 
then rederr("Iteration Failed to Converge");
\end{reduce}

\subsection{Post-process}

Post-process writes out the multi-continuum, multi-modal,
evolution, 
\begin{reduce}
for i:=0:mm-1 do write "U",i
    ,if alpha=1 then "_t = " else "_tt= " ,dudt(i);
\end{reduce}
Also write out the instructions to draw \LaTeX\ graphs 
of the spatial structure of the invariant manifold, up 
to coefficients of the second order.
\begin{reduce}
vars:={{uu(0,0,0),uu(0,1,0),uu(0,2,0),uu(0,0,2)}
      ,{uu(1,0,0),uu(1,1,0),uu(1,2,0),uu(1,0,2)}};
nd:=3$
array xs(nd); xs(0):=""$ xs(1):="x"$ xs(2):="xx"$ xs(3):="yy"$
if mm=1 then out "hceg2Manifold1.tex"
else if mm=2 then out "hceg2Manifold2.tex"
else out "hceg2Manifold3.tex";
writeTheCase();
write"\begin{tabular}{@{}c@{}c@{}}";
for k:=0:mm-1 do begin
    write if evenp(k+mm-1) then "& " 
    else if k>0 then "\\ " else ""
    ,"\begin{tikzpicture}
    \makeatletter\let\gob\@gobble\makeatother
    \begin{axis}[small
        ,mark repeat=32,mark phase=17
        ,xlabel={microscale $\theta$}
        ,ylabel={component in $u(\theta)$}
        ,legend pos=south west
        ,legend style={font=\footnotesize}]";
    for l:=0:nd do begin
        ukl:=( map(coeffn(~a,part(vars,k+1,l+1),1),u) 
            where d=>1);
        write "\addplot+[] plot coordinates {";
        for j:=1:n do write 
            "    (",q(j),",",ukl(j,1),")\gob"$
        write "};
        \addlegendentry{$U_{",k,xs(l),"}$}"$
    end;
    write "\end{axis}
    \end{tikzpicture}";
end;
write"\end{tabular}";
if mm=1 then shut "hceg2Manifold1.tex"
else if mm=2 then shut "hceg2Manifold2.tex"
else shut "hceg2Manifold3.tex";
\end{reduce}

Finish of script.
\begin{reduce}
end;
\end{reduce}


\section{Computer algebra constructs 2-D elastic homogenisation}
\label{cas2d}


This Reduce-algebra code constructs any chosen
multi-continuum, micromorphic, invariant manifold
homogenisation of the 2-D heterogeneous elasticity example
of~\eqref{EEstress,EEaccel} discussed by \cref{Selastic2d}. 
Specify the number of macroscale multi-continuum
modes~\(M\), and the microscale cell size \(n_x\times n_y\).
\begin{reduce}
mm:=3;
nx:=ny:=10; 
\end{reduce}
Parameter~\verb|ordd| sets the order of the error residual,
and \verb|maxit| the maximum number of iterations.
\begin{reduce}
ordd:=3;   
maxit:=99;
\end{reduce}
Turn on rounded to get numerical coefficients, and only
report two significant digits.
\begin{reduce}
on rounded;  print_precision 2$
\end{reduce}
Coefficients smaller than the following tolerance are set to
zero.  This enables the simple-minded iteration to
terminate, but gives a small numerical error in
coefficients.  
\begin{reduce}
tolerance:=1e-7;
\end{reduce}

Improve printed output:
\begin{reduce}
on div; off allfac; on revpri;
factor d;
\end{reduce}
Load function \verb|zeroSmall| that zeros all numbers
smaller than the set tolerance.  Only works with
\verb|roundbf| off.
\begin{reduce}
in "zeroSmall.txt"$
\end{reduce}

Let the heterogeneity be \(1\)-periodic in both~\(x,y\), and
the same periodicity of say~20 in the sub-cell lattice. 
These lattices are half-spaced in~\(x,y\) and so double
size---provides for staggered micro-grid.
\begin{reduce}
dx:=1/nx;  dy:=1/ny;
matrix xx(2*nx,2*ny),yy(2*nx,2*ny),ones(2*nx,2*ny);
for i:=1:2*nx do for j:=1:2*ny do xx(i,j):=(i-0.5)*dx/2;
for i:=1:2*nx do for j:=1:2*ny do yy(i,j):=(j-0.5)*dy/2;
for i:=1:2*nx do for j:=1:2*ny do ones(i,j):=1;
zeros:=0*ones$
\end{reduce}
Define sets of staggered grid indices: horizontal~\(u\) is
in evens\(\times\)evens, whereas vertical~\(v\) is in
odd\(\times\)odds.
\begin{reduce}
odds:= for i:=1 step 2 until 2*nx collect i;
evns:= for i:=2 step 2 until 2*nx collect i;
\end{reduce}

\paragraph{Define heterogeneity}
Define the heterogeneity in terms of Young's modulus and
Poisson parameter.  Recall \(2\sin A\sin B = \cos(A-B)
-\cos(A+B)\). The following sub-cell heterogeneity, for
\(n_x=10\), leads to a eigenvalue gap ratio~\(1.8\) which is
enough for linear problems.
\begin{reduce}
ee := 1/(2*pi)*(map(cos,(xx-yy)*pi)-map(cos,(xx+yy)*pi))$
write ee := ee+0.01/pi*ones;
write nu:=0.4*ones+0*xx+0*yy;
harmonicMeanE:=(4*nx*ny)
    /(for i:=1:2*nx sum for j:=1:2*ny sum 1/ee(i,j));
\end{reduce}
Compute corresponding Lame parameters.
\begin{reduce}
matrix lla(2*nx,2*ny),lmu(2*nx,2*ny);
for i:=1:2*nx do for j:=1:2*ny do begin
    lla(i,j):=nu(i,j)*ee(i,j)/(1+nu(i,j))/(1-2*nu(i,j)); 
    lmu(i,j):=ee(i,j)/2/(1+nu(i,j)); end;
\end{reduce}
Account for the cell-periodicity of indices, assumes
\(n_x=n_y\)
\begin{reduce}
procedure c(i); mod(i-1,2*nx)+1;
\end{reduce}

\subsection{Compute the cell operator}
Then the cell operator returns matrix of residuals of
\pde{}s on the (odd,odd) and (even,even) lattice points of
the matrix, zero otherwise. Take 1\,s for nx=6 and 4\,s for
nx=8.
\begin{reduce}
matrixproc rescell(uv); begin 
  scalar pde; pde:=0*uv;  
  foreach i in evns do for each j in odds do
    uv(i,j):=lmu(i,j)*( (uv(c(i+1),j)-uv(c(i-1),j))/dx
                       +(uv(i,c(j+1))-uv(i,c(j-1)))/dy );
  foreach i in odds do for each j in evns do
    uv(i,j):=(lla(i,j)+2*lmu(i,j))*(uv(c(i+1),j)-uv(c(i-1),j))/dx
                         +lla(i,j)*(uv(i,c(j+1))-uv(i,c(j-1)))/dy;
  foreach i in evns do for each j in evns do
    pde(i,j):= (uv(c(i+1),j)-uv(c(i-1),j))/dx
              +(uv(i,c(j+1))-uv(i,c(j-1)))/dy ;
  foreach i in odds do for each j in evns do
    uv(i,j):=         lla(i,j)*(uv(c(i+1),j)-uv(c(i-1),j))/dx
        +(lla(i,j)+2*lmu(i,j))*(uv(i,c(j+1))-uv(i,c(j-1)))/dy;
  foreach i in odds do for each j in odds do
    pde(i,j):= (uv(c(i+1),j)-uv(c(i-1),j))/dx
              +(uv(i,c(j+1))-uv(i,c(j-1)))/dy ;  
  return pde;
end;
\end{reduce}

\subsection{Finds eigenmodes via cell Jacobian and SVD}
The displacements~\(u,v\) fit into half of a
\(2n_x\times2n_y\) matrix, so the vector of displacements
would be of length~\(2n_xn_y\), and the Jacobian is
\(2n_xn_y\)-square. For \(n_x=n_y=6\) takes 2.4\,s to
compute the \textsc{svd} for nx=6, 20\,s for nx=8.
\begin{reduce}
write "Forming Jacobian";
matrix jac(2*nx*ny,2*nx*ny);
l:=0$
foreach li in odds do foreach lj in odds do  
for lo:=0:1 do begin
    l:=l+1; 
    uv:=zeros; uv(li+lo,lj+lo):=1;
    lop:=rescell(uv);
    k:=0; 
    foreach ki in odds do foreach kj in odds do 
        for ko:=0:1 do jac(k:=k+1,l):=lop(ki+ko,kj+ko);
end;
jacSym:=max(abs( jac-tp jac ));
if zeroSmall(jacSym)neq 0 then
    rederr("Jacobian is not symmetric");
showtime;
\end{reduce}

Use the \textsc{svd} to get eigenvalues and eigenvectors of
the symmetric Jacobian.
\begin{reduce}
load_package linalg;
Write "Computing SVD of Jacobian (aka eigenvalues)";
usv:=svd(jac)$
uvc:=part(usv,1)$
ss:=part(usv,2)$
showtime;
\end{reduce}
Normalise all eigenvectors to have max-abs value equal one.
\begin{reduce}
for l:=1:2*nx*ny do begin
    tmp := max abs(for k:=1:2*nx*ny collect uvc(k,l));
    for k:=1:2*nx*ny do uvc(k,l):=uvc(k,l)/tmp;
end;
\end{reduce}

Find the smallest \(M+1\) singular values (assuming the
first is not one of the smallest): here, eigenvalues are the
negative of these.
\begin{reduce}
array ls(mm+1);
for m:=1:mm+1 do ls(m):=1;
procedure ssmin(l,m); 
    if ss(l,l)<ss(ls(m),ls(m)) then <<
        for q:=mm+1 step -1 until m+1 do ls(q):=ls(q-1);
        ls(m):=l; >>
    else if m<mm+1 then ssmin(l,m+1);
for l:=2:2*nx*ny do ssmin(l,1);
lss := for m:=1:mm+1 collect ls(m); 
smallest := for m:=1:mm+1 collect ss(ls(m),ls(m)); 
\end{reduce}
Override the two zero singular values and their vectors as
we know these are equivalent to neutral horizontal and
vertical displacements.
\begin{reduce}
l1:=ls(1)$ l2:=ls(2)$
ss(l1,l1):=ss(l2,l2):=0$
for l:=1:nx*ny do begin
    uvc(2*l,l1):=1; uvc(2*l-1,l1):=0;
    uvc(2*l,l2):=0; uvc(2*l-1,l2):=1;
end;
\end{reduce}
Form gravest~\(M\) eigenvectors and eigenvalues into arrays
of spatial matrices.
\begin{reduce}
clear evc;
array evl(mm),evc(mm);
for m:=1:mm do begin
    write evl(m-1):=-ss(ls(m),ls(m));
    ev:=zeros;
    k:=0; 
    foreach ki in odds do foreach kj in odds do 
      for ko:=0:1 do ev(ki+ko,kj+ko):=uvc(k:=k+1,ls(m));
    evc(m-1):=ev;
end;
\end{reduce}

Write approximations to the slow manifold model of the
embedding \pde~\eqref{Eempde} in terms of `modal'
fields~\(U_i(t,x,y)\), denoted by~\verb|uu(i)|, that evolves
according to \(\DD t{U_i}=\verb|dudtt(i)|\) for whatever
\verb|dudtt(i)| happens to be.
\begin{reduce}
array dudtt(mm);
operator uu;  depend uu,x,y,t;
let { df(uu(~i),t,2)=>dudtt(i)
    , df(uu(~i),t,2,x)=>df(dudtt(i),x)
    , df(uu(~i),t,2,y)=>df(dudtt(i),y)
    , df(uu(~i),t,2,x,~p)=>df(dudtt(i),x,p)
    , df(uu(~i),t,2,y,~p)=>df(dudtt(i),y,p)
    };
\end{reduce}

\subsection{Iteration systematically constructs multi-modal model}
\label{seccaisc2d}

Let's iteratively construct the standard, higher-order,
and/or multi-continuum homogenizations~\eqref{EE2dBase},
\eqref{EE2dO4} and~\eqref{EE2dM3}. For \(M=3\) modes and
\(n_x=6\) takes about 1\,s per iteration.
\begin{reduce}
write "
Iteratively Construct
---------------------";
\end{reduce}

Define spatial shift operators of~\cref{EdiffOps} in terms
of derivatives in~\(x,y\). From the operator identity that
\(E=\exp(h\partial_x)\) \cite[e.g.,][p.65]{npl61}, and
using~\verb|d| to count the number of derivatives. Let
\(E_x,E_y\)~denote forward half-shifts, and
\(F_x,F_y\)~denote backward half-shifts.
\begin{reduce}
operator ex,fx,ey,fy;
let { ex(~f)=>f+for n:=1:ordd-1 sum d^n*df(f,x,n)*(+dx/2)^n/factorial(n)
    , fx(~f)=>f+for n:=1:ordd-1 sum d^n*df(f,x,n)*(-dx/2)^n/factorial(n)
    , ey(~f)=>f+for n:=1:ordd-1 sum d^n*df(f,y,n)*(+dy/2)^n/factorial(n)
    , fy(~f)=>f+for n:=1:ordd-1 sum d^n*df(f,y,n)*(-dy/2)^n/factorial(n)
    };
\end{reduce}

Form the extended matrix to use to solve for updates.
\begin{reduce}
write "Finding LU-decomposition for updates";
matrix zerom(mm,mm);
tmp := get_columns(uvc,for m:=1:mm collect ls(m))$
zz := matrix_augment(foreach z in tmp collect z)$ 
jaczz := matrix_augment( jac,zz )$
zz := matrix_augment( (tp zz),zerom )$
jaczz:=matrix_stack(jaczz,zz)$
\end{reduce}
Perform LU decomposition
\begin{reduce}
in_tex "lu_decomp.tex"$
in_tex "lu_backsub.tex"$
lu:=lu_decomp(jaczz)$
showtime;
\end{reduce}
If it got turned on, need to turn off \verb|roundbf| in
order for \verb|zeroSmall| to work its magic.
\begin{reduce}
off roundbf;
\end{reduce}

Start the iteration from the base invariant subspace
approximation for the field and its evolution.
\begin{reduce}
uv:=for m:=0:mm-1 sum evc(m)*uu(m)$
for m:=0:mm-1 do dudtt(m):=evl(m)*uu(m);
\end{reduce}
Seek solution to the specified orders of errors.
\begin{reduce}
for it:=ordd:ordd do let d^it=>0;
for it:=1:maxit do begin write "
**** ITERATION ",it;
\end{reduce}

Zero the stress entries in~\verb|uv| before computing time
derivatives:
\begin{reduce}
foreach i in evns do foreach j in odds do uv(i,j):=0;
foreach i in odds do foreach j in evns do uv(i,j):=0;
\end{reduce}
Compute the residual of the embedding
equations~\cref{EEembedE}, via the stresses: 
\begin{reduce}
  pde:=-df(uv,t,t);  
  foreach i in evns do for each j in odds do
    uv(i,j):=lmu(i,j)*( (ex uv(c(i+1),j) -fx uv(c(i-1),j))/dx
                       +(ey uv(i,c(j+1)) -fy uv(i,c(j-1)))/dy );
  foreach i in odds do for each j in evns do
    uv(i,j):=(lla(i,j)+2*lmu(i,j))*(ex uv(c(i+1),j) -fx uv(c(i-1),j))/dx
                         +lla(i,j)*(ey uv(i,c(j+1)) -fy uv(i,c(j-1)))/dy;
  foreach i in evns do for each j in evns do
    pde(i,j):=pde(i,j) 
        +(ex uv(c(i+1),j) -fx uv(c(i-1),j))/dx
        +(ey uv(i,c(j+1)) -fy uv(i,c(j-1)))/dy ;
  pde:=pde;
  foreach i in odds do for each j in evns do
    uv(i,j):=         lla(i,j)*(ex uv(c(i+1),j) -fx uv(c(i-1),j))/dx
        +(lla(i,j)+2*lmu(i,j))*(ey uv(i,c(j+1)) -fy uv(i,c(j-1)))/dy;
  foreach i in odds do for each j in odds do
    pde(i,j):= pde(i,j)
        +(ex uv(c(i+1),j) -fx uv(c(i-1),j))/dx
        +(ey uv(i,c(j+1)) -fy uv(i,c(j-1)))/dy ;  
\end{reduce}
Trace print the length of the residual.
\begin{reduce}
  pde:=zeroSmall(pde);
  if length(pde(1,1))<10 then write pde11:=pde(1,1);
  write maxlengthpde:=max( map(length,pde) );
\end{reduce}
On the first iteration check the invariant subspace.
\begin{reduce}
if it=1 then if sub(d=0,pde) neq zeros then
    rederr("Invariant subspace not satisfied");
\end{reduce}

To compute update, form residual into column vector, then
apply pre-computed LU factorisation.
\begin{reduce}
matrix res(2*nx*ny+mm,1);
k:=0; 
foreach ki in odds do foreach kj in odds do 
    for ko:=0:1 do res(k:=k+1,1):=pde(ki+ko,kj+ko);
upd:=lu_backsub(lu,res); 
\end{reduce}
Unpack update into corrections
\begin{reduce}
k:=0; 
foreach ki in odds do foreach kj in odds do 
    for ko:=0:1 do uv(ki+ko,kj+ko):=uv(ki+ko,kj+ko)-upd(k:=k+1,1);
for m:=0:mm-1 do dudtt(m):=dudtt(m)+upd(k:=k+1,1);
\end{reduce}

Finish the loop when the residuals of the \pde\ are zero to
the specified error,
\begin{reduce}
  showtime;
  if pde=zeros then write "Success: ",it:=it+100000;
end;
if pde neq zeros then rederr("Iteration Failed to Converge");
\end{reduce}

\subsection{Post-process}

Check amplitudes of at least the mean displacement modes
(more robust would be to check these inside the loop, and
update accordinglyl).
\begin{reduce}
print_precision 3$
meanu:= zeroSmall( (foreach i in evns sum 
                    foreach j in evns sum uv(i,j))/nx/ny );
meanv:= zeroSmall( (foreach i in odds sum 
                    foreach j in odds sum uv(i,j))/nx/ny );
if meanu neq uu(0) then rederr("U0 amplitude not preserved");
if meanv neq uu(1) then rederr("U1 amplitude not preserved");
\end{reduce}

Report effective material constants of the second-order
derivative terms when have chosen two-mode bi-continuum.
\begin{reduce}
if mm=2 then begin
    write mu_h:=coeffn(dudtt(0),df(uu(0),y,2),1);
    write lambda_2:=coeffn(dudtt(0),df(uu(1),x,y),1)-mu_h;
    write lambda_1:=coeffn(dudtt(0),df(uu(0),x,2),1)-2*mu_h;
    write errors:={ mu_h-coeffn(dudtt(1),df(uu(1),x,2),1)
        , lambda_2-coeffn(dudtt(1),df(uu(0),x,y),1)+mu_h
        , lambda_1-coeffn(dudtt(1),df(uu(1),y,2),1)+2*mu_h};
end;
\end{reduce}

Neglect any small coefficients (e.g., less than~0.001), and
print out the evolution \pde{}s, to up to two decimal
places.
\begin{reduce}
tolerance:=0.001;
for m:=0:mm-1 do write dudtt(m):=zeroSmall( dudtt(m) );
uv:=map(zeroSmall,uv)$
\end{reduce}

Print out the terms grouped by order.
\begin{reduce}
array ord(ordd);
for p:=0:ordd-1 do write
    ord(p):=for m:=0:mm-1 collect coeffn(dudtt(m),d,p);
\end{reduce}

Finish of the script.
\begin{reduce}
end; 
rederr("Post script: should not occur");
\end{reduce}


%
%
%
%
%



\section{General solution of modal fractional differential equation}
\label{Sgsmfde}

To fillin details for \cref{SSSfde}, this appendix seeks a general solution~\(u(t)\) of the fractional differential equation (\fde)
\begin{equation}
\alphaD u(t)=\lambda u(t)+q(t),
\label{Efde}
\end{equation} 
for parameter \(\lambda\) that is typically real and negative, and where \(q(t)\) is some given forcing.
In terms of functions to be defined next, we show a general solution is of the form
\begin{equation}
u(t):=\sum_{k=0}^{\alphahat} c_k\mu^{-k}\e_\alpha^{(-k)}(\mu t) 
-\mu^{1-\alpha}\e_\alpha^{(1)}(\mu t)\star q(t),
\quad\text{for }\mu:=(-\lambda)^{1/\alpha}.
\label{Egensol}
\end{equation}
The free constants~\(c_k\) may be identified as the initial values \(c_k=u^{(k)}(0^+)\).
Useful large time asymptotics, from~\eqref{Easymz} and for negative real~\(\lambda\), are 
\begin{align}&
\e_\alpha^{(0)}(t)\sim \frac{t^{-\alpha}}{\Gamma(1-\alpha)}\,,
&&
\e_\alpha^{(-1)}(t)\sim \frac{t^{1-\alpha}}{\Gamma(2-\alpha)}\,,
&&
\e_\alpha^{(1)}(t)\sim \frac{t^{-1-\alpha}}{\Gamma(-\alpha)}\,.
\label{EfdeAsym}
\end{align}
\cref{FfdeSolns1} plots the solution component~\(\e_\alpha^{(0)}(t)\) for various \(0<\alpha<1\) with their large-time approximations.
Further, for various \(1<\alpha<2\), \cref{FfdeSolns2} plots the two solution components~\(\e_\alpha^{(0)}(t)\) and~\(\e_\alpha^{(-1)}(t)\)\,.

The following adapts and generalises some of the report by \cite{Gorenflo2008}, cited by the acronym~\textsc{gm}.

We define \alphaD\ in the Caputo sense \GM{(1.17)}.
First define the fractional integral operator as the following convolution \GM{(1.2)}: for every \(\alpha>0\)
\begin{equation}
\cJ^\alpha f(t) := \frac1{\Gamma(\alpha)}t^{\alpha-1}\star f(t) := \frac1{\Gamma(\alpha)}\int_0^t (t-\tau)^{\alpha-1}\star f(\tau)\d\tau\,.
\label{Efi}
\end{equation} 
Secondly for non-integer~$\alpha$, and defining the integer $\alphahat:=\lceil\alpha\rceil$, the fractional derivative in \fde~\eqref{Efde} is
\begin{equation}
\alphaD u(t):= \cJ^{\alphahat-\alpha} u^{(\alphahat)}(t),
\end{equation} 
Then solutions~\eqref{Egensol} of the \fde~\eqref{Efde} are expressed in terms of the function $\e_\alpha^{(0)}(t):=E_{\alpha,1}(-t^\alpha)$ \GM{(3.11)}, where in turn we define the Mittag--Leffler function 
$E_{\alpha,\beta}(z):=\sum_{k=0}^\infty\frac{z^k}{\Gamma(\alpha k+\beta)}$ \GM{(A.1)}.
Special cases are $E_{0,1}=1/(1-z)$, 
$E_{1/2,1}=\e^{z^2}\erfc(-z)$, 
$E_{1,1}=\e^z$, 
$E_{2,1}=\cosh\sqrt z$.
Various relations and derivatives are useful \GM{(A.8)--(A.10)}, such as
\begin{equation}
\de z{}E_{\alpha,\beta}(z)=\frac1{\alpha z}\left[ 
E_{\alpha,\beta-1}(z) -(\beta-1)E_{\alpha,\beta}(z) \right].
\label{EdEab}
\end{equation}
Two asymptotic expansions may be of interest \GM{(A.22)--(A.23); and \cite{Haubold2011}, (6.10)--(6.11)}:  for $0<\alpha<2$ and as $|z|\to\infty$,
\begin{align}
E_{\alpha,\beta}(z)&
\sim -\sum_{k=1}^\infty \frac{z^{-k}}{\Gamma(\beta-\alpha k)}
+\begin{cases}
\tfrac1\alpha z^{(1-\beta)/\alpha}\exp(z^{1/\alpha})
,& |\arg z|<\alpha\pi/2\,,
\\0\,,&|\arg(-z)|<\alpha\pi/2\,.
\end{cases}
\label{Easymz}
\end{align}
For integer $k>0$, define $\e_\alpha^{(k)}(t)$ to be the $k$th~derivative of~$\e_\alpha^{(0)}(t)$, and $\e_\alpha^{(-k)}(t)$ to be the $k$th~integral from zero of~$\e_\alpha^{(0)}(t)$, that is, $\e_\alpha^{(-k)}(t):=1\star\e_\alpha^{(-k+1)}(t) = \int_0^t \e_\alpha^{(-k+1)}(\tau)\d\tau$\,.
For examples of each, and using~\eqref{EdEab}: 
\(\e_\alpha^{(1)}(t) = \frac1tE_{\alpha,0}(-t^\alpha)\); and
\(\e_\alpha^{(-1)}(t) = tE_{\alpha,2}(-t^\alpha)\) (e.g., differentiate the right-hand side).

\begin{lemma}
The fractional differential equation~\eqref{Efde} has general solution~\eqref{Egensol}.
\end{lemma}

\begin{subequations}\label{EfdePs}%
\begin{proof}
We show that solutions~\eqref{Egensol} satisfy the \fde~\eqref{Efde}.
Take Laplace transform,~\Lap{\cdot}, of \eqref{Egensol}, and since \(\Lap{f(at)}=\tfrac1a\tilde f(s/a)\),
\begin{equation}
\tilde u(s) = 
\sum_{k=0}^{\alphahat} c_k\mu^{-k}\frac1\mu\Lap{\e_\alpha^{(-k)}}(s/\mu) 
-\mu^{1-\alpha}\frac1\mu\Lap{\e_\alpha^{(1)}}(s/\mu) \tilde q(s).
\label{EfdeP1}
\end{equation}
Now \(\tilde\e_\alpha^{(0)}(t) = \frac{s^{\alpha-1}}{s^\alpha+1}\) \GM{(3.11)}.  
Hence by the integration rule, for integer~\(k>0\), \(\Lap{\e_\alpha^{(-k)}(t)}=\frac{s^{\alpha-k-1}}{s^\alpha+1}\), and by the differentiation rule \(\Lap{\e_\alpha^{(1)}(t)} = s\tilde\e_\alpha^{(0)} -\e_\alpha^{(0)}(0^+) = s\frac{s^{\alpha-1}}{s^\alpha+1} -1 = -\frac1{s^\alpha+1}\)  \GM{(3.14)}.
Hence the transform~\eqref{EfdeP1} becomes
\begin{align}
\tilde u(s) &
= \sum_{k=0}^{\alphahat} c_k\mu^{-k-1}\frac{(s/\mu)^{\alpha-k-1}}{(s/\mu)^\alpha+1} 
+\mu^{-\alpha}\frac1{(s/\mu)^\alpha+1} \tilde q(s)
\nonumber\\&
= \sum_{k=0}^{\alphahat} c_k\frac{s^{\alpha-k-1}}{s^\alpha+\mu^\alpha} 
+\frac1{s^\alpha+\mu^\alpha} \tilde q(s)
\nonumber\\&
= \sum_{k=0}^{\alphahat} c_k\frac{s^{\alpha-k-1}}{s^\alpha-\lambda} 
+\frac1{s^\alpha-\lambda} \tilde q(s).
\end{align} 
Multiplying by \(1-\lambda/s^\alpha=(s^\alpha-\lambda)/s^\alpha\) gives
\begin{align}
&(1-\lambda/s^\alpha)\tilde u(s)
=\sum_{k=0}^{\alphahat} c_k\frac{1}{s^{k+1}} 
+\frac1{s^\alpha} \tilde q(s),
\nonumber\\\text{that is,}\quad&
\tilde u(s)
-\sum_{k=0}^{\alphahat} c_k\frac{1}{s^{k+1}} 
=\lambda\frac1{s^\alpha}\tilde u(s)
+\frac1{s^\alpha} \tilde q(s).
\end{align}
Take the inverse Laplace transform:
\begin{align}
u(t) - \sum_{k=0}^{\alphahat} c_k\frac{t^k}{k!} 
=\lambda\frac{t^{\alpha-1}}{\Gamma(\alpha)}\star u(t)
+\frac{t^{\alpha-1}}{\Gamma(\alpha)}\star q(t)
=\lambda\cJ^\alpha u(t)
+\cJ^\alpha q(t).
\label{EfdeP2}
\end{align}
Since \(\alphaD u=\cJ^{-\alpha}\left[ u(t) - \sum_{k=0}^{\alphahat} u^{(k)}(0^+)\frac{t^k}{k!} \right]\) \GM{(1.20)},
identifying \(c_k=u^{(k)}(0^+)\), and applying the operator \(\cJ^{-\alpha}\) to~\eqref{EfdeP2} immediately gives that the~\(u(t)\) of~\eqref{Egensol} satisfies the \fde~\eqref{Efde}, as required.
\end{proof}
\end{subequations}

\begin{figure}\centering
\caption{\label{FfdeTraj}elementary emergence of a slow subspace.  Plot six representative trajectories of the `multiscale' \fde\ system in~\((x(t),y(t))\): \(\alphaD x=-0.03\,x\) and \(\alphaD y=-y\).  The six subplots are for six different \(\alpha\in(0,2)\).  The trajectories show, especially for~\(\alpha\approx 1\), that solutions of the system reasonably quickly reach a neighbourhood of the \(x\)-axis, the slow subspace \(y=0\).}
\includegraphics{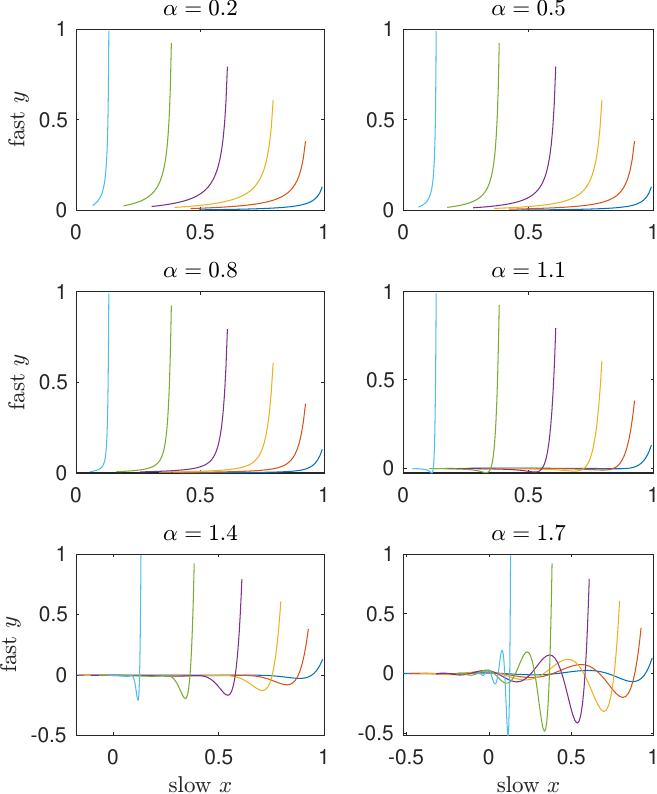}
\end{figure}

\fi

\end{document}